\documentclass[12pt]{amsart}
\usepackage{amssymb}
\usepackage{amscd}
\usepackage{pstricks}
\usepackage{pst-node}
\newtheorem{thm}{Theorem}[section]
\newtheorem{cor}[thm]{Corollary}
\newtheorem{lem}[thm]{Lemma}
\newtheorem{prop}[thm]{Proposition}
\theoremstyle{definition}
\newtheorem{defn}[thm]{Definition}
\newtheorem{example}[thm]{Example}

\newtheorem{notation}[thm]{Notation}

\newtheorem{rem}[thm]{Remark}

\theoremstyle{remark}

\numberwithin{equation}{section}

\newcommand{\norm}[1]{\left\lVert{#1}\right\rVert}
\newcommand{\GL}{\mathit{GL}}

\newcommand{\C}{\mathbb{C}}
\newcommand{\Z}{\mathbb{Z}}

\newcommand{\trans}[1]{{}^t\kern-.2em{#1}}
\newcommand{\ytrans}[1]{{}^t\kern-.11em{#1}}
\newcommand{\Trans}[1]{{}^T\kern-.2em{#1}}
\newcommand{\lsup}[2]{{}^{#1}\kern-.1em{#2}}

\newcommand{\Id}{\operatorname{Id}}
\newcommand{\cal}[1]{\mathcal{#1}}

\newcommand{\Ker}{\operatorname{Ker}}
\newcommand{\Coker}{\operatorname{Coker}}
\newcommand{\sign}{\operatorname{sign}}
\newcommand{\Mas}{\mathbf{Mas}}
\newcommand{\M}{\mathbf{\M}}

\newcommand{\Hi}{\mathbb{H}}
\newcommand{\J}{\mathbb{J}}
\newcommand{\graph}{\operatorname{graph}}
\newcommand{\lla}{\lambda}
\newcommand{\LLa}{\Lambda}
\newcommand{\I}{{\rm Im}}
\newcommand{\Log}{\rm Log}
\newcommand{\Tr}{\rm Tr}
\newcommand{\bg}{{\bf \gamma}}

\renewcommand{\tilde}[1]{\widetilde{#1}}
\renewcommand{\hat}[1]{\widehat{#1}}
\DeclareFixedFont{\bgn}{OT1}{cmr}{m}{n}{20.74}
\DeclareFixedFont{\bgi}{OT1}{cmr}{m}{it}{20.74}

\newcommand{\bigzerou}{\smash{\lower1.7ex\hbox{\bgi O}}}

\newcommand{\inner}[2]{\left\langle{#1},{#2}\right\rangle}

\makeatletter
\def\eqnarray{%
   \stepcounter{equation}%
   \def\@currentlabel{\p@equation\theequation}%
   \global\@eqnswtrue
   \m@th
   \global\@eqcnt\z@
   \tabskip\@centering
   \let\\\@eqncr
   $$\everycr{}\halign to\displaywidth\bgroup
       \hskip\@centering$\displaystyle\tabskip\z@skip{##}$\@eqnsel
      &\global\@eqcnt\@ne \hfil$\displaystyle{{}##{}}$\hfil
      &\global\@eqcnt\tw@ $\displaystyle{##}$\hfil\tabskip\@centering
      &\global\@eqcnt\thr@@ \hb@xt@\z@\bgroup\hss##\egroup
        \tabskip\z@skip
      \cr}
\makeatother
\makeatletter
\def\varin{\mathrel{\mathpalette\@varin\relax}}
\def\@varin#1{%
   \hbox{\setbox\z@\hbox{\m@th$#1\cup$}%
       \def\reserved@a{bold}%
       \dimen@\ifx\reserved@a\math@version .3\else .2\fi\p@
       \kern.5\wd\z@\kern-\dimen@
       \vrule\@width2\dimen@\@height1.08\ht\z@\@depth\z@
       \kern-\dimen@\kern-.5\wd\z@
       \box\z@}}
\makeatother

\begin{document}



\title [Maslov index] 
{Fredholm-Lagrangian-Grassmannian\\
 and\\
the Maslov index}



\author{Kenro Furutani}

\address{Kenro Furutani \endgraf 
Department of Mathematics \endgraf 
Faculty of Science and Technology \endgraf 
Science University of Tokyo \endgraf 
2641 Noda, Chiba (278-8510)\endgraf 
Japan \endgraf}
\email{furutani@ma.noda.sut.ac.jp}

\bigskip
\bigskip

\begin{abstract}
We explain the topology of the space, so called, 
Fredholm-Lagrangian-Grassmannain and the 
quantity ``{\it Maslov index}'' for paths in this space
based on the standard theory of Functional Analysis. 
Our standing point is to define the Maslov index for arbitrary paths
in terms of the fundamental spectral property of the
Fredholm operators, which was first recognized 
by J. Phillips and used to define the ``{\it Spectral flow}''.  
We tried to make the arguments to be all elementary and we summarize
basic facts for this article from Functional Analysis in the Appendix.
\end{abstract}

\maketitle
\tableofcontents


\section*{Introduction}

The purpose of this article is to develop a unified 
theory of the topology of the space, 
``{\it Fredholm-Lagrangian-Grassmannian}'',
and the theory of the Maslov index for arbitrary
paths in this space. Many of the contents of this article 
are treated in the
papers \cite{BF1}, \cite{BF2}, \cite{BFO} and \cite{FO2}.  
Also there are already many papers written which treat with 
more or less similar topics with this paper
(\cite{Le}, \cite{So}, \cite{Du}, \cite{CLM1}, \cite{CLM2},
\cite{Go1}, \cite{RS}, \cite{Sw}, \cite{DK} and others).  
Even so our method to treat with this topics, especially the
treatment of the Maslov index
is different from other papers and so,
the whole theory should be rewritten
in a complete form for being well understood and
will provide a reasonable framework of this subject.   
We would like to emphasis here that 
the method for defining 
the Maslov index for paths with fixed end points
is quite natural and elementary following the basic 
spectral property of the Fredholm operators and valid 
for both in the finite and infinite dimensions.  
We believe that this point must be important,
and should be known widely,
since in the application it naturally appears 
the requirement to treat with
such an integer for not only loops, but also 
paths of Lagrangian subspaces in an intrinsic way.  
Here neither we need any generic arguments
which was assumed in the paper \cite{RS}, nor we  
rely on quantities, ``{\it Leray index}'' and ``{\it Kashiwara index}''
which are only defined for the finite dimensional cases
(\cite{Go1}, \cite{CLM1}) and our method provides simple
and clear theory for dealing with the Maslov index for paths. 

There are many places in which the 
Maslov index and related quantities
appear, and so here we do not mention them, since they are
explained and treated in many articles cited above according to their
subjects. Here I only concentrate
to explain the basic theory of the topology of the
Fredholm-Lagrangian-Grassmannian and the Maslov index for paths
from the point of view of the standard theory of Functional Analysis
and in the elementary ways.

The main method in this article is in the analysis
of operators on Hilbert spaces and the arguments should be
carefully carried out, simply because it is in the infinite dimension. 
There are many parts which are similar to finite dimensional cases, but
also there are many parts which are not just a generalization
of the finite dimensional cases.  We will make clear
the differences of the infinite 
dimensional case from finite dimensional
cases.

We avoid to base on a general theory of the infinite dimensional
manifold theory 
and try
the treatments as elementally as possible and to be self-contained. 
However we must recognize several 
highly non-trivial facts like
\begin{enumerate}
\item
Kuiper's Theorem (\ref{thm:Kuiper} in Appendix)
\item
Palais's Theorem (\ref{thm:Palais} in Appendix)
\item
The spaces of certain class of Fredholm operators are
identified as 
classifying spaces for $K$ and $KO$-groups.
\end{enumerate}

In $\S 1$ we just begin from the very basic facts in symplectic
Hilbert space and the space of their Lagrangian subspaces.
Especially we explain the``{\it Souriau map}'' precisely and
give a proof for determining the fundamental group of 
the Fredholm-Lagrangian-Grassmannian.
$\S 2$ we define the Maslov index for paths, H\"ormander index
and construct the universal covering space of the 
Fredholm-Lagrangian-Grassmannian. Also we discuss the Maslov index
with the relations between certain bilinear forms.
In $\S 3$ we just summarize the finite dimensional cases
and extend the quantity ``Kashiwara index''(``cross index'') to any 
triples of unitary operators.
$\S 4$ we treat with polarized symplectic Hilbert spaces and
prove a symplectic reduction theorem in the infinite
dimensions. Finally in $\S 5$ we discuss examples in this framework
and a formula relating with ``{\it Spectral Flow}''.








\section{Symplectic Hilbert space and Lagrangian subspace}

We start from the very definition of the symplectic Hilbert spaces
and their isotropic, involutive and Lagrangian subspaces and
operations among them.

\subsection{Symplectic Hilbert space}

Let $(H\,,\,<\bullet\,,\,{\bullet}>\,,\,\omega)$ be a 
(real and separable) Hilbert space with an inner product
$<{\bullet}\,,\,{\bullet}>$ and we assume $H$ has a symplectic form
$\omega(\bullet\,,\,\bullet)$, 
i.e., a non-degenerate, skew-symmetric bounded bilinear form.


Here we mean that the bilinear form $\omega$ is non-degenerate
in such a sense that the linear map
\begin{equation}
{\omega}^{\#}:H \rightarrow H^{*} \quad(\text{= dual space})
\end{equation}
\begin{equation*}
{\omega}^{\#}(x)(y)=\omega(x\,,\,y)
\end{equation*}
gives 
an isomorphism between the Hilbert space $H$ and its dual space $H^{*}$.
In finite dimensional cases, the injectivity of the map $\omega^{\#}$
implies that it is an isomorphism, but in the infinite dimension
this does not hold automatically.  In this case we call the Hilbert
space as a symplectic Hilbert space.

In the theory below we do not replace the symplectic form $\omega$
after once it was introduced in the real Hilbert space $H$, but we may
always assume that there exists an orthogonal transformation
$J:\;H\rightarrow H$ such that $\omega(x,y)=<{Jx}\,,\,{y}>$ for any $x,y\in H$ and
$J^2=-\Id$ by replacing the inner product with another one which defines
an equivalent norm on $H$. 

We give the proof of this fact in the appendix.
  
So we can assume from the beginning the following relations:
\begin{equation*}
\begin{array}{l}
\text{$\trans{J}=-J,\quad <{Jx}\,,\,{Jy}>=<{x}\,,\,{y}>$\quad and}\\
\omega(Jx,Jy)=\omega(x,y)\quad\text{for all $x,y\in H$.}
\end{array}
\end{equation*}
Here $\trans{J}$ denotes the transpose of $J$ with respect to the
inner product $<{\bullet}\,,\,{\bullet}>$.  In this case we call
these quantities, the symplectic form $\omega$, the 
inner product $<\bullet\,,\,\bullet>$ and the almost complex structure $J$
are compatible each other.

\begin{example}\label{example:1}
Let $E$ be a real separable Hilbert space and $E^*$ its dual space.
We denote the identification between $E$ and $E^*$, 
by $\mathcal{D} :E \rightarrow E^*$;
$E\ni x \mapsto {\mathcal D}(x)(\bullet) = <\bullet\,,\,x> \in E^*$ (Riesz
representation Theorem).
Then we can introduce an inner product on the dual space $E^*$ through
the map ${\mathcal D}$ in an obvious way and then the direct sum 
$H = E\oplus E^*$ has a naturally defined
skew-symmetric bilinear form 
\begin{equation*}
\omega : H\times H \rightarrow \mathbb{R},
\end{equation*}
\begin{equation*}
\omega(x\oplus\phi\,,\,y\oplus\psi) = \psi(x) - \phi(y) = 
<J(x\oplus\phi)\,,\,y\oplus\psi>,
\end{equation*}
where the almost complex
structure $J : H \rightarrow H$ is given as 
\begin{equation*}
J(x\oplus\phi) = {\mathcal D}^{-1}(\phi)\oplus -{\mathcal D}(x).
\end{equation*}
\end{example}

\begin{example}\label{example:2}
Let $A$ be a densely defined closed symmetric operator on a Hilbert
space $L$.  Let $\mathfrak{D}(A)$ (respectively $\mathfrak{D}(A^*)$) be the 
domain of $A$ (respectively $A^*$) and we impose the graph inner product on
$\mathfrak{D}(A^*)$: $<x\,,\,y>^{\cal{G}} = <x\,,\,y>
+<A^*(x)\,,\,A^*(y)>$.  
Then $\mathfrak{D}(A^*)$
becomes a Hilbert space and $\mathfrak{D}(A)$ is a closed subspace
in $\mathfrak{D}(A^*)$ with respect to this graph norm.  Let
${\bf{\beta}}$
be the factor space ${\bf{\beta}}$ = $\mathfrak{D}(A^*)/\mathfrak{D}(A)$.  We can introduce a
non-degenerate anti-symmetric bilinear form $\omega$ on ${\bf{\beta}}$ by
\begin{equation}
\omega([x]\,,\,[y]) = <A^*(x)\,,\,y> - <x\,,\,A^*(y)>,
\end{equation}
where we denote by $[x]$, the class of $x\in \mathfrak{D}(A^*)$ in ${\bf{\beta}}$.

It will be apparent of the well-definedness of the form $\omega$ just
from the definition of the adjoint operator.

We will note the non-degeneracy of the form $\omega$: The factor space
${\bf{\beta}}$ is identified with the orthogonal complement $\mathfrak{D}(A)^{\perp}$ 
of $\mathfrak{D}(A)$ in $\mathfrak{D}(A^*)$ with respect
to the graph inner product.  
It 
is characterized as follows.
\begin{equation*}
\mathfrak{D}(A)^{\perp}=\{x\in \mathfrak{D}(A^*)\,|\,A^*(x) \in \mathfrak{D}(A^*)\,
~\text{and}~A^*(A^*(x)) = -x\}.
\end{equation*}
{}From this characterization we know at once that
$A^*$ restricted to $\mathfrak{D}(A)^{\perp}$ is an orthogonal
transformation
into itself and defines an almost complex structure on
$\mathfrak{D}(A)^{\perp}$
and moreover we have
\begin{equation*}
\omega([x]\,,\,[y]) = <A^*(x)\,,\,y> - <x\,,\,A^*(y)>
\end{equation*}
\begin{equation*}
=<A^*(x)\,,\,y> +<A^*(A^*(x))\,,\,A^*(y)>
=<A^*(x)\,,\,y>^{\cal{G}}.
\end{equation*}
This equality shows that our Hilbert space ${\bf{\beta}}$ with
the symplectic form $\omega$ above together with the almost complex structure
$A^*$ (after being identified with the orthogonal complement
$\mathfrak{D}(A)^{\perp}$)
is  a symplectic Hilbert space with a compatible symplectic form,
inner product and the almost complex structure.
\end{example}

We will deal with this example in $\S 5$ 
together with a homotopy
invariant, so called, 
``{\it {Spectral flow}}'' 
of a family of selfadjoint Fredholm operators.

\begin{example}\label{example:3}
Let $\pi : E \rightarrow M$ be a real vector bundle on a manifold
$M$ with a bundle map of 
almost complex structures $J :E\rightarrow E\,, J^2 =-Id$.  By
introducing a suitable inner product on $E$ and a (smooth) measure
on $M$ we have a Hilbert space $L_2(M\,,\,E)$ of $L_2$-sections 
of $E$ with a symplectic form defined by the bundle map $J$ in an
obvious way. 
\end{example} 

When we regard the real Hilbert space $H$ as a complex Hilbert
space through the almost complex structure $J$ with the Hermitian
inner product 
$<{\bullet}\,,\,{\bullet}>_{J}
=<{\bullet}\,,\,{\bullet}>-\sqrt{-1}\omega(\bullet,\bullet)$, we denote it by
$H_J$,
and we denote the group of unitary transformations on $H_J$ by 
\begin{equation*}
\cal{U}(H_J) 
=\{U\in\cal{B}(H)\;|\; UJ=J\,U \;\text{and}\; {}^tUU =U\,{}^tU=\Id \},
\end{equation*}
where $\cal{B}(H)$ denotes the space of bounded linear operators on the real
Hilbert space $H$.

For a subspace $\mu$ in $H$, let us denote by $\mu^{\circ}$ the
annihilator
of $\mu$ with respect to $\omega$:
\begin{equation}
\mu^{\circ} = \{x\in H\,|\, \omega(x, y)= 0 \,\,\text{for all}\,\,y
\in \mu \},
\end{equation}
and we denote the orthogonal complement (with respect to the fixed inner
product $<\bullet\,,\,\bullet >$ on $H$) of $\mu$ by $\mu^{\perp}$.

Note that we know easily by the definition that for any subspace $\mu$
the annihilator  $\mu^{\circ}$ is closed by the similar way to prove 
the closedness of the orthogonal compliment $\mu^{\perp}$. 
Also by the non-degeneracy assumption of the
symplectic form, we have the 
idempotentness of the operation $\mu \mapsto \mu^{\circ}$: 

\begin{prop}
$(\mu^{\circ})^{\circ} =\overline{\mu}$.
\end{prop}

\begin{proof}
By the definition of the annihilator it will be apparent that $\overline{\mu}
\subset (\mu^{\circ})^{\circ}$. 
Let $z_0 \in (\mu^{\circ})^{\circ}$ and assume that $z_0 \not\in\overline{\mu}$, then
there is a bounded linear functional $f$ on $H$ such that
$f= 0$ on $\mu$ and $f(z_0)\ne 0$.  By the non-degeneracy
assumption of the symplectic form $\omega$, we have an element $u_0 \in
H$ such that
$f(x) =\omega(x\,,\,u_0)$.  Then $u_0 \in \mu^{\circ}$, 
but $\omega(z_0\,,\,u_0)\ne 0$. This is a contradiction.
So there are no such $z_0$. 
\end{proof}

The following properties will be proved easily.
\begin{prop}
Let $\mu$, $\nu$ be subspaces in $H$, then
\begin{align}
(\mu +\nu)^{\circ} &= \mu^{\circ}
    \bigcap\nu^{\circ}\\
(\mu\bigcap\nu)^{\circ}& = {\overline{\mu^{\circ}+\nu^{\circ}}}
\end{align}
\end{prop}

As in the same way with finite dimensional cases we characterize
a subspace $\mu\in H$ in the  
\begin{defn}\label{defn:four subspaces}
\begin{enumerate}
\item
isotropic, if $\mu\subset \mu^{\circ}$
\item
Lagrangian, if $\mu^{\circ}= \mu$
\item
coisotropic (or involutive), if $\mu^{\circ} \subset \mu$
\item
symplectic, if $\mu$ is closed and $\mu +\mu^{\circ} = H$ (= direct sum).
\end{enumerate}
\end{defn}

By the compatibility assumption among the symplectic form 
$\omega$, the inner
product $<\bullet\,,\,\bullet>$ and
the almost complex structure $J$ the following properties hold:
\begin{prop}\label{p:isolaginvo}
\begin{enumerate}
\item
if $\mu$ is isotropic, then $J(\mu)$ is also isotropic and
$\mu\perp J(\mu)$.
\item
if $\mu$ is Lagrangian, then $\mu$ is a closed subspace,
$J(\mu)$ is also Lagrangian and
$J(\mu) = \mu^{\perp}$.  Conversely let $\mu$ be a closed subspace
and assume that $\mu^{\perp} =J(\mu)$, then $\mu$ is a 
Lagrangian subspace.
\item
if $\mu$ is coisotropic, the $J(\mu)$ is also coisotropic.
\end{enumerate}
\end{prop}

If $\mu$ is symplectic, then $\mu +\mu^{\circ}$ is a direct sum,
however it is not always orthogonal.  Then
\begin{prop}\label{prop:symplectic subspace}
If $\mu$ is symplectic, then the restriction of the map 
$\omega^{\#}$ to each  of $\mu$ and 
$\mu^{\circ}$ is isomorphic with $\mu^*$ and $(\mu^{\circ})^*$ 
respectively.
So, by replacing the inner product with a suitable one so that 
we can assume that  $\mu$ and $\mu^{\circ}$ are orthogonal and then
each is a symplectic Hilbert space with the compatible structure.
\end{prop}

\begin{proof}
If we embed
$\mu^*$ into $H^*$ by extending $f\in\mu^*$ to $\tilde{f}$ 
being zero on $\mu^{\circ}$, then
for any $f$ there is an element $a+b \in \mu+\mu^{\circ}$ such that
$\omega^{\#}(a+b) = \tilde{f}$ and from the assumption, $b$ must be zero, that
is, we have $(\omega_{|\mu\times\mu})^{\#} =
(\omega^{\#})_{|\mu}$. 
Hence
$\mu$ is a symplectic Hilbert space.  So is $\mu^{\circ}$.  Then
the rests of the proposition will follow easily {}from Proposition
\ref{prop:appen1} in the Appendix.
\end{proof}

\begin{rem}
\begin{enumerate}
\item
Let $E$ be a finite dimensional subspace in $H$ 
such that $E\bigcap E^{\circ}=\{0\}$,
then $E$ is symplectic in the above sense of 
Definition \ref{defn:four subspaces} (d), that is $E\oplus E^{\circ} =H$
\item
Let $\lla$ be a Lagrangian subspace and $L$ be a closed subspace in
$\lla$. Put $H_1 =L+J(L)$ and 
$H_2= L^{\perp}\bigcap\lla+J(L^{\perp}\bigcap\lla)$, then $H_1$  and $H_2$
are symplectic, of course with the compatibility assumption of the
symplectic structure on $H$.
\end{enumerate}
\end{rem}

\subsection{Lagrangian-Grassmannian}

Let $\Lambda (H)$ denote the space of all Lagrangian subspaces of $H$.
We call this space {\it Lagrangian Grassmannian} of the symplectic Hilbert
space $H$.

\begin{rem}
\begin{enumerate}
\item
Let $\lambda \in \Lambda (H)$, then by the above Proposition \ref{p:isolaginvo}
we have an orthogonal decomposition $H = \lambda\oplus J(\lambda)$ and
by identifying the dual space of $\lambda$ with $J(\lambda)$ we know
that any
symplectic Hilbert space is isomorphic with the Example \ref{example:1}.
\item
In the symplectic Hilbert space a maximum isotropic subspace
is always a Lagrangian subspace.  For the symplectic Banach space
(this is defined by the same way as symplectic Hilbert spaces)
a maximal isotropic subspace need not be a Lagrangian subspace,
moreover there is a symplectic Banach space which has no Lagrangian
subspace (see \cite{KS}).  This fact says that a symplectic Banach
space is not necessarily isomorphic with a standard one of the form
$V\oplus V^*$ with a reflexive Banach space $V$.
In this article we do not treat with the symplectic Banach space.
\end{enumerate}
\end{rem}

We denote
by $\cal{P}_{\lambda}$ the orthogonal projection operator in $H$ onto
the subspace $\lambda$. With this correspondence
we embed $\Lambda(H)$ into $\cal{B}(H)$ as
a closed subset (see Corollary \ref{cor:closedness} below for the closedness): 

\begin{equation}
\begin{array}{r@{\,}l}
\cal{P}:\Lambda (H) & \rightarrow \cal{B}(H)\\
\lambda &\mapsto \cal{P}_\lambda
\end{array}
\end{equation}

{\it Then it will be natural to introduce the metric $d$
on the space $\Lambda(H)$ as the difference of
the norm of the corresponding projection operators :
$d(\lambda\,,\,\mu) = \norm{\cal{P}_{\lambda}-\cal{P}_{\mu}}$.
Henceforth 
we regard the space $\Lambda(H)$ equipped with this metric always.}

A projection operator in the image of the map $\cal{P}$ is
characterized by the following 

\begin{prop}\label{orthogonal-projection-and-Lag}
Let $P$ be an orthogonal projection operator in $H$.  Then the image
$P(H)$ is a Lagrangian subspace, if and only if $J= J\circ P + P\circ
J$.
\end{prop}
\begin{proof}
Let an orthogonal projection operator $P$ satisfy the relation $J=
JP+PJ$, then we have $\omega(P(x),~P(y)) = <J\circ P(x),~P(y)> =
<J(x)-P\circ J(x),~P(y)> = <J(x),~P(y)> -<J(x),~P(y)>=0$.  So $P(H)$
is
an isotropic subspace.  
Let assume for any $x \in H$ $\omega(P(x),~y) =0$, then
$<J(x)-P\circ J(x),~y>=0 $. So we have $<J(x),~y-P(y)> =0$ for any
$x\in$ $H$.  Hence $y=P(y)$, and so $P(H)^{\circ} = P(H)$, that is
$P(H)$ is a Lagrangian subspace.

Now assume that $P(H)$ is a Lagrangian subspace.  Then we have for
$x\in P(H)$, $J(x) = P\circ J(x) +J\circ P(x)$, since $P\circ J(x)=0$ 
and for $x \in \Ker(P)$ we have $J(x) = P\circ J(x) +J\circ P(x)$, since 
$P\circ J(x)= J(x)$.
\end{proof}

As a Corollary of this proposition
we have

\begin{cor}\label{cor:closedness}
The subspace consisting of orthogonal projections whose image is a
Lagrangian subspace is closed in the Banach space $\cal{B}(H)$.
\end{cor}

The group $\cal{U}(H_J)$ acts on $\Lambda(H)$ in an obvious way. Then

\begin{prop}
The action $\cal{U}(H_J)\times\Lambda(H) \to \Lambda(H)$ is continuous.
\end{prop}
\begin{proof}
{}From the relation 
\begin{equation}
\cal{P}_{U(\mu)}= U\circ \cal{P}_{\mu}\circ U^{-1},
\end{equation} 
we have
\begin{align*}
&\cal{P}_{U(\mu)}-\cal{P}_{V(\nu)}\\
&=U\circ \cal{P}_{\mu}\circ U^{-1} - V\circ \cal{P}_{\mu}\circ V^{-1}\\
&=U\circ(\cal{P}_{\mu}-\cal{P}_{\nu})\circ U^{-1}
+(U-V)\circ\cal{P}_{\nu}\circ U^{-1}+V\circ\cal{P}_{\nu}\circ(U^{-1}-V^{-1}),
\end{align*} 
and this formula shows the continuity of the action.
\end{proof}

By fixing an $\ell \in  \Lambda(H)$ we have a surjective map $\rho_{\ell}$:
\begin{equation}\label{eq:fiber-map}
\begin{array}{r@{\,}l}
\rho_{\ell}:\cal{U}(H_J) & \rightarrow \Lambda(H)\\
U &\mapsto U({\ell}^{\perp}).
\end{array}
\end{equation}

\begin{thm}\label{thm:trivial}
The map (\ref{eq:fiber-map})
defines a principal fiber bundle with the structure group
$\cal{O}(\ell)$, the group of orthogonal 
transformations on $\ell$, and
by Kuiper's Theorem \ref{thm:Kuiper} it is a trivial bundle and
the space $\Lambda(H)$ itself is also contractible.
\end{thm}

\begin{rem}
Of course the triviality of this bundle
is not true for the finite dimensional case.
\end{rem}

\begin{cor}\label{cor:quotient topology}
The map $\rho_{\ell}$ is an open map and the topology
on $\Lambda(H)$ coincides with the quotient topology of ${\cal U}(H_J)$
by the map $\rho_{\ell}$.
\end{cor}

Theorem \ref{thm:trivial} is proved 
if we have local sections of the map
$\rho_{\ell}$.  Here
we construct local sections in two ways. Because both of the arguments contain
several interesting properties of the space $\Lambda(H)$ and
relating properties of projection operators.

\smallskip 

{\bf I}. {\it{First method}.}\\
We begin from a lemma:
\begin{lem}  Let $P$ and $Q$ be two projection operators on the
  Hilbert space $H$, and assume that $\norm {P- Q} < 1$.
Put 
\begin{enumerate}
\item $A = (1-P)(1-Q)+PQ$, 
\item $B=(1-Q)(1-P)+QP$, 
\item $C= 1-(P-Q)^2$ 
\item $D= \sum\limits_{n=0}^{\infty}\alpha_n (P-Q)^{2n}$, where
$(1-x)^{-1/2} =\sum\limits_{n=0}^{\infty}\alpha_n x^n$ is the Taylor
expansion.  
\end{enumerate}
Then we have
\begin{enumerate}
\item  $A B = B A = C$,
\item $D^2 C = C D^2 = I$,
\item  $P(P-Q)^2 =(P-Q)^2 P$, $Q (P=Q)^2 = (P-Q)^2Q$ and
\item $DP=PD$, $DQ=QD$.
\end{enumerate}
\end{lem}
\begin{proof}
All these will be proved by direct calculations. Note that
all of the operators $A$, $B$, $C$ and $D$ are, as a result, 
invertible and $AC =CA\,,\,CB=BC$ and $DC=CD$. 
\end{proof}

Now put $W = DA$, then
\begin{prop}
\begin{enumerate}
\item $W$ is invertible and the inverse is given by $W^{-1}= BD$,
\item WQ =PW.
\end{enumerate}
Hence we have $W(Q(H))= P(W(H))$.  Moreover if both of $P$ and $Q$ are
orthogonal
projections, the  the operator $W$ is unitary, that is the ranges of
the projections $P$ and $Q$ are transformed each other by a unitary
operator $W$.
\end{prop}
\begin{proof}
$WQ = D((1-P)(1-Q)+PQ)Q = DPQ$ and
$PW = PDA = DPA = DP((1-P)(1-Q)+PQ)= DPQ$. 
Since $(DA)(BD) = DCD =1$ and $(BD)(DA) = BD^2A = BC^{-1}A = 1$
the operator $W$ is invertible.  Also since $W^* =A^*D^*$, 
if both of $P$ and $Q$ are orthogonal
we have $W^* = BD=W^{-1}$, that is, $W$ is unitary and give a unitary
equivalence of the projections $P$ and $Q$.
\end{proof}

Let $\mu \in \Lambda(H)$ and ${\bf{V}}_{\mu}$
$= \{ \nu\, |\, \norm {\cal{P}_{\nu}- \cal{P}_{\mu}} \,<\, 1\}$,
an open neighborhood of $\mu$,
where $\cal{P}_{\nu}$ denote the orthogonal projection operator with the
image $\nu$.  

Now we give a local section 
$s_{\mu}^{(1)} :{\bf{V}}_{\mu} \rightarrow \cal{U}(H_J)$
of the map 
$\rho_{\ell}:\cal{U}(H_J)  \rightarrow \Lambda(H)$.  

We fix a unitary operator $V_0$ such that
$V_0(\ell^{\perp})=\mu$ and define
\begin{equation}\label{local-section-1}
s_{\mu}^{(1)} : {\bf{V}}_{\mu} \ni \nu \rightarrow W^{-1}_{\nu}\circ V_0,
\end{equation}
where we denote $W_{\nu} =
(1-(\cal{P}_{\mu}-\cal{P}_{\nu})^2)^{-1/2}
((1-\cal{P}_{\mu})(1-\cal{P}_{\nu})+\cal{P}_{\mu}\cal{P}_{\nu})$.
The continuity of this section will be apparent from the expression.

\smallskip

{\bf II}. {\it{Second method}.}
\smallskip

Let $\lambda \in \Lambda(H)$.

\begin{notation}\label{not:transversal}
${\bf O}_{{\lambda}}$
= $\{ \mu\in\Lambda(H)\,|\,\mu \,\,
\text{is transversal to}\, \lambda \}$.
Note that we mean {\it transversal} by the condition: $\lla + \mu = H$. 
\end{notation}

The subset
${\bf O}_{{\lambda}^{\perp}}$ is an open neighborhood of $\lambda$.
We denote by $\widehat{\cal{B}}(\lambda)$ the space of selfadjoint bounded
operators on the real Hilbert space $\lambda$.  Then we have a
bijection 
\begin{equation*}
G_{\lambda} :\widehat{\cal{B}}(\lambda) \rightarrow {\bf O}_{{\lambda}^{\perp}}
\end{equation*}
defined by
\begin{equation*}
G_{\lambda}:\widehat{\cal{B}}(\lambda)\ni A \rightarrow
G_{\lambda}(A)=\{x+JA(x)\,|\,x\in\lambda \}\in {\bf O}_{{\lambda}^{\perp}}.
\end{equation*}

By the identification $H_J \cong\lambda\otimes\mathbb{C}$
we regard $A\in\widehat{\cal{B}}(\lambda)$ as a selfadjoint operator on
$H_J$.
Let $A =\int\limits_{-\infty}^{\infty} \,t\, dE_t(A)$ be the
spectral decomposition of the selfadjoint operator $A$ 
with the spectral measure $\{E_t(A)\}_{t\in\mathbb{R}}$. We define
a unitary operator $U_A$ by
\[U_A 
=\int \sqrt{\frac{1+\sqrt{-1}t}{1-\sqrt{-1}t}}\,
dE_t(A),
\]
then 
\begin{equation}
\int(1+\sqrt{-1}t)d E_t(A)
= U_A \circ \int\sqrt{1+ t^2}d E_t(A).
\end{equation}
Since $\int\sqrt{1+t^2}dE_t(A)(\lambda) 
= (\Id+A^2)^{1/2}(\lambda)= \lambda$,
we have
\begin{equation}
U_A(\lambda)= \int(1+\sqrt{-1}t)dE_t(A)(\lambda) = G_{\lla}(A).
\end{equation}

Note that $U_A^{~2} = (\sqrt{-1}\Id-A)(\sqrt{-1}\Id+A)^{-1}$ is the Cayley
transformation of the operator $A$.

Now fix a unitary operator $V$ 
such that $V(\ell^{\perp})=\lla^{\perp}$, then
the correspondence 
\begin{equation}\label{local-section-2}
s_{\lla^{\perp}}^{(2)} : {\bf O}_{{\lambda}^{\perp}} 
\ni \mu \rightarrow U_A\circ J\circ V
\end{equation}
gives a local section of the map
\begin{equation*}
\begin{array}{r@{\,}l}
\rho_{\ell}:\cal{U}(H_J) & \rightarrow \Lambda(H)\\
U &\mapsto U(\ell^{\perp}).
\end{array}
\end{equation*}

We must show the continuity of this section $s_{\lla^{\perp}}^{(2)}$.  
This is proved by showing two facts: 
\begin{enumerate}
\item 
the continuity of the correspondence 
\[\widehat{\cal{B}}(\lambda)\ni A \mapsto U_A 
= \int\sqrt{\frac{1+\sqrt{-1}t}{1-\sqrt{-1}t}}dE_t(A) \in \cal{U}(H_J)
\]
with respect to the norm topology and, 
\item
the map $G_{\lambda}$ 
is an isomorphism between the spaces $\widehat{\cal{B}}(\lambda)$ and
${\bf O}_{\lambda^{\perp}}$.
\end{enumerate}

The first one follows from a more general
\begin{prop}[\cite{AS}]
Let $H$ be a Hilbert space (real or complex) and $f$ be a continuous
function defined on $\mathbb R$, then the map 
$\widehat{\cal{B}}(H)\ni A \rightarrow f(A) \in \cal{B}(H)$
is continuous.  
Here the operator $f(A) =\int f(t)dE_t(A)$ is defined by the spectral 
decomposition $A=\int t\, dE_t(A)$ of the operator $A$.
\end{prop}

\begin{proof}
Let $\{p_n(t)\}$ ($n=1,2,\dots$) be polynomials which converge uniformly
to the continuous function $f$ on an interval $[-N\,,\,N]$, then
for the operator $A$ whose spectrum $\sigma(A)$ is contained in the
open interval $(-N\,,\,N)$,\\
the operator $ p_n(A)=\sum_{k\geq 0}^{N_{n}}c_{k}^{n}A^{k}$ 
is also expressed as
  \[
 p_n(A) = \int\limits_{-N}^{+N} p_n(t) dE_t(A).
\]
So we know that $\{p_n(A)\}$ converges to $\int f(t) dE_t(A)$ in the sense of operator norm.  The correspondence
$A \mapsto p_n(A)$ is apparently continuous on the open subspace 
$\{A\in\widehat{\cal{B}}(H)\,|\, \sigma(A)\subset (-N\,,\,N)\}$
in $\widehat{\cal{B}}(H)$ and so the map
$\widehat{\cal{B}}(\lambda)\ni A\mapsto f(A)\in \cal{B}(\lambda)$ is continuous
on each such open subspace 
$\{A\in\widehat{\cal{B}}(H)\,|\, \sigma(A)\subset(-N\,,\,N)\}$.
Hence we have the desired result. 
\end{proof}

\begin{prop}\label{prop:local chart}
The map $G_{\lambda}:A \mapsto G_{\lambda}(A) = \{x+JA(x)\,|\, x\in \lambda\}$
is an isomorphism between the spaces $\widehat{\cal{B}}(\lambda)$ and
${\bf O}_{\lambda^{\perp}}$.  Hence it gives a local chart of $\Lambda(H)$
\end{prop}

We prove a characterization of an orthogonal projection operator 
corresponding to a Lagrangian subspace
in ${\bf O}_{\lambda^{\perp}}$. 
\begin{lem}\label{lem:transversal}
Let $\cal{P}_{\mu}$ be an orthogonal projection operator onto a
Lagrangian subspace $\mu\in \Lambda(H)$.  Then 
$\mu \in {\bf O}_{\lambda^{\perp}}$, if and only if
 $L_{\mu}=\cal{P}_{\mu} + 1-\cal{P}_{\lambda} = \cal{P}_{\mu}
 +\cal{P}_{\lambda^{\perp}}$  is an isomorphism.  
\end{lem}
\begin{proof}
If $L_{\mu}=\cal{P}_{\mu}+\cal{P}_{\lambda^{\perp}}$ is an isomorphism, then
since  $H =(\cal{P}_{\mu}+\cal{P}_{\lambda^{\perp}})(H)\subset \mu+\lambda^{\perp}$
we know at once $\mu\in {\bf O}_{\lambda^{\perp}}$. 

Conversely let us
assume  $\mu$ and $\lambda^{\perp}$ are transversal.  Then there is a
bounded operator $A\in\widehat{\cal{B}}(\lambda)$ such that 
$\mu=G_{\lambda}(A) 
=\{x+JA(x)\,|\, x\in\lambda\}$, the graph of the operator $A$.  Note
that
the boundedness of the operator $A$ is proved by 
the closed graph Theorem 
and the selfadjointness of $A$ comes from the fact that $\mu$ is 
a Lagrangian subspace.
These arguments are same with that of finite dimensional cases.
Now we solve the equation
\begin{equation}\label{eq:surjective}
L_{\mu}(u+J(v))
=(\cal{P}_{\mu}+\cal{P}_{\lambda^{\perp}})(u+J(v)) =x+J(y),
\end{equation}
for any given $x$, $y\in \lambda$ by $u$ and $v,\in \lambda$.  
Since $\cal{P}_{\mu}(u+J(v))= x+JA(x)$,
$u+J(v) = x+JA(x) +J(a+JA(a))\in \mu+\mu^{\perp}$ with an element
$a\in\lambda$, 
and $\cal{P}_{\lambda^{\perp}}(u+J(v)) = J(v)=J(y-A(x))$ we have
\begin{align*}\label{eq:solution}
a&=x\\
u&=x-A(a)=x-A(y)+2A^2(x)\\
v&=y-A(x).
\end{align*}
This implies that the operator 
$L_{\mu}=\cal{P}_{\mu} +\cal{P}_{\lambda^{\perp}}$
is an isomorphism of $H$.  Note that we have in general
$\Ker(\cal{P}_{\mu}+\cal{P}_{\nu})
=\Ker(\cal{P}_{\mu})\bigcap \Ker(\cal{P}_{\nu})=
\mu^{\perp}\bigcap\nu^{\perp}=
J(\mu\bigcap\nu)$
(see the proof of Proposition \ref{prop:Fred-pair}).
\end{proof}

\begin{rem}
In Proposition \ref{prop:Fred-pair} we will give a generalization 
of this property after introducing the notion of "\textit{Fredholm pair}".
\end{rem}

\begin{proof}[Proof of Proposition \ref{prop:local chart}]

Let $\mu$ and $\nu$ be transversal with $\lambda^{\perp}$, then
\begin{equation}\label{eq:norm difference}
\norm{L_{\mu}-L_{\nu}} =\norm{\cal{P}_{\mu}-\cal{P}_{\nu}}. 
\end{equation}
So we have 
\[ L_{\nu}\,^{-1} 
= \sum\limits_{k=0}^{\infty}
(L_{\mu}\,^{-1}(L_{\mu}-L_{\nu}))^k\cdot (L_{\mu}\,^{-1})
\]
for such $\nu$ that $\norm{L_{\mu}\,^{-1}}\norm{\cal{P}_{\mu}-\cal{P}_{\nu}}<1$,
and we have
\[\norm{L_{\nu}\,^{-1}} \leq 
\sum\limits_{k=0}^{\infty}
\norm{L_{\mu}\,^{-1}}^{k+1}(\norm{\cal{P}_{\mu}-\cal{P}_{\nu}})^k
= \norm{L_{\mu}\,^{-1}}
\frac{1}{1-\norm{L_{\mu}\,^{-1}}\norm{\cal{P}_{\mu}-\cal{P}_{\nu}}}.
\]
Hence we have 
\[
\norm{L_{\mu}\,^{-1} - L_{\nu}\,^{-1}}
\leq 
\norm{L_{\mu}\,^{-2}}
\frac{1}{1-\norm{L_{\mu}\,^{-1}}\norm{\cal{P}_{\mu}-\cal{P}_{\nu}}}
\norm{\cal{P}_{\mu}-\cal{P}_{\nu}}.
\]

Now by putting $x=0$ in the equation (\ref{eq:surjective})
we have 
\begin{equation}\label{eq:inverse}
L_{\mu}\,^{-1}(J(y)) = -A(y) +J(y)
\end{equation} 
and  
we have the inequality
\begin{align*}
&\norm{A_{\mu}(y)-A_{\nu}(y)}\\
&\leq \norm{L_{\mu}\,^{-1}(J(y))- L_{\nu}\,^{-1}(J(y))}\\
&\leq
\frac{\norm{L_{\mu}\,^{-2}}}{1-\norm{L_{\mu}\,^{-1}}\norm{\cal{P}_{\mu}-\cal{P}_{\nu}}}
\norm{\cal{P}_{\mu}-\cal{P}_{\nu}}\norm{y},
\end{align*}

The last inequality shows that the map 
$G_{\lambda}^{-1}:{\bf O}_{\lambda^{\perp}}
\rightarrow \widehat{\cal{B}}(\lambda)$
is continuous.

The continuity of the map 
$G_{\lambda}:\widehat{\cal{B}}(\lambda) \rightarrow {\bf O}_{\lambda^{\perp}}$
is proved more easily:
let $\mu\in {\bf O}_{\lambda^{\perp}}$, that is, 
$\mu$ and $\lambda^{\perp}$ are transversal. Then we can express in two
ways of the element $x+J(y) \in H$:
\[
\lambda+\lambda^{\perp} \ni x+J(y) = a+JA(a)+J(b+JA(b))\in
\mu+\mu^{\perp}.
\]

By solving this equation we have
\begin{align*}
a=&(\Id+A^2)^{-1}(x+A(y))\\
b=&(\Id+A^2)^{-1}(y-A(x)),
\end{align*}
so 
\begin{align}\label{eq:projection}
&\cal{P}_{G_{\lambda}(A)}(x+J(y))\\
&= (\Id+A^2)^{-1}(x+A(y))
+JA((\Id+A^2)^{-1}(x+A(y)))\notag.
\end{align}

{}From this expression of the projection $\cal{P}_{G_{\lambda}(A)}$
and by a standard argument we have 
\begin{equation*}
\norm{\cal{P}_{G_{\lambda}(A)}- \cal{P}_{G_{\lambda}(B)}}
\leq N(\norm{A},\norm{B})\norm{A-B}.
\end{equation*}
Here we denote by $N(s,t)$ a polynomial of degree three of two variables
and note that for any $A\in\widehat{\cal{B}}(\lambda)$
$\norm{(\Id+A^2)^{-1}}\leq 1$.

Consequently we have proved both of the continuities of the map $G_{\lambda}$
and its inverse $G_{\lla}^{~-1}$, in other words,
we have proved that the map $G_{\lambda}$ 
gives a local chart of the space $\Lambda(H)$.
\end{proof}
\begin{rem}
{}From the proof above we see easily that the map $G_{\lambda}$
is not isometric.
\end{rem}
\smallskip

\begin{proof}[Proof of Theorem \ref{thm:trivial}]
It will be clear that a unitary operator  
which preserves the subspace
$\ell^{\perp}$ comes from an orthogonal transformation
on $\ell$ as the complexification of it.
So we have proved Theorem \ref{thm:trivial} together
with the help of local sections (\ref{local-section-1}, \ref{local-section-2}).  
\end{proof}

\begin{cor}The Lagrangian-Grassmannian
$\Lambda(H)$ is an infinite dimensional differentiable manifold
modeled on the Banach space of bounded selfadjoint operators.
\end{cor}

\begin{proof}
Since we have an open covering 
$\{{\bf O}_{\lambda^{\perp}}\}_{\lambda \in\Lambda(H)}$ of
the Lagrangian Grassmannian $\Lambda(H)$, each
of which is isomorphic to a Banach space 
$\widehat{\cal{B}}(\lambda)$,
it will be enough to show the coordinate transformations of 
this covering are "{\it differentiable}" in a suitable sense.
Of course the Banach spaces  
$\widehat{\cal{B}}(\lambda)$ are
all isomorphic to a typical one.  

Let $G_{\lambda }:\hat{\cal{B}}(\lambda) 
\rightarrow {\bf O}_{\lambda^{\perp}}$
be the map in Proposition \ref{prop:local chart}, then 
{}from the expression of $\cal{P}_{G_{\lambda}(A)}$
(see (\ref{eq:projection})) 
we know 
the compositions of maps from $\widehat{\cal{B}}(\lambda)$
to $\cal{B}(H)$,
$$
\widehat{\cal{B}}(\lambda) \ni A\mapsto G_{\lambda }(A)
\mapsto \cal{P}_{G_{\lambda}(A)} \mapsto \cal{P}_{G_{\lambda} (A)}
+\cal{P}_{\mu^{\perp}}\in \cal{B}(H)
$$
is a "differentiable" map.

If $G_{\lambda }(A)\in {\bf O}_{\mu^{\perp}}$, that is,
$G_{\lambda }(A)=G_{\mu}(B)$ with an operator $B\in
\widehat{\cal{B}}(\mu)$,
then from the relation 
(see (\ref{eq:inverse}))
\begin{align*}
&(\cal{P}_{G_{\mu}(B)}+\cal{P}_{\mu^{\perp}})^{-1}(J(y))\\
&=(\cal{P}_{G_{\lambda }(A)}+\cal{P}_{\mu^{\perp}})^{-1}(J(y))\\
&=- B(y)+J(y) \ (y\in\mu),
\end{align*}
it will be apparent that the coordinate transformation
$ G_{\mu}^{-1}\circ G_{\lambda}$ : 
$A \mapsto B = J - (\cal{P}_{G_{\lambda}(A)}
+\cal{P}_{\mu^{\perp}})^{-1}\circ J$ is
a differentiable map between open sets in Banach 
space $\widehat{\cal{B}}(\lambda )$ 
and $\widehat{\cal{B}}(\mu)$.
\end{proof}

%

\subsection{Fredholm pairs and Fredholm operators}

Theorem \ref{thm:trivial} says that in the infinite
dimension we must work in a smaller space than the whole space of
Lagrangian
subspaces $\Lambda(H)$ to obtain a similar quantity to the Maslov
index in the finite dimensional case.  In this section we introduce
a notion, so called, Fredholm pairs and discuss relations of Fredholm operators and 
Fredholm pairs (\cite{Ka}).

Let $\ell_1$ and $\ell_2$ be two closed subspaces in $H$, then
first of all we recall the definition of $\ell_1$ and $\ell_2$ 
being a Fredholm pair.
\begin{defn}\label{def:Fred-pair}
We call two closed subspaces $\ell_1$ and $\ell_2$ being a Fredholm pair
if,  
\begin{enumerate}
\item 
\begin{equation*}
\dim (\ell_1 \bigcap \ell_2) \quad \text{is finite},
\end{equation*}
and
\item  
\begin{equation*}
\ell_1 +\ell_2 \quad\text{is closed and of finite codimensional in $H$}.
\end{equation*}
\end{enumerate}
\end{defn}

We give a relation of two notions ``Fredholm pair'' and ``Fredholm
operator'').

\begin{prop}\label{prop:fred and fred}
Let $\cal{P}_1: H \rightarrow H$ be the orthogonal projection operator 
with the image $\cal{P}_1(H) = \ell_1^{\perp}$.  Then $(\ell_1\,,\,\ell_2)$
is a Fredholm pair, if and only if, the restriction
$\cal{P}_1|_{\ell_2}$ of $\cal{P}_1$ to the space $\ell_2$
is a Fredholm operator, and
\begin{equation}\label{index-identity}
ind ~\cal{P}_1|_{\ell_2}= \dim \Ker\cal{P}_1|_{\ell_2}
- \dim \ell_1^{\perp}/\cal{P}_1(\ell_2) 
= \dim(\ell_1\bigcap\ell_2)-\dim(H/(\ell_1+\ell_2)).
\end{equation}
\end{prop}

\begin{proof}
In the algebraic sense we have $\Ker(\cal{P}_1|_{\ell_2}) =
\ell_1\bigcap\ell_2$  
and $H/(\ell_1+\ell_2) = \ell_1^{\perp}/\cal{P}_1(\ell_2)$
by the definition of the operator $\cal{P}_1|_{\ell_2}$.
Also we have that the closeness of $\ell_1+\ell_2$ 
and 
$\cal{P}_1|_{\ell_2}(\ell_2)$ is equivalent (A little bit general 
fact is proved in the next Lemma \ref{lem:closeness}).  
These prove the equivalence and we have (\ref{index-identity}).
\end{proof}

\begin{lem}\label{lem:closeness}
Let $T : H \rightarrow H'$ be a bounded surjective operator from a
Hilbert
space $H$ to a Hilbert space $H'$ and let $L$ 
be a closed subspace containing   $\Ker (T)$.  Then $T(L)$ is closed.
\end{lem}

\begin{proof}
Let $\pi$ be the orthogonal projection operator in $H$ with the image
= $\pi(H) = \Ker(T)$, and let $\tilde{T}$ be an isomorphism between
$H$ and $H'\oplus \Ker(T)$ defined by $\tilde{T}(x) =
T(x)\oplus\pi(x)$.
Then we have $\tilde{(T)}(L)$ is a closed subspace and we know that
$\tilde{(T)}(L) = T(L)\oplus \Ker(T)$.  This implies the closeness of
$T(L)$ in $H'$
\end{proof}

Next we generalize Lemma \ref{lem:transversal}, by which we give a
characterization of two Lagrangian subspaces being a Fredholm pair.

\begin{prop}\label{prop:Fred-pair} 
Let ${\mu},{\nu}\in \Lambda(H)$ and
let $\cal{P}_{\mu}$ (resp. $\cal{P}_{\nu}$) denote the orthogonal
projection operator of $H$ onto ${\mu}$ (resp. $\nu$). Then
$\cal{P}_{\mu}+\cal{P}_{\nu}$ is a Fredholm operator, if and only if
$({\mu},{\nu})$ is a Fredholm pair. 
\end{prop}

\begin{proof}
First we show 
$$\Ker (\cal{P}_{\mu}+\cal{P}_{\nu}) 
=\mu^{\perp}\bigcap\nu^{\perp}
$$ (see the end of the proof of Lemma
\ref{lem:transversal}). Since 
${\cal{P}}_{\mu}(x)+{\cal{P}}_{\nu}(x)=0$ implies that 
\[
  <x,{\cal{P}}_{\mu}(x)> 
= <x,-{\cal{P}}_{\nu}(x)> = - \|{\cal{P}}_{\nu}(x)\|^2 
= \|{\cal{P}}_{\mu}(x)\|^2 \,. 
\]
Hence ${\cal{P}}_{\mu}(x)={\cal{P}}_{\nu}(x)=0$, which shows that 
\begin{equation}\label{eq:kernel-iso}   
\Ker({\cal{P}}_{\mu}+{\cal{P}}_{\nu})={\mu}^\perp\bigcap
{{\nu}}^\perp=J(\mu)\bigcap J(\nu) = J({\mu}\bigcap {\nu})= (\mu+\nu)^{\perp}. 
\end{equation}

Now let
$\cal{P}_{\mu}+{\cal{P}}_{\nu}$ be a Fredholm operator. 
Then,
since $({\cal{P}}_{\mu}+{\cal{P}}_{\nu})(H)\subset {\mu}+{\nu}$,
and 

$\I({\cal{P}}_{\mu}+{\cal{P}}_{\nu})$ is closed and of finite

codimensional, so ${\mu}+{\nu}$ must be also
closed and of finite codimensional.
Hence together with the isomorphism (\ref{eq:kernel-iso}) we have proved that 
$({\mu},{\nu})$ is a Fredholm pair.

\medskip 
Next assume that $({\mu},{\nu})$ is a Fredholm pair, and we prove
$\cal{P}_{\mu}+{\cal{P}}_{\nu}$ is a Fredholm operator. 

{}From Proposition \ref{prop:fred and fred} we have 
$\cal{P}_{\mu}(\nu^{\perp})$ 
(resp. $\cal{P}_{\nu}(\mu^{\perp})$) 
is a finite codimensional closed 
subspace in $\mu$ (resp. $\nu$). Since $\dim (\mu\bigcap\nu)< \infty$,
in the direct sum $\mu\oplus\nu$ 
the subspace $\cal{P}_{\mu}(\nu^{\perp})\oplus\cal{P}_{\nu}(\mu^{\perp})$
+ $\{x\oplus -x\,|\, x\in\mu\bigcap\nu\}$ is still closed.
Consequently the subspace 
$\cal{P}_{\mu}(\nu^{\perp})+\cal{P}_{\nu}(\mu^{\perp})$ is closed
and finite codimensional in $\mu+\nu$.
Hence 
the image $(\cal{P}_{\mu}+\cal{P}_{\nu})(H)$ is a 
finite codimensional closed subspace in $\mu+\nu$, because it includes
the finite codimensional closed subspace 
$\cal{P}_{\mu}(\nu^{\perp})+\cal{P}_{\nu}(\mu^{\perp})$.
In fact it
coincides with $\mu+\nu$, since it is closed and
$(\cal{P}_{\mu}+\cal{P}_{\nu})(H)^{\circ}= \mu\bigcap\nu$.
Now we have proved that 
$\Ker (\cal{P}_{\mu}+{\cal{P}}_{\nu}) = J(\mu\bigcap\nu)$ and 
$\I(\cal{P}_{\mu}+{\cal{P}}_{\nu}) =\mu+\nu$, which
shows the operator $\cal{P}_{\mu}+{\cal{P}}_{\nu}$ 
is a Fredholm operator.   
\end{proof}

\subsection{Fredholm-Lagrangian-Grassmannian}

 We fix a Lagrangian subspace $\lambda$ and introduce 
a subspace of $\Lambda(H)$,
so called, {\it Fredholm-Lagrangian-Grassmannian} with respect to
$\lambda$. 

\begin{defn}\label{def:FLG}
The Fredholm-Lagrangian-Grassmannian of $H$ with respect to a
Lagrangian subspace $\lambda$
is defined as
\begin{equation}
\cal{F}\Lambda_{\lambda}(H)=\{\mu\in\Lambda (H)\;|\;(\mu,\lambda)\;
\text{is a Fredholm pair}\}.
\end{equation}
\end{defn}

\begin{defn}
We call the subset 
\begin{equation}
{\mathfrak M}_{\lambda}(H) = \{\mu\in\cal{F}\Lambda_{\lambda}(H)\;|\; \mu\bigcap\lambda\ne
\{0\}\}
\end{equation}
the Maslov cycle with respect to $\lambda$.
\end{defn}

\begin{notation}
$\cal{F}\Lambda_{\lambda}^{\,(0)}(H)=
\{\theta\in\cal{F}\Lambda_{\lambda}(H)\,|\, \theta 
~\text{is transversal to}~\lambda\}=\cal{F}\Lambda_{\lambda}(H)
\backslash{\mathfrak M}_{\lambda}(H)$ (= ${\bf O}_{\lambda}$, see 
Notation \ref{not:transversal}).
\end{notation}

\begin{rem}
\begin{enumerate}
\item
In the finite dimensional case, the subset ${\mathfrak M}_{\lambda}(H)$ is
a singular cycle whose homology class is a generator of the codimension one
homology group $H_{\frac{n(n+1)}{2}-1}(\Lambda(H),\mathbb{Z})$, 
where we put $\dim H = 2n$.
\item
As we proved in Proposition \ref{prop:local chart}
the subset $\cal{F}\Lambda_{\lambda}(H)\backslash {\mathfrak M}_{\lambda}(H)$
= $\cal{F}\Lambda_{\lambda}^{\,(0)}(H)$ 
is isomorphic to the space of bounded selfadjoint operators on
$\lambda^{\perp}$.
\end{enumerate}
\end{rem}

First
we study how the Fredholm-Lagrangian-Grassmannian $\cal{F}\Lambda_{\lambda}(H)$
depends on the space ${\lambda}$. 
In the finite dimensional case, it is clear that $\cal{F}\Lambda_{\lambda}(H)$ =
$\Lambda(H)$.
In the infinite dimension, $\cal{F}\Lambda_{\lambda}(H)$ is 
an open subset of $\Lambda(H)$. Openness follows {}from 
Proposition \ref{prop:open} and 
Proposition \ref{prop:Fred-pair}, and it can not include $\lambda$
itself. However we can prove

\begin{prop}
$\lambda$ can be approximated by a
sequence in $\cal{F}\Lambda_{\lambda}(H)$, i.e., 
$\lambda\in\partial\mathcal{F}\Lambda_{\lambda}(H)$(= the boundary). 
\end{prop}
\begin{proof}
Let $ A:\lambda \rightarrow $
 $\lambda$  be a bounded selfadjoint operator 
and assume that $A$ is an isomorphism. 
Then for all $\epsilon>0$, the
Lagrangian subspace
\[
G_{\epsilon\cdot A} = \{x+\epsilon J A(x) \mid x\in \lambda\}
\]
is transversal with both of $\lambda$ and $\lambda^{\perp}$.
Since $\epsilon  A$ converges to $0$ in $\widehat{\cal{B}}(\lambda)$
when $\epsilon \rightarrow 0$,
we know that the orthogonal projection operator
$\cal{P}_{G_{\epsilon\cdot A}}$ 
onto the graph of $\epsilon\cdot J \circ  A$ converges to $\cal{P}_\lambda$.
Hence we have
\begin{equation}\label{e:closure_point}
\lambda\in\overline{\cal{F}\Lambda_{\lambda}(H)} 
\setminus \cal{F}\Lambda_{\lambda}(H).
\end{equation}
\end{proof}

Let $\lla$ and $\mu$ in $\LLa(H)$ and assume that 
\begin{equation}\label{eq:mod-compact}
\mu =U(\lla)~\text{with}
\end{equation}
$$
U=\Id + K\in \cal{U}(H_J) 
~\text{is of the form} ~\Id+ ~\text{{\em compact operator}},
$$
then 
\begin{prop}
$$
\cal{F}\LLa_{\lla}(H) =\cal{F}\LLa_{\mu}(H). 
$$
\end{prop}
\begin{proof}
By Proposition \ref{prop:Fred-pair}, $\nu\in\cal{F}\LLa_{\lla}(H)$ 
if and only if 
$\cal{P}_{\lla} + \cal{P}_{\nu}$ is a Fredholm operator.  
{}From the assumption
$\cal{P}_{\mu} + \cal{P}_{\nu}$ =$\cal{P}_{U(\lla)} + \cal{P}_{\nu}$ 
= $U\circ\cal{P}_{\lla}\circ U^{-1} + \cal{P}_{\nu}$
=$(\Id + K)\circ\cal{P}_{\lla}\circ (\Id +K^*) + \cal{P}_{\nu}$ 
=$\cal{P}_{\lla} + \cal{P}_{\nu}+ ~\text{compact operator}$.
Hence if $\nu\in\cal{F}\LLa_{\lla}(H)$, then 
$\nu\in \cal{F}\LLa_{\mu}(H)$.  Since $(\Id+K)^{-1} =\Id+K^*$,
by the same way
we have $\cal{F}\LLa_{\mu}(H)\subset\cal{F}\LLa_{\lla}(H)$. 
\end{proof}

\begin{defn}
We denote by $\cal{U}_{\text{res}}(H_J)$ the subgroup of
$\cal{U}(H_J)$
consisting of such operators that
$$
\cal{U}_{\text{res}}(H_J)=\{\Id + ~\text{compact operator}\}.    
$$
\end{defn}

\begin{cor}
The group $\cal{U}_{\text{res}}(H_J)$ acts on $\cal{F}\LLa_{\lambda}(H)$
and $\cal{U}_{\text{res}}(H_J)(\lambda)\subset
\partial(\cal{F}\LLa_{\lambda}(H))$,
that is, the orbit of the element $\lambda$ is also included in the
boundary of $\cal{F}\LLa_{\lambda}(H)$.
\end{cor}

As a special case of the relation (\ref{eq:mod-compact})
we introduce an equivalence relation on the space $\Lambda(H)$:
\begin{defn}\label{def:eqivalence}
We call $\lambda$ and $\mu$ $\in\Lambda(H)$ almost coincide, if
\[
\dim  {\lambda}/({\lambda} \bigcap {\mu}) < +\infty.
\]
and denote
\[
  {\lambda}\sim {\mu}, 
\]
when two Lagrangian subspaces $ {\lambda}$ and ${\mu}$ almost coincide.
\end{defn}

It will be easy to prove that this is in 
fact an equivalence relation.  Note that in this case
\[ 
\dim {\lambda}/({\lambda} \bigcap {\mu}) = \dim {\mu}/({\lambda} \bigcap {\mu}),
\]
and in fact
\begin{prop}\label{prop:compact-equivalence}
Let $\lla \sim\mu$, then there exists a unitary operator $U$
of the form $U=\Id +~\text{compact operator}$ such that
$\mu =U(\lla)$.
\end{prop}
\begin{proof}
Since $\lla$ and $\mu$ are Lagrangian subspaces,
the sum of the complex subspaces spanned 
by $\lla\bigcap\mu$ and $\lla\bigcap(\lla\bigcap\mu)^{\perp}$
is an orthogonal sum of $H_J$, 
and so $(\lla\bigcap(\lla\bigcap\mu)^{\perp})\otimes\mathbb{C}$ =
$(\mu\bigcap(\lla\bigcap\mu)^{\perp})\otimes\mathbb{C}$ in $H_J$.  
Hence we can find such an unitary operator $U$ that is identity on 
the subspace $(\lla\bigcap\mu)\otimes\mathbb{C}$. Hence we can take
$U=Id + K$ with $K$ being a finite rank operator.   
\end{proof}

\begin{prop} \label{prop:equivalence}
Let ${\lambda}\in\Lambda(H)$ and let $W \subset {\lambda}$ be a finite
codimensional closed
subspace in ${\lambda}$. Then for ${\mu}\in \Lambda(H)$, the pair 
$({\lambda},{\mu})$ is a Fredholm pair, if and only if,
$(W,{\mu})$ is a Fredholm pair. 
\end{prop}

We denote by $\cal{F}\Lambda_W(H)$
\begin{equation}
\cal{F}\Lambda_W(H)= 
\{\mu \in\Lambda (H)\,|\, (W\,,\,\mu) ~\text{is a Fredholm pair}\}. 
\end{equation}

\begin{proof}[Proof of Proposition \ref{prop:equivalence}]
We prove $\cal{F}\Lambda_W(H)$ = $\cal{F}\Lambda_{\lambda}(H)$.

Let $\mu\in\cal{F}\Lambda_W(H)$. Then, since the map 
\[
(\lambda\bigcap\mu)/(W\bigcap\mu)\rightarrow \lambda/ W
\]
is injective, we have
\[
\dim(\lambda\bigcap\mu) \leq \dim (\lambda/ W) + \dim (W\bigcap\mu),
\]
and the space $\lambda +\mu$ is a finite dimensional extension of the
closed subspace $W+\mu$.  Hence $\lambda$ and $\mu$ is a Fredholm
pair.
 
Now let $\mu\in\cal{F}\Lambda_{\lambda}(H)$.
In the short exact sequence
\[
  0\rightarrow  \lambda\bigcap {\mu} \stackrel{j}{\rightarrow}
{\lambda}\oplus {\mu} \stackrel{\tau}{\rightarrow}
{\lambda}+{\mu} \rightarrow  0,
\]
where $j(a)= a\oplus -a\in H\oplus H$ and $\tau(a\oplus
b)=a+b$, we have
\[
  \tau^{-1}(W+{\mu})=W\oplus {\mu} + j({\lambda}\bigcap {\mu}). 
\]
Hence $W+\mu$ must be closed in $\lambda+\mu$, so is in $H$. 
Also we have at once $\dim W\bigcap\mu <\infty$.  These
proves the coincidence
$\cal{F}\Lambda_W(H)$ = $\cal{F}\Lambda_{\lambda}(H)$.
\end{proof}

\begin{cor}\label{cor:coincidence} 
If ${\lambda}\sim {\mu}$, then $\cal{F}\Lambda_{\lambda}(H)$ = 
$\cal{F}\Lambda_{\mu}(H)$.
\end{cor}
\begin{proof}
By Proposition \ref{prop:compact-equivalence} we already know this, 
but also by
putting $W =\lambda\bigcap\mu$ in the proof of 
Proposition \ref{prop:equivalence} we can prove the coincidence.
\end{proof}

\begin{rem}
Since in the proof of the above proposition
we did not use any particular properties of Lagrangian subspaces, 
the above coincidence holds for any Fredholm pair $(L_1,L)$ and
$(L_2,L)$, where $L_2$ is a finite codimensional closed subspace in
$L_1$.  
\end{rem}

Finally we note
\begin{prop}\label{prop:existence of transversal subspaces}
We have an open covering ${\cal F}\LLa_{\lla}(H)$ = 
$\bigcup\limits_{\mu\sim\lla}{\bf O}_{\mu}$, and each ${\bf O}_{\mu}$
is open dense in ${\cal F}\LLa_{\lla}(H)$. Hence 
$\bigcap\limits_{i=1}^{\infty}{\bf O}_{\mu_i}$ (each $\mu_i \sim\lla$)
is not empty. In other words, for any given countable number of
Lagrangian subspaces $\{\mu_i\}_{i=1}^{\infty}$ 
each of which is equivalent to a fixed Lagrangian
subspace $\lla$ there exists a Lagrangian subspace which is
transversal to each $\mu_i$. 
\end{prop}


\subsection{Souriau map and the universal Maslov cycle}
When we fix a $\lambda \in \Lambda(H)$ then we have
an identification 
\begin{equation}\label{eq:decom}
\begin{array}{l@{}l}
H_{J}=\lambda\oplus\lambda^{\perp}= & \lambda\oplus J\lambda\cong
\lambda\otimes\C\\
& x+Jy\mapsto x\otimes 1+y\otimes\sqrt{-1}.
\end{array}
\end{equation}

We denote by $\tau_{\lambda}$ the complex conjugation 
under this identification:
\[ 
\tau_{\lambda}(x+J(y))=x-J(y), x, y \in \lambda.
\]

It will be easy to show the following relation:
\begin{equation}\label{eq:conjugation and projection}
\tau_{\lla}= 2\cal{P}_{\lla} -\Id.
\end{equation}

Any $U \in \cal{U}(H_J)$ can be expressed as
\[
U=X+\sqrt{-1}Y
\]
with $ X,Y \in \cal{B}(\lambda)$ in such a way that
\begin{eqnarray*}
U(x \otimes 1+y \otimes \sqrt{-1})&=&(X(x)-Y(y))\otimes 1+(X(y)+Y(x))\otimes\sqrt{-1}\\
                                 &=&X(x)-Y(y)+J(X(y)+Y(x)),
\end{eqnarray*}
and $X$, $Y$ satisfy the relations :  
\[
X{\,^tY}=Y{\,^tX}, \ {\,^tY}X={\,^tX}Y.
\]
\[
X{\,^tX}+Y{\,^tY}=\Id, \ {\,^tX}X+{\,^tY}Y=\Id.
\]

For $\lla\in\LLa(H)$ we denote by $\theta_{\lla}$ an anti-group
isomorphism $\cal{U}(H_J) \to \cal{U}(H_J)$ defined by
\begin{equation}\label{eq:lifted conjugation}
  \theta_{\lla}(U) = \tau_{\lla}\circ U^* \circ \tau_{\lla},
\end{equation}
then $\theta_{\lla}(U) = {\,^tX}+\sqrt{-1}{\,^tY}$.

Note that if $Y\ne 0$ $\theta_{\lla}(U)\ne {\,^tU}$,
where we mean ${\,^tU}$ is a transposed operator when we regard $U$
as a real linear operator.


Then we have 
\begin{equation}
\{U\in \cal{U}(H_J)\,|\, \theta_{\lla}(U)=U^{-1}\}
=\cal{O}(\lla) (= \text{orthogonal group on}~\lla).
\end{equation}
Hence the map $\cal{U}(H_J) \to \cal{U}(H_J)$, 
$U \mapsto U\circ\theta_{\ell}(U)$
induces a continuous map (see Corollary \ref{cor:quotient topology})

\begin{equation}\label{eq:Souriau-map}
\begin{array}{r@{\,}l}
\cal{S}_{\ell}:\LLa(H) \longrightarrow \cal{U}(H_J)\\
\mu =U(\ell^{\perp}) \mapsto \cal{S}_{\ell}(\mu)= & U\circ\theta_{\ell}(U).
\end{array}
\end{equation}

We call this map as ``{\it Souriau map}'' henceforth.

{}From the relation (\ref{eq:conjugation and projection}) we have
an expression of the Souriau map in terms of projection operators
corresponding to Lagrangian subspaces:
\begin{prop}\label{prop:expression of Souriau map}
$\cal{S}_{\ell}(\mu) = (\Id-2\cal{P}_{\mu})(2\cal{P}_{\ell}-\Id)
=-\tau_{\mu}\circ\tau_{\ell}$.
\end{prop}

\begin{cor}
Let $\lambda,~\mu,~\nu$ be three Lagrangian subspaces, then
\begin{equation}
\mathcal{S}_{\mu}(\nu)\circ
\mathcal{S}_{\lambda}(\mu)=-\mathcal{S}_{\lambda}(\nu).
\end{equation}
\end{cor}

{}From the relations (\ref{eq:conjugation and projection}),
(\ref{eq:lifted conjugation}) and 
Proposition \ref{prop:expression of Souriau map}
we have 
\begin{prop}
The maps 
\begin{equation}\label{eq:Souriau-map2}
\begin{array}{r@{\,}l}
\qquad\cal{U}(H_J)\times\Lambda(H) & \rightarrow \cal{U}(H_J)\\
(\ell\,,\,U)\qquad &\mapsto U\circ \theta_{\ell}(U)
=U\circ(2{\cal P}_{\lla}-\Id)\circ U^*\circ (2{\cal P}_{\lla}-\Id),
\end{array}
\end{equation}
and 
\begin{equation}\label{eq:Souriau-map3}
\begin{array}{r@{\,}l}
\LLa(H)\times\LLa(H) & \rightarrow \cal{U}(H_J)
\qquad\qquad\qquad\qquad\qquad\qquad\\
(\ell\,,\,\mu)\qquad &\mapsto \cal{S}_{\ell}(\mu)=
(\Id-2\cal{P}_{\mu})(2\cal{P}_{\ell}-\Id)
\end{array}
\end{equation}
are continuous.
\end{prop}

By Proposition \ref{prop:expression of Souriau map}, 
\begin{equation}
U\circ \cal{S}_{\ell}(\mu)\circ U^* =\cal{S}_{U(\ell)}(U(\mu)),
\end{equation}
that is, the following commutative diagram:
\begin{prop}
\begin{equation}
\begin{CD}
\LLa(H) @> {\cal{S}_{\ell}} >>\cal{U}(H_J)\\
@V{U}VV                                  @VV{\text{Ad}_U}V\\
\LLa(H)  @>> {\cal{S}_{U(\ell)}} >\cal{U}(H_J).\\
\end{CD}
\end{equation}
\end{prop}

In particular, 
when $U\in\mathcal{U}_{res}(H_J)$
we have a commutative diagram :
\begin{prop}
\begin{equation}
\begin{CD}
\mathcal{F}\LLa_{\ell}(H)=\mathcal{F}\LLa_{U^{-1}(\ell)}(H) 
@> {\cal{S}_{U^{-1}(\ell)}} >>\cal{U}(H_J)\\
@V{U}VV                                  @VV{\text{Ad}_{U}}V\\
\mathcal{F}\LLa_{\ell}(H)\qquad @>> {\cal{S}_{\ell}}>\cal{U}(H_J).\\
\end{CD}
\end{equation}
\end{prop}

Here we remark the adjoint operator of an anti-linear operator:
let $T$ be an anti-linear operator on a complex Hilbert space $H$ with
a Hermitian inner product $(\bullet,\bullet)$, then the adjoint
operator $T^*$ is defined by the relation
$(T(z),~w) = (T^*(w),~z)$ ($z\,,\,w\in H$).  Then
$T^{\ast}$  is again an anti-linear operator and
we have a composition formula with a linear or anti-linear
operator $L$ : 
$(T\circ L)^* = L^*\circ T^*$. 

Now $\tau_{\lla}$ is anti-linear and
we have by a direct calculation
$$
\tau_{\lla}^* =\tau_{\lla}, 
~\text{that is} ~(\tau_{\lla}(z),w)_J = (\tau_{\lla}(w),z)_J.
$$

{}From this fact $\theta_{\lla}^{~2} =\Id$, in other words, 
$\theta_{\lla}$ is an anti-linear involution on ${\cal B}(H_J)$.

By the above remark and 
the expression of the Souriau map 
(Proposition \ref{prop:expression of Souriau  map})
we have
\begin{prop}\label{prop:symmetry of Souriau map}
\begin{equation}
\cal{S}_{\lla}(\mu)^* =\cal{S}_{\mu}(\lla). 
\end{equation}
\end{prop}

{\it We call the restriction of the Souriau map to $\cal{F}\LLa_{\lla}(H)$
also {\em Souriau map} always}.

Now for a fixed $\lla$, we put 
$\cal{U}_{\lambda}(H_{J})=\rho_{\lla}^{-1}(\cal{F}\Lambda_{\lambda}(H))$, 
where $\rho_{\lla} :\cal{U}(H_{J})\to \cal{F}\Lambda_{\lambda}(H)$,
$\rho_{\lla}(U)=U(\lla^\perp)$. Then

\begin{prop}\label{prop:real part}
Let $ U \in \cal{U}(H_J)$, then 
$U=X+\sqrt{-1}Y \in \cal{U}_{\lambda}(H_J)$, 
if and only if, $X\in \cal{B}(\lla)$ is a Fredholm operator.
\end{prop}
\begin{proof}
Let $U \in \cal{U}(H_J)$, and put $\mu =U(\lla^{\perp})$.  Then the inclusion 
map $\lambda^{\perp} \to  H=\lambda+\lambda^{\perp}$ induces the isomorphism
\[
(\lambda+\lambda^{\perp})/(\lambda+\mu)\cong 
\lla^{\perp}/J(X(Y^{-1}(\lla))) =\lambda^{\perp}/J(X(\lambda))
\cong \lla/X(\lla).
\]
Also 
\[
\lambda \cap U(\lambda^{\perp}) \cong \Ker X.
\]
These shows the assertions.
\end{proof}

Let $\mu\in \cal{F}\LLa_{\lla}(H)$ and $U(\lla^{\perp})=\mu$,
then by the definition of the Souriau map 
$\cal{S}_{\lla}(\mu)= U\circ \theta_{\lla}(U)$ and 
{}from above Proposition \ref{prop:real part}, we have
\begin{prop}\label{prop:fred unitary}
\begin{equation}
U\circ\theta_{\lla}(U) +\Id
\end{equation}
is a Fredholm operator.
\end{prop}
\begin{proof}
Let $U = X+\sqrt{-1}Y$, with $X\,,\, Y \in \cal{O}(\lla)$, then
$$
U\circ\theta_{\lla}(U) +\Id = 2X\circ\theta_{\lla}(U),
$$
and this shows the Fredholmness of the operator 
$U\circ\theta_{\lla}(U) +\Id$. 
\end{proof}
 
Let $\mu \in \LLa(H)$, then from the relation that 
an element $z=x+J(y)$ ($x,y\in \lla$) is in $\mu$ if and only if
$-z=W_{\mu}(\tau_{\lambda}(z))$, we have 

\begin{equation}\label{eq:kernel}
\Ker(W_{\mu}+\Id)=(\mu\cap\lambda)\otimes\C\cong(\mu\cap\lambda)\oplus
J(\mu\cap\lambda).
\end{equation}
Hence
\begin{prop}\label{prop:dimension}
For any $\mu\in\cal{F}\Lambda_{\lambda}(H)$ and any $U\in\cal{U}_{\lambda}
(H_J)$ with $\mu=U(\lambda^{\perp})$, 
\begin{equation*}
\dim_{\mathbb{R}}\,(\mu\cap\lambda)=\dim_{\C}\,\Ker(W_{\mu}+\Id).
\end{equation*}
\end{prop}

Let us now consider the space
\begin{equation}\label{eq:2.8}
\cal{U}_{\cal{F}}(H_J)=\{U\in\cal{U}(H_J)\;|\;\text{$U+\Id$ is a Fredholm
operator}\}
\end{equation}
and a subset
\begin{equation}\label{eq:2.9}
\cal{U}_{\mathfrak{M}}(H_J)=\{U\in\cal{U}_{\cal{F}}(H_J)\;|\;\Ker(U+\Id)\ne
\{0\}\}
\end{equation}
which by the preceding Proposition \ref{prop:dimension} we can regard 
as a kind of the universal Maslov cycle:
\begin{prop}
For any $\lambda$,
$\cal{S}_{\lambda}^{-1}(\cal{U}_{\mathfrak{M}}(H_J)) = {\mathfrak M}_{\lambda}(H)$.
\end{prop}

Now we state the fundamental property for discussing the Maslov index in
the infinite dimension:

\begin{thm}\label{thm:pi_1}
\begin{enumerate}
\item
\begin{equation*}
\pi_1(\cal{F}\Lambda_{\lambda}(H))\simeq\mathbb{Z},\;
\end{equation*}
\item
\begin{equation*}
\pi_1(\cal{U}_{\cal{F}}(H_J))\simeq\mathbb{Z}.
\end{equation*}
\item
The induced map 
\begin{equation*}
(\cal{S}_{\lambda})_{*}:\;\pi_1(\cal{F}\Lambda_{\lambda}(H))
\rightarrow\pi_1(\cal{U}_{\cal{F}}(H_J)) 
\end{equation*}
is an isomorphism.
\end{enumerate}
\end{thm}

We give the proof of this Theorem in the next subsection 
by the method of the finite dimensional reduction.

\subsection{Proof of Theorem \ref{thm:pi_1} (a)}

\begin{notation}\label{n:notation}
\begin{enumerate}
\item
Let $\lambda\in\Lambda(H)$.
We denote by $Sub_{fin}(\lambda)$ the
set of all closed subspaces 
$W\not=\lambda$ of $\lambda$ of finite codimensions.  
\item
Let $W$ be a closed isotropic subspace such that 
$\dim W^{\circ}/W <\infty$. We denote by
$\Lambda(W,H)$ for the set of Lagrangian subspaces of $H$
which contains $W$.
\end{enumerate}
\end{notation}

We prove

\begin{thm}\label{thm:finite reduction-1}
Let $\lambda\in\Lambda(H)$ and $W\in Sub_{fin}(\lambda)$.

(a) The inclusions 
\[
\cal{F}\Lambda_{W}^{\,(0)} 
=\{\theta\in\cal{F}\Lambda_{\lambda}(H)\mid\theta\cap
W=\{0\}\,\} \hookrightarrow \cal{F}\Lambda_{\lambda}(H),
\]
define an isomorphism
\[
ind-\lim_{W\to \{0\}} \pi_1(\cal{F}\Lambda_{W}^{\,(0)}(H)) 
\stackrel{\sim}{\rightarrow}
\pi_1(\cal{F}\Lambda_{\lambda}(H)).
\]

(b)  There is a natural isomorphism
\[
\pi_1(\cal{F}\Lambda_{W}^{\,(0)}(H)) \stackrel{\sim}{\rightarrow} 
\pi_1(\Lambda(W^{\circ}/W)) \cong \Z \]
for each $W\in Sub_{fin}(\lambda)$.
\end{thm}

\bigskip

{\it By combining (a) and (b) we obtain}
Theorem \ref{thm:pi_1} (a).


\bigskip

The proof of Theorem \ref{thm:finite reduction-1} will follow from two
Propositions below which will be of independent interest. 
First we shall prove

\begin{prop}\label{p:compact}
Let $K\subset\cal{F}\Lambda_{\lambda}(H)$ be a compact
set. Then there exists a $W\in
 Sub_{fin}(\lambda)$ such that $\mu\cap W= \{0\}$ for all
$\mu\in K$.
\end{prop}

\begin{proof} Let $\mu_0\in K$. Then the sum of the orthogonal
projections $\cal{P}_{\lambda}+\cal{P}_{\mu_0}$ 
is a Fredholm operator by Proposition \ref{prop:Fred-pair} 
and we have
\[
\Ker (\cal{P}_{\lambda}+\cal{P}_{\mu_0}) = J(\lambda\cap\mu_0).
\]
Let
\[
h=\Bigl(J(\lambda\cap\mu_0)\Bigr)^\perp = \lambda+\Bigl(\lambda^\perp\cap
\bigl(J(\lambda\cap\mu_0)\bigr)^\perp\Bigr) .
\]
Then the operator $\cal{P}_{\lambda}+\cal{P}_{\mu_0}$ is injective on $h$
and its range $\lambda+\mu_0$ is closed. Hence there exists an open
neighborhood $U$ of $\mu_0$ in $\cal{F}\Lambda_{\lambda}(H)$ such that
$\cal{P}_{\lambda}+\cal{P}_{\mu}$ is injective on $h$ for all $\mu\in K\cap
U$. Since $K$ compact, a finite set
$U_1,\dots, U_N$ of such neighborhoods covers the whole of
$K$. Then
\[
W=\bigcap_{j=1}^N\Bigl((\lambda\cap\mu_j)^\perp\cap\lambda\Bigr)
\]
satisfies our requirement for suitable choices of $\mu_j\in
U_j\cap K$.
\end{proof}

\bigskip

The next proposition gives a property of $\cal{F}\Lambda_{\lambda}(H)$
relating with the finite dimensional reduction of the Maslov index.

\begin{prop}\label{p:fibre_bundle}
Let $W\in Sub_{fin}(\lambda )$, then the mapping
\[
\begin{matrix}
\rho_W: & \cal{F}\Lambda_W^{\,(0)}(H) & \rightarrow & \Lambda (W^{\circ}/W) \\
\ & \mu & \mapsto & \bigl((\mu\cap W^{\circ})+W\bigr)/W
\end{matrix}
\]
defines a fiber bundle.
\end{prop}

The proof of this proposition is given by proving lemmas below.

\begin{lem}
\begin{enumerate}
\item  
Let $H,\lambda,W$ be as above and let $\theta\in\Lambda(W,H)$,
i.e., $\theta$ is a Lagrangian subspace including $W$. 
Then  
\[
U_{\theta}=\{\mu\in\cal{F}\Lambda_W^{\,(0)}(H) \mid \mu\cap\theta=\{0\}\,\} \]
is an open subset of the total space 
$\cal{F}\Lambda_W^{\,(0)}(H)$ and we
have \[
\bigcup\limits_{\theta \in\Lambda(W,H)} U_{\theta} 
= \cal{F}\Lambda_W^{\,(0)}(H).
\]
\item 
Let $\bar\theta=\theta/W\in\Lambda(W^{\circ}/W)$. Then the
set 
\[
U_{\bar\theta}=\{{\ell}\in \Lambda(W^{\circ}/W) \mid
{\ell}\cap\bar\theta=\{0\}\,\}
\]
is an open subset of the Lagrangian Grassmannian manifold
$\Lambda(W^{\circ}/W)$
and the union of
all such subsets covers $\Lambda(W^{\circ}/W)$.
\item  The mapping 
\[
{\rho}_W:U_{\theta}\to U_{\bar\theta}
\]
is surjective.
\end{enumerate}
\end{lem}

\begin{proof} Since $\mu\in U_{\theta}$ is transversal with $\theta$,
openness of $U_{\theta}$ follows from Lemma\ref{lem:transversal}.
For a given $\mu\in\cal{F}\Lambda_W^{\,(0)}(H)$ one finds easily a
$\theta = W+ L\in\Lambda (H)$ with  
$\theta\cap\mu=\{0\}$, by taking a suitable Lagrangian subspace $L$ in
$(\lambda\cap W^{\perp})\oplus J(\lambda\cap W^{\perp})$. 
That gives the claimed open covering
and (b) and (c) can be seen easily.
\end{proof}

\medskip

Again let $W\in Sub_{fin}(\lambda)$ and $\theta
\in\mathcal{F}\Lambda_{\lambda}(H),~ \theta \supset W$ and
we decompose $H$ into four mutually orthogonal subspaces:
\begin{equation}\label{e:decomposition}
\begin{matrix}
H & = &\                  & \theta &\ &+&  
\                    & J(\theta) &\ \\
\    & = & W^\perp\cap\theta &+      & W&+&
J(W^\perp\cap\theta) & + & J(W).
\end{matrix}
\end{equation}


\begin{lem}\label{l:alpha_gamma}
Let $\mu\in U_{\theta}$. Then
there exist linear mappings
\begin{align*}
a &: J(W^\perp\cap\theta) \to W^\perp\cap\theta\\
g  &: J(W^\perp\cap\theta) \to W
\end{align*}
such that each $z\in\mu\cap W^0$ can be written in the form \[
z=x+a(x)+g(x) \quad\text{ with $x\in J(W^\perp\cap\theta)$}.
\] \end{lem}

\begin{proof}
Since $\mu$ intersects $\theta$ transversally, there is a map
$A:J(\theta)\to\theta$ such that $A\circ J$ self--adjoint on
$\theta$ and $\mu = \{u+Au\mid u\in J(\theta)\}$. We decompose
$u=x+y$ with  $x\in J(W^\perp\cap\theta)$ and $y\in J(W)$
according to the decomposition of $J\theta$ in
\eqref{e:decomposition}. With regard of that decomposition,
the mapping $A$ can be written as a $2 \times 2$  matrix
$\begin{pmatrix} a & b\\ g & d \end{pmatrix}$. More
explicitly, we have
\[
Au=a(x) + b(y)  + g(x) + d(y),
\]
where
\begin{align*}
a &: J(W^\perp\cap\theta) \to W^\perp\cap\theta\\
b &: J(W) \to W^\perp\cap\theta\\
g  &: J(W^\perp\cap\theta) \to W\\
d  &: J(W) \to W\,.
\end{align*}
We notice that
\begin{equation}\label{e:relations}
a\circ J,\,d\circ J \text{are selfadjoint, and} 
^{\,t}\,(b\circ J)= g\circ J.
\end{equation}

Now, let $z\in\mu\cap W^0$. It can be written as
\[
z=u+Au= x+y+ a(x) + b(y) +g(x) +d(y).
\]
{}From the decomposition \eqref{e:decomposition} it follows that
the component $y$ in $J(W)$ must vanish. So,
\[
z= x+ a(x) + g(x).
\]
\end{proof}

\begin{cor}\label{c:footpoint}
Let $\lambda,W,\theta$ be as above. Let $\mu=\{u+Au\mid u\in J(\theta)\}
\in U_{\theta}$ with $A= \begin{pmatrix} a & b\\ g & d
\end{pmatrix}$ with respect to the decompositions $J(\theta)=
J(W^\perp\cap\theta)+J(W)$ and $\theta=
W^\perp\cap\theta\,+\,W$. As before, we identify $W^{\circ}/W$ with
$(W^\perp\cap\theta) + J(W^\perp\cap\theta)$. Then
\begin{equation}\label{e:footpoint_alpha}
{\rho}_W(\mu) = \{x+ a(x)\mid x\in J(W^\perp\cap\theta)\}.
\end{equation}
In particular, two $\mu,\mu'\in U_{\theta}$ belong to the same
fiber, i.e.,
${\rho}_W(\mu)= {\rho}_W(\mu')$, if and only if, $a=a'$.
\end{cor}

Now we prove Proposition \ref{p:fibre_bundle}.
\begin{proof}
We define a local trivialization on $U_{\bar\theta}$:
\begin{equation}\label{e:CD_trivialization}
\begin{matrix}
U_{\bar\theta}\times F & 
\quad\stackrel{\tau}{\rightarrow}& U_{\theta} \\
\quad &\ &\ \\
\ & \pi\searrow\quad & \quad\downarrow {\rho}_W \\
\quad &\ &\ \\
\ &\ &U_{\bar\theta}
\end{matrix}
\end{equation}
Here, $\pi$ denotes the projection onto the first component.
We take
\[
F = \cal{B}(JW,W^\perp\cap\theta) + \cal{B}_{sa}(JW,W)
\]
where $\cal{B}(JW,W^\perp\cap\theta)$ denotes the vector space of
bounded operators from $JW$ to $W^\perp\cap\theta$ and
$\cal{B}_{sa}(JW,W)$ the vector space of bounded operators from
$JW$ to $W$ which become selfadjoint on $W$ by combing
with $J$. For a fixed point  $L\in U_{\bar\theta}$ and a point
in the
fiber $(b\,,\,d)\in F$, we define
\[
\tau(L;b,g) =\{u+Au\mid A=\begin{pmatrix} a_L &b\\ g_b
&d \end{pmatrix},\ u\in J\theta\},
\]
with the decomposition $J\theta=J(W^\perp\cap\theta)+JW\}$.
The operator $a_L:J(W^\perp\cap\theta)\to W^\perp\cap\theta$
with $a_L\circ J$ selfadjoint is uniquely determined by the
condition
\[
L=\{x+a_L(x)\mid x\in J(W^\perp\cap\theta)\}.
\]
As a consequence, we get
$\tau$ is surjective and injective. By the definition of $a_L$ from
$L$ we get the commutativity of the diagram
\eqref{e:CD_trivialization}.
\end{proof}

Before proving Theorem \ref{thm:finite reduction-1} we remark the
a commutative diagram (\ref{e:CD_lim}).

Let us consider 
two spaces $W, W'\in Sub_{fin}(\lambda)$ with $ W'\subset W$. 
So
\[
\cal{F}\Lambda_{W}(H) = \cal{F}\Lambda_{W'}(H) 
= \cal{F}\Lambda_{\lambda}(H)
\]
and
\[
\cal{F}\Lambda_{W}^{\,(0)} \subset \cal{F}\Lambda_{W'}^{\,(0)}(H) 
\subset \cal{F}\Lambda_{\lambda}(H)\,.
\]

Recall that $\Lambda(W,H)$ denotes the set of Lagrangian
subspaces of $H$ which contain $W$, and then this space is
isomorphic with the Lagrangian Grassmannian
$\Lambda(W^{\circ}/W)$ in an obvious way:
\[
\begin{matrix}
\Lambda(W,H) & \stackrel{\sim}{\rightarrow}& \Lambda(W^{\circ}/W)\\
\theta & \mapsto & \theta/W,
\end{matrix}
\]
and a corresponding isomorphism for $W'$. Now let $C:I\to
\cal{F}\Lambda_{\lambda}(H)$ be a curve which is transversal to $W$. So, it
gives us the curve $C:I\to \cal{F}\Lambda_{W}^{\,(0)}(H)$. 
Then we have the
following commutative diagram :
\begin{equation}\label{e:CD_lim}
\begin{CD}
I              @>C>>     \cal{F}\Lambda_{\lambda}(H) \\
@V{C}VV                 @AA{\cup}A\\
\cal{F}\Lambda_{W}^{\,(0)}(H)@>>{\hookrightarrow}>
\cal{F}\Lambda_{W'}^{\,\,(0)}(H)\\
@V{{\rho}_W}VV             @VV{{\rho}_{W'}}V\\
\Lambda(W^{\circ}/W)     @>>>      \Lambda( W'^{\,\circ}/ W')\\
@A{\cong}AA                 @AA{\cong}A\\
\Lambda(W,H)@>>\stackrel{{\bf h}_{W,W'}}{}>\Lambda(W',H)\,.
\end{CD}
\end{equation}
\begin{proof}[Proof of Theorem \ref{thm:finite reduction-1}] 

By Proposition \ref{p:compact},  it will not be difficult to see
the mapping
\[
ind\lim_{W\to \{0\}} \pi_1(\cal{F}\Lambda_{W}^{\,(0)}(H)) 
\to \pi_1(\cal{F}\Lambda_{\lambda}(H))
\]
is naturally isomorphic(and also it is isomorphic 
for all homotopy groups, but we do not treat with higher homotopy groups).

To see (b), we just notice that the maps ${\bf h}_{W,W'}$
in the above commutative diagram gives us isomorphisms of their
fundamental groups (\cite{Ar}) together with the exact sequence
\[
\begin{CD}
\{0\}=\pi_1(F) @>>> \pi_1(\cal{F}\Lambda_{W}^{\,(0)}(H))
          @>>{{\bf p}_W}_*> \pi_1(\Lambda(W^{\circ}/W)) @>>> \pi_0(F)=\{0\}.
\end{CD}
\]
\end{proof}


\subsection{Proof of Theorem \ref{thm:pi_1} (b) and (c)}


In this section first we explain  
the space $\cal{U}_{\cal{F}}(H_J)$ in
the framework of the complexified symplectic Hilbert space
(Proposition \ref{prop:introduce})
and give a proof of the isomorphisms:
\begin{prop}\label{prop:homotopy}
$\pi_{1} ({\cal F}\LLa_{\lambda}(H))
\underset{(\cal{S}_{\lambda})_*}{\xrightarrow{\sim}}
\pi_1(\cal{U}_{\cal{F}}(H_J))
\xrightarrow{\sim}\Z$.
\end{prop}
{\it Then these will give a proof of} Theorem \ref{thm:pi_1} (b) and (c).

Let $H$ be a separable symplectic Hilbert space with compatible 
symplectic form $\omega$, an inner product $\inner{\bullet}{\bullet}$
and an almost complex structure  $J$,
$\omega(x,y) =\inner{J(x)}{y},\quad J^2 = -Id$.

The complexification $H\otimes\C$ 
of the real Hilbert space is installed
with the Hermitian inner product as usual and we denote by
${\LLa}^{\C}(H\otimes \C)$ the space of complex Lagrangian
subspaces in $H\otimes \C$:
\begin{equation*}
\Lambda^{\C}(H\otimes\C)=\{l\;|\;
~l ~\text{is a complex subspace such that} ~l^{\perp}=J(l)\}.
\end{equation*}
Then a subgroup of the unitary operators in $H\otimes \C$,
we denote it by
$\mathcal{U}_0(H\otimes \C)$, 
consisting of those operators $U$ that 
$U(l)^{\perp}=J(U(l))$ for any $l\in \LLa^{\C}(H\otimes \C)$
acts on ${\LLa}^{\C}(H\otimes \C)$ transitively.
This condition for $U \in\mathcal{U}_0(H\otimes \C)$ is
equivalent
to say that it commutes with the complexified almost
complex structure $J$.

Taking the complexification of $\lambda\in \Lambda(H)$ gives us a
natural embedding
$\Lambda(H)\rightarrow\Lambda^{\C}(H\otimes\C)$,
and its restriction to $\mathcal{F}\Lambda_{\lambda}(H)$ has the
image in 
$\mathcal{F}\Lambda_{\lambda\otimes \C}^{\,\C}(H\otimes \C)$, 
a subspace of 
$\Lambda^{\C}(H\otimes\C)$ consisting of those subspaces which are Fredholm
pairs with $\lambda \otimes \C$.
We denote this map by $\otimes\mathbb{C}$.

When we consider an operator $U\in \mathcal{U}(H_J)$ 
as a real operator and take its
complexification, we denote it by $U^{\C}$, then
$U^{\C}$ is in $\mathcal{U}_0(H\otimes \C)$ and 
we have $U(\mu)\otimes\C$ = ${U^{\C}}(\mu\otimes\C)$, 
$\mu\in \Lambda(H)$, and the following diagram is a commutative:

\begin{equation}
\begin{CD}
\mathcal{U}(H_J) @>{U\mapsto U^{\C}}>>\mathcal{U}(H\otimes\C)\\
@V{\rho_{\ell}}VV  @VV{\rho_{\ell\otimes\C}}V\\
\Lambda(H) @>>{\otimes\mathbb{C}}>\Lambda^{\C}(H\otimes\C).\\
\end{CD}
\end{equation}

Let $E_{\pm}=\{z\in H\otimes\C\;|\; J(z)=\pm \sqrt{-1}z\}$, then 
\begin{equation*}
H\otimes\C=E_{+}\oplus E_{-},
\end{equation*}
is an orthogonal decomposition of $H\otimes\C$.
If $U\in\cal{U}_0(H\otimes\C)$, then
$U(E_{\pm})=E_{\pm}$.
Hence we have an isomorphism
\begin{equation*}
\cal{U}_0(H\otimes\C)\cong \cal{U}(E_{+})\times \cal{U}(E_{-}),
\end{equation*}
where $\cal{U}(E_{+})$ denotes the group of unitary operators 
on $E_{+}$, and so on.
Also the space $\Lambda^{\C}(H\otimes\C)$ is identified
with the space of graphs  of unitary operators 
$U$ $\in$ $\cal{U}(E_{+},E_{-}),\,\, U:E_{+}\rightarrow E_{-}$. 

Let ${\mathfrak{K}}: H_J \rightarrow E_{+}$, $u \mapsto u\otimes 1
- J(u)\otimes\sqrt{-1}$ and 
${\mathfrak{k}}:H_J \rightarrow E_{-}$, 
$ u \mapsto u\otimes 1 + J(u)\otimes\sqrt{-1}$,
be an isomorphism and an anti-isomorphism respectively, 
then, the following diagram is commutative:
\begin{lem}
\begin{equation*}
\begin{CD}
H_J @>{\tau_{\lambda}}>>H_J\\
@V{\mathfrak{K}}VV  @VV{\mathfrak{k}}V\\
E_{+} @>>{T_{\lambda}}>E_{-},\\
\end{CD}
\end{equation*}
where $\tau_{\lambda}$ is the complex conjugation defined
through the identification $H_J \cong \lambda\otimes\C$, and
the graph of the unitary operator $T_{\lambda}$ is $\lambda\otimes\C$,
$\lambda\otimes\C = \{x+T_{\lambda}(x)\,|\,x\in E_{+}\}$.
\end{lem}

Now we have 
\begin{prop}\label{prop:introduce}
Let 
$\varPhi_{\lambda} : \cal{U}_{\cal{F}}(H_J)\rightarrow 
\cal{F}\Lambda_{\lambda\otimes\C}^{\C}(H\otimes\C)$
be a map defined by
$\varPhi_{\lambda}(V) = \text{the graph of the unitary operator}\, 
\, -{\mathfrak{k}}\circ V \circ \tau_{\lambda}\circ
{\mathfrak{K}}^{-1} \in \cal{U}(E_+,E_{-})$.
Then 
$\varPhi_{\lambda}$ is an isomorphism 
and the diagram is commutative:
\begin{equation}\label{varPhi}
\begin{psmatrix}
[name=FL] \cal{F}\Lambda_{\lambda}(H)  & [name=FLC]
\cal{F}\Lambda_{\lambda\otimes\C}^{\C}(H\otimes\C)\\
 & [name=U] \cal{U}_{\cal{F}}(H_J)
\psset{arrows=->}
\ncline{FL}{FLC}^{\otimes\mathbb{C}}\ncline{FL}{U}_{\cal{S}_{\lambda}}
\ncline{U}{FLC}>{\varPhi_{\lambda}}.
\end{psmatrix}
\end{equation}
\end{prop}

\begin{proof}
We will be enough to prove the commutativity of the diagram.
Let $U \in \cal{U}_{\lambda}(H_J)$.  Since $U^{\C}|_{E_{\pm}}$ can be
identified with $U$ through the map ${\mathfrak{K}}$ and 
${\mathfrak{k}}$
respectively, and we
have  $U^{\C}(\lambda^{\perp}\otimes\C)$ = $\{U(x)-U\circ
T_{\lambda}(x)\,|\,x\in E_{+}\}$ = $\{x-U\circ T_{\lambda}\circ
U^{-1}(x)\,|\, x\in E_{+}\}$.  By the above lemma
${\mathfrak{k}}\circ U\circ T_{\lambda}
\circ U^{-1}\circ {\mathfrak{K}}^{-1}$
= ${\mathfrak{k}}\circ U\circ \tau_{\lambda}\circ U^{-1}\circ
\tau_{\lambda}\circ \tau_{\lambda}\circ {\mathfrak{K}}^{-1}$ =
${\mathfrak{k}}\circ U\circ 
\theta_{\lambda}(U) \circ\tau_{\lambda}\circ {\mathfrak{K}}^{-1}$,
which gives the commutativity of the diagram.
\end{proof}

Let $W$ be a closed finite codimensional subspace in 
$\lambda\otimes\C$ and we denote 
by $\cal{F}\Lambda_W^{(0)}(H\otimes\C)$ 
a subspace of $\cal{F}\Lambda_{\lambda\otimes\C}^{\C}(H\otimes\C)$
consisting of those subspaces $l$ which do not intersect with $W$
except $0$.
Let $H_W$ =
$J(W^{\perp}\cap(\lambda\otimes\C))+W^{\perp}\cap(\lambda\otimes\C)$,
and $\Lambda(H_W)$ be the similar space as $\Lambda(H\otimes\C)$
(note that $H_W$ is invariant under the map $J$).  $\Lambda(H_W)$ is
identified with the space of unitary operators on 
$W^{\perp}\cap(\lambda\otimes\C)$.  Let 
${\pi}_W : 
\cal{F}\Lambda_W^{(0)}(H\otimes\C) \ni l \rightarrow 
(l\,\cap\,(J(W^{\perp}\cap\lambda\otimes\C)+\lambda\otimes\C)+W)\cap
W^{\perp} \in \Lambda(H_W)$, 
and then ${\pi}_W : \cal{F}\Lambda^{(0)}_W(H\otimes\C) 
\rightarrow \Lambda(H_W)$ is a fiber bundle 
with the contractible fiber.   A typical fiber 
= $\pi^{-1}_W(J(\lambda\otimes\C\cap W^{\perp}))$ is isomorphic to
the space $\widehat{\cal{B}}(W) \times \cal{B}(W, \lambda\otimes\C\cap
W^{\perp})$, where $\widehat{\cal{B}}(W)$ is the space of selfadjoint
operators on $W$ and $\cal{B}(W, W^{\perp}\cap(\lambda\otimes\C))$ 
is the space of bounded operators from $W$ to 
$W^{\perp}\cap(\lambda\otimes\C)$.
Unfortunately
for any pair of such subspaces $W_1$ and $W_2$ 
satisfying $W_1\subset W_2$
there are no natural map $\Lambda(H_{W_2}) \rightarrow
\Lambda(H_{W_1})$ which
makes the diagram
\begin{equation*}
\begin{CD}\cal{F}\Lambda_{W_2}^{(0)}(H\otimes\C)@>>>
\cal{F}\Lambda_{W_1}^{(0)}(H\otimes\C)\\
@VV{\pi_{W_2}}V @VV{\pi_{W_1}}V\\
\Lambda(H_{W_2}) @>>> \Lambda(H_{W_1}) 
\end{CD}
\end{equation*}
commutative. However if we define a map $s_W : \Lambda(H_W)
\rightarrow \cal{F}\Lambda_W^{(0)}(H\otimes\C)$ 
by $s_W(l) = l+J(W)$, then $\pi_W\circ s_W =Id$ and 
we have the following commutative diagram:
\begin{equation*}
\begin{CD}
\Lambda(H_{W_2}) @>>{\bf{i}}_{W_1,W_2}> \Lambda(H_{W_1}) \\
@VV{s_{W_2}}V @VV{s_{W_1}}V\\
\cal{F}\Lambda_{W_2}^{(0)} @>>> \cal{F}\Lambda_{W_1}^{(0)},
\end{CD}
\end{equation*}
where the map 
${\bf{i}}_{W_1,W_2}:\Lambda(H_{W_2}) \rightarrow \Lambda(H_{W_1})$
is defined as ${\bf{i}}_{W_1,W_2}(l) = l+J(W_2 \cap W_1^{\perp})$.

Then 
for any compact subset $K$ in
$\cal{F}\Lambda_{\lambda\otimes\C}^{\C}(H\otimes\C)$
we can find such a finite codimensional subspace $W$ in $\lambda\otimes\C$ 
that for any $l$ in $K$, $l\cap W = \{0\}$, so  
$\bigcup\cal{F}\Lambda_W^{(0)}(H\otimes\C)$ 
= $\cal{F}\Lambda_{\lambda\otimes\C}^{\C}(H\otimes\C)$.
Hence 
$\lim\limits_{W\rightarrow\{0\}}\pi_{k}(\cal{F}\Lambda_W^{(0)}(H\otimes\C))$ 
= $\pi_k(\cal{U}_{\cal{F}}(H_J))$ =
$\lim\limits_{W\rightarrow \{0\}}\pi_k(\Lambda(H_W)$. 
These show that the homotopy groups of
$\cal{U}_{\cal{F}}(H_J)$ coincides with the stable homotopy groups
of the unitary group, which together gives the proof of
Proposition \ref{prop:homotopy}, and finally gives us a proof of
Theorem \ref{thm:pi_1} (b) and (c).

\section{Maslov index in the infinite dimension}
In the last section we proved that the fundamental group
of the Fredholm-Lagrangian-Grassmannian is isomorphic to
$\mathbb{Z}$.
So in this section we define an integer, 
so called, 
the {\it Maslov index}, for arbitrary continuous paths in the
Fredholm-Lagrangian-Grassmannian $\cal{F}\Lambda_{\lambda}(H)$.
In particular it gives us an explicit isomorphism between the
fundamental group of the Fredholm-Lagrangian-Grassmannian 
and $\mathbb{Z}$.
We base on a spectral property of the Fredholm operator
to define the {\it Maslov index}, so that our method 
is valid for both of finite and infinite dimensional cases. 

\subsection{Maslov index for continuous paths}

Let
\begin{equation*}
\begin{array}{rcc}
{\bf d}:\;I=[0,1]&\rightarrow &\cal{U}_{\cal{F}}(H_J)\\
t&\mapsto & {\bf d}(t)
\end{array}
\end{equation*}
be a continuous path in $\cal{U}_{\cal{F}}(H_J)$.
First we prove
\begin{lem}\label{lem:partition and pnumber}
There exist a partition $0=t_0<t_1<\dots<t_N=1$ of the interval $I$
and positive numbers 
$\varepsilon_j\,(j=1,\dots,N)$ with $0<\varepsilon_j<\pi$
such that
\begin{equation}\label{eq:partition}
e^{\sqrt{-1}(\pi \pm \varepsilon_j)}\in \rho ({\bf d}(t))
\end{equation}
and
\begin{equation}\label{eq:pnumber}
\sum_{|\theta|\le\varepsilon_j}\dim
\Ker({\bf d}(t)-e^{i(\pi+\theta)})\,<\, \infty
\end{equation}
for $t_{j-1}\le t\le t_j$.
\end{lem}

Note here $\rho({\bf d}(t))$ denote the resolvent set of the operator
${\bf d}(t)$.

\begin{proof}
Since ${\bf d}(t)+\Id $ is a Fredholm operator, for each  
$t\in I=[0,1]$ we can find a positive number 
$\varepsilon _{t} > 0$ such that
\[
\{e^{\sqrt{-1}(\pi +\theta )}|\ 0 < |\theta | \le \varepsilon_{t}\}
\subset 
\rho({\bf d}(t)),
\]
because $-1$ is an isolated eigenvalue of ${\bf d}(t)$ 
with finite multiplicity.
So there exist positive numbers $\delta^{\pm}_t >0$
such that the projection operator $P_s$ defined by
\begin{equation}\label{eq:eigen-projection}
P_s=\frac{1}{2\pi\sqrt{-1}}
\int_{|u+1|=\varepsilon _{t}} (u-{\bf d}(s))^{-1}du
\end{equation}
has the constant rank equal to $\dim\Ker ({\bf d}(t)+\Id)$
for $s\in [t-\delta_{t}^-,t+\delta_{t}^+]$,
because $\{P_s\}$ is a norm continuous family :
\begin{equation}\label{eq:const rank}
\dim P_s(H_J) = \dim\Ker ({\bf d}(t)+\Id) \,,\,
s\in [t-\delta_{t}^-,t+\delta_{t}^+]. 
\end{equation}

Note that
the continuity of the family of projection operators $\{P_s\}$
is proved by using the ``resolvent equation'': 
$$(u-{\bf d}(s))^{-1} -(u- {\bf d}(t))^{-1} =
(u-{\bf d}(s))^{-1}({\bf d}(s)-{\bf d}(t))(u-{\bf d}(t))^{-1}
$$
The continuity of $\{{\bf d}(t)\}$
is reflected by this equation to the continuity of the family
$\{P_s\}$.

Hence we have an open covering 
$\{(t-\delta_{t}^-,t+\delta_{t}^+)\}_{t \in I}$ 
of the interval $I$ and positive numbers 
$\{\varepsilon _{t}\}_{t \in I}$ 
such that for $s \in [ t-\delta_{t}^-,t+\delta_{t}^+]$
\[
\sum _{|\theta | \le \varepsilon_{t}} 
\dim \Ker ({\bf d}(s)-e^{\sqrt{-1}(\pi +\theta )})
=\dim \Ker({\bf d}(t)+\Id)
\]
\[
e^{\sqrt{-1}(\pi \pm \varepsilon _{t})} \in \rho({\bf d}(s)).
\]
Now we can choose enough number of points $\{s_i\}_{i=0}^{N-1}$ 
satisfying following properties:
\[
0=s_0 <s_1 < \cdots < s_{N-1}=1 ~\text{such that} 
\]
\[
   s_{i-1} <   s_i-\delta_{s_i}^-,
\]
\[
s_{i-1}+\delta_{s_{i-1}}^+ < s_{i},
\]
\[ 
s_{i}-\delta_{s_{i}}^- < s_{i-1}+\delta_{s_{i-1}}^+.
\]
Here if necessary, we need to replace $\delta_{s_i}^{\pm}$
by a smaller one (but then the number of the points $\{s_i\}$ will
increase).  
Then finally we define the point $t_k$ ($k=0,\dots,N)$ in such a way
that 
\[
t_0 = s_0 =0, t_1 = \delta_0^+, t_2 = s_1+\delta_{s_1}^+,
t_3= s_2+\delta_{s_2}^+, ~\cdots,
\]
\[
t_{N-1} = s_{N-2}+\delta_{s_{N-2}}^+, t_N=s_{N-1}=1,
\]
and on the each interval $[t_{k-1},t_k]$ 
we take  the positive number
$\varepsilon_k =\varepsilon_{s_{k}}$, 
then we have a desired partition of the interval $I$ and
positive numbers satisfying 
(\ref{eq:partition}) and (\ref{eq:pnumber}).
\end{proof}

We now define a quantity, we call it  ``{\it unitary Maslov index}'',
and denote by $\mathbf{M}(\{{\bf d}(t)\})$ 
of a continuous curve 
$\{{\bf d}(t)\}_{t\in I}\subset \cal{F}\Lambda_{\lla}(H)$. 

\begin{defn}\label{def:unitary Maslov index}
Let $\{t_j\}_{j=0}^N$ be the partition of the interval $I$
and $\{\varepsilon\}_{j=1}^N$ positive numbers satisfying
(\ref{eq:partition}) and (\ref{eq:pnumber}) 
as in the above lemma, then we define
\begin{equation}\label{eq:unitary Maslov index}
\mathbf{M}(\{{\bf d}(t)\})=\sum_{j=1}^N (k(t_j,\varepsilon_j)
-k(t_{j-1},\varepsilon_{j}))
\end{equation}
with
\begin{equation*}
k(t,\epsilon_j)=\sum_{0\le\theta\le\varepsilon_j}\dim
\Ker({\bf d}(t)-e^{i(\pi+\theta)})
\end{equation*}
for $t_{j-1}\le t\le t_j$.
\end{defn}

In order that the definition has a meaning, we need to prove  

\begin{prop}\label{prop:invariance}
The definition of the quantity $\mathbf{M}(\{{\bf d}(t)\})$
does not depend on the choices of the 
partition $\{t_j\}_{j=0}^N$ of the interval $I$
nor on the positive numbers $\{\varepsilon_j\}_{j=1}^{N}$
satisfying (\ref{eq:partition}) and (\ref{eq:pnumber}).
\end{prop}

This follows from the following lemma:
Let $\{t_j\}_{j=0}^N$ be the partition of the interval $I$ and 
$\{\tilde{\varepsilon}_{j}\}_{j=1}^N$ another positive numbers
satisfying (\ref{eq:partition}) and (\ref{eq:pnumber}). 
\begin{lem}\label{lem:invariance}
The two integers coincide defined in 
(\ref{eq:unitary Maslov index})
one in terms of the partition $\{t_j\}_{j=0}^N$ and positive numbers
$\{{\varepsilon}_{j}\}_{j=1}^N$ and other
in terms of the ``same partition'' $\{t_j\}_{j=0}^N$ and
``different positive numbers'' $\{\tilde{\varepsilon}_{j}\}_{j=1}^N$.
\end{lem}

\begin{proof}
Since both of
$e^{\sqrt{-1}(\pi+\varepsilon_j)}$ and
$e^{\sqrt{-1}(\pi+\tilde{\varepsilon}_j)} \in \rho({\bf c}(t))$ 
on each small interval $[t_{j-1},t_j]$, the difference of the
dimensions $k(t,\varepsilon_j)-k(t,\tilde{\varepsilon}_j)$
is constant on the interval $[t_{j-1},t_j]$.  
Hence we have
$$
 k(t_j,\varepsilon_j)-k(t_{j-1},\varepsilon_j)
=k(t_j,\tilde{\varepsilon}_j)-k(t_{j-1},\tilde{\varepsilon}_j),
$$
which proves the lemma.
\end{proof}

\begin{proof}[Proof of Proposition of \ref{prop:invariance}]
By adding a suitable number of points both in the partitions 
$\{t_j\}$ and $\{\tilde{t}_l\}$,
we may  assume that $t_{j-1} < \tilde{t}_j < t_j$ for each $j$.
Then from Lemma \ref{lem:invariance} we have
\begin{align}
&k(t_j,\varepsilon_j)-k(t_{j-1},\varepsilon_j)\label{eq:diff1}\\
&=k(t_j,\varepsilon_j)- k(\tilde{t}_j, \varepsilon_j)
+k(\tilde{t}_j, \varepsilon_j)- k(t_{j-1},\varepsilon_j)\\
&=k(t_j,\tilde{\varepsilon}_{j+1})- 
k(\tilde{t}_j,\tilde{\varepsilon}_{j+1})
+k(\tilde{t}_j,\tilde{\varepsilon}_{j})
-k(t_{j-1},\tilde{\varepsilon}_{j})\label{eq:diff2},
\end{align}
which gives us the coincidence of the two integers by adding
(\ref{eq:diff1}) and (\ref{eq:diff2}) with respect to $j$. 
\end{proof}

\begin{notation}
\begin{enumerate}
\item
Let $\{{\bf d}_1(t)\}_{t\in [0,1]}$ and 
$\{{\bf d}_2(t)\}_{t\in [0,1]}$ 
be continuous curves with the relation ${\bf d}_1(1)={\bf d}_2(0)$,
then we denote the catenation of these two curves by 
${\bf d}_1*{\bf d}_2$:
\begin{equation*}
({\bf d}_1* {\bf d}_2)(t)=
\left\{\begin{array}{ll}
{\bf d}_1(t/2) & ~\text{for} ~0\le t \le 1/2,\\
{\bf d}_2(2t-1) & ~\text{for} ~1/2 \le t \le 1.
\end{array}\right.
\end{equation*}
\item The curve
$-{\bf d}$ denotes 
the curve defined by $-{\bf d}(t)={\bf d}(1-t),~t\in I.$
\end{enumerate}
\end{notation}

This unitary Maslov index has the following properties:
\begin{thm}\label{thm:properties of unitary Maslov index}
\begin{enumerate}
\item 
Additivity under the catenation of the paths, and
\item
Modulo sign and additive constants, 
it is only a homotopy invariant of curves
in $\cal{U}_{\cal{F}}(H_J)$ with fixed endpoints and distinguishes
the homotopy classes.
\end{enumerate}
\end{thm}

\begin{proof}
(a) The additivity follows from the very definition of the 
quantity $\mathbf{M}\{{\bf d}(t)\}$.

(b) Let $\{{\bf w}(s,t)\}_{(s,t)\in I\times I}\subset\cal{U}_{\cal{F}}(H_J)$
be a continuous 
two-parameter family. By the similar continuity arguments
for the projection operator (\ref{eq:eigen-projection}) in the proof
of Lemma \ref{lem:partition and pnumber}, for each $s\in I$ 
there are a positive number $c_s > 0$, the
partition $\{t_j\}$
of the interval and the positive numbers $\{\varepsilon_j\}$
such that (\ref{eq:partition}) and (\ref{eq:pnumber}) hold
for $t_{j-1}\le t\le t_j$ and $ |s'-s|\le c_s$:
\begin{equation*}
e^{\sqrt{-1}(\pi \pm \varepsilon_j)}\in \rho({\bf w}(s',t))
\end{equation*}
and
\begin{equation*}
\sum_{|\theta|\le\varepsilon_j}\dim
\Ker({\bf w}(s',t)-e^{i(\pi+\theta)})\,<\, \infty.
\end{equation*}
So on the each small rectangle $[t_{j-1},t_j]\times [s,s+c_s]$, 
$v\in [s,s+c_s]$ 
\begin{align*}
&\sum_{0\le \theta\le\varepsilon_j}
            \dim\Ker({\bf w}(s+v,t_j)-e^{i(\pi+\theta)})\\
&\qquad\qquad\qquad-\sum_{0\le \theta\le\varepsilon_j}
             \dim\Ker({\bf w}(s+v,t_{j-1})-e^{i(\pi+\theta)})\\
&\qquad+\sum_{0\le\theta\le\varepsilon_j}
             \dim\Ker({\bf w}(s+v,t_{j-1})-e^{i(\pi+\theta)})\\
&\qquad\qquad\qquad\quad -\sum_{0\le\theta\le\varepsilon_j}
             \dim\Ker({\bf w}(s,t_{j-1})-e^{i(\pi+\theta)})\\
=&\sum_{0\le\theta\le\varepsilon_j}
             \dim\Ker({\bf w}(s+v,t_j)-e^{i(\pi+\theta)})\\
&\qquad\qquad\qquad -\sum_{0\le
  \theta\le\varepsilon_j}\dim\Ker({\bf w}(s,t_j)-e^{i(\pi+\theta)})\\
&\qquad+\sum_{0\le\theta\le\varepsilon_j}
             \dim\Ker({\bf w}(s,t_j)-e^{i(\pi+\theta)})\\
&\qquad\qquad\qquad\quad -\sum_{0\le \theta\le\varepsilon_j}
             \dim\Ker({\bf w}(s,t_{j-1})-e^{i(\pi+\theta)}).
\end{align*}
By adding above equalities with respect to $j$
we have in general (locally with respect to the parameter $s$)
\begin{align*}
&\mathbf{M}(\{{\bf w}(s+v,0)\}_{0\le v\le c_s}) 
+\mathbf{M}(\{{\bf w}(s+c_s,t)\}_{t\in I})\\
=&\mathbf{M}(\{{\bf w}(s,t)\}_{t\in I})
+\mathbf{M}(\{{\bf w}(s+v,1)\}_{0\le v \le c_s}),
\end{align*}
and then on the rectangle $I\times I$
\begin{align*}
&\mathbf{M}(\{{\bf w}(s,0)\}_{s\in I}) 
+\mathbf{M}(\{{\bf w}(1,t)\}_{t\in I})\\
=&\mathbf{M}(\{{\bf w}(0,t)\}_{t\in I})
+\mathbf{M}(\{{\bf w}(s,1)\}_{s\in I}).
\end{align*}

Now here we assume that 
${\bf w}(s,0) \equiv {\bf w}(0,0)$ and
 ${\bf w}(s,1)\equiv {\bf w}(0,1)$ ($s\in I$),    
hence 
$\mathbf{M}(\{{\bf w}(0,t)\}_{t\in I})
=\mathbf{M}(\{{\bf w}(1,t)\}_{t\in I})$,
and this shows the homotopy invariance of the integer
$\mathbf{M}(\{{\bf w}(t)\})$.  

The uniqueness (mod additive constant and signature) 
follows from the fact that $\pi_1(\cal{U}_{\cal{F}}(H))
\cong \mathbb{Z}$.
\end{proof}

The space $U_{\cal F}(H_J)$ is closed under the adjoint
operation, so we have    
\begin{prop} 
$\mathbf{M}(\{{\bf w}(t)\})$ = $ - \mathbf{M}(\{{\bf w}(t)^*\})$.   
\end{prop}


Using this ``unitary Maslov index'' 
we give a functional analytic
definition of the infinite version of the 
{\it Maslov index} for arbitrary continuous paths in the
Fredholm-Lagrangian-Grassmannian.  

Let $\mu:\;I\rightarrow\cal{F}\Lambda_{\lambda}(H)$ be a continuous path
in $\cal{F}\Lambda_{\lambda}(H)$ (so that $\cal{S}_{\lambda}\circ\mu$ is a
continuous path in $\cal{U}_{\cal{F}}(H_J)$).

\begin{defn}\label{def:M-index-path}
We define the Maslov
index of the curve $\{\mu(t)\}$ with respect to $\lambda$ by
\begin{equation*}
\mathbf{Mas}(\{\mu(t)\},\lambda)
=\mathbf{M}(\{\cal{S}_{\lambda}(\mu(t))\}).
\end{equation*}
\end{defn}

By Theorem \ref{thm:pi_1}, the Maslov index inherits the all
properties of the ``unitary Maslov index''.

In the case that 
$\cal{F}\Lambda_{\lambda }(H)=\cal{F}\Lambda_\mu(H)$, 
but $\lambda \neq \mu$, then Maslov cycles 
$\mathfrak{M}_{\lla}(H)$ and $\mathfrak{M}_{\mu}(H)$ do not
coincide. Hence Maslov indexes for a path with respect to 
$\mathfrak{M}_{\lla}(H)$ and $\mathfrak{M}_{\mu}(H)$ 
will not coincide in general. However
for loops, as in the finite dimensional case we have 
\begin{prop}\label{prop:maslov-coin}
Let $\lambda ,\mu  \in \Lambda(H)$ and assume that $\mu=U_1(\lla)$ with
a unitary operator $U_1\in\mathcal{U}_{res}(H_J)$.
Then for any continuous loops $\{{\bf c}(t)\}_{t\in [0,1]}$ in
$\cal{F}\LLa_{\lla}(H) =\cal{F}\LLa_{\mu}(H)$ 
their Maslov indexes coincide: 
\begin{equation*}
\mathbf{Mas}(\{{\bf c}(t)\},\lambda)
=\mathbf{Mas}(\{{\bf c}(t)\},\mu).
\end{equation*}
\end{prop}

\begin{proof}
Let $\{U_s\}_{s\in [0,1]}$ be a continuous curve in 
$\mathcal{U}_{res}(H_J)$
which joins $\lla$ and $\mu$, that is, $U_0 =\Id$ and $U_1(\lla)=\mu$.
Note then for each $s\in[0,1]$, $\cal{F}\LLa_{U_s(\lla)}(H)$ =
$\cal{F}\LLa_{\lla}(H)$.

We define a map ${\bf h}:I\times I \to \cal{U}_\cal{F}(H_J)$:
\begin{equation*}
{\bf h}(s,t)=\left\{\begin{array}{ll}
\cal{S}_{U_{2ts}(\lla)}({\bf c}(t)) 
& ~\text{for}~(s,t)\in [0,1]\times [0,1/2]\\
\cal{S}_{U_{(2-2t)s}(\lla)}({\bf c}(t)) & 
~\text{for} ~(s,t)\in [0,1]\times [1/2,1].\\
\end{array}\right.
\end{equation*}
Then $\{{\bf h}(s,t)\}$ is a homotopy between the loop
$\{\mathcal{S}_{\lambda}({\bf c}(t))\}$ 
and the loop $\{{\bf h}(1,t)\}$ with
the fixed common initial and end point 
$\mathcal{S}_{\lambda}({\bf c}(0))=
\mathcal{S}_{\lambda}({\bf c}(1))
= {\bf h}(s,0)={\bf h}(s,1)$, $s\in [0,1]$.
Hence 
\begin{equation*}
\mathbf{Mas}(\{{\bf c}(t)\},\lambda)
=\mathbf{M}(\{\cal{S}_{\lambda}({\bf c}(t))\})
=\mathbf{M}(\{{\bf h}(1,t)\}).
\end{equation*}

By the same way for the loops $\{ {\bf h}(1,t)\}$ and
$\mathcal{S}_{\mu}({\bf c}(t))$ we can construct a homotopy
in $\mathcal{U}_{\mathcal{F}}(H_J)$ between them
and these shows the coincidence of the two Maslov indexes.


\end{proof}

\begin{cor}
Let $\{{\bf c}(t)\}_{t\in [0,1]}$ be a continuous path
in $\mathcal{F}\Lambda_{\lambda}(H)$ such that 
${\bf c}(0),~{\bf c}(1)\notin\mathfrak{M}_{\lambda}$ and let
$\{U_s\}_{s\in[0,1]}\subset \mathcal{U}_{res}(H_J)$ 
be a continuous family with $U_0= Id$. We assume that
${\bf c}(0), {\bf c}(1)\notin\mathfrak{M}_{U_s(\lambda)}(H)$
for all $s\in[0,1]$, then for all $s$
$$
\mathbf{Mas}(\{{\bf c}(t)\},\lambda)
=\mathbf{Mas}(\{{\bf c}(t)\},U_s(\lambda)).
$$
\end{cor}

\subsection{H\"ormander index in the infinite dimension}

Let $\lambda ,\mu \in \Lambda (H)$ and assume that 
$\mu =U(\lla)$ with a unitary operator $U$ of the form 
$\Id + ~\text{compact operator}$,
and let $\{{\bf c}(t)\}_{t\in[0,]}$ be a continuous curve in 
$\cal{F}\Lambda _\lambda (H)=\cal{F}\Lambda_\mu (H)$.

\begin{prop}
The difference
\[
\mathbf{Mas}(\{{\bf c}(t)\},\lambda)-\mathbf{Mas}(\{{\bf c}(t)\},\mu)
\]
depends only on the end points.
\end{prop}

\begin{proof}
Let $\{\tilde{{\bf c}}(t)\}$ be another path 
with $\tilde{{\bf c}}(0)={\bf c}(0),\tilde{{\bf c}}(1)={\bf c}(1)$,
then
\begin{eqnarray*}
 \mathbf{Mas}(\{{\bf c} \ast (-\tilde{{\bf c}})(t)\},\lambda)
&=&\mathbf{Mas}(\{{\bf c}(t)\},\lambda)
-\mathbf{Mas}(\{\tilde{{\bf c}}(t)\},\mu)\\
& =& \mathbf{Mas}(\{{\bf c}\ast (-\tilde{{\bf c}})(t)\},\mu)\\
&=&\mathbf{Mas}(\{{\bf c}(t)\},\mu)
-\mathbf{Mas}(\{\tilde{{\bf c}}(t)\},\mu)
\end{eqnarray*}
by Proposition \ref{prop:maslov-coin}.
Hence we have the desired result.
\end{proof}

Using of this fact
we can define an infinite dimensional version of an 
integer, called {\it H\"ormander index}, for  four Lagrangian
subspaces:
\begin{defn}
For $\mu=U(\lla)\in \Lambda (H)$ 
with $U = \Id + ~\text{compact operator}$, 
then we call the difference 
\[
\mathbf{Mas}(\{{\bf c}(t)\},\lambda)-\mathbf{Mas}(\{{\bf c}(t)\},\mu)
\]
the {\it H\"ormander index} in the infinite dimension and denote it by 
\begin{equation}
\sigma({\bf c}(0),{\bf c}(1);\lambda ,\mu).
\end{equation}
\end{defn}

Let $\mu =U(\lla)$ be as above and let 
$\ell_0, \ell_1, \ell_3 \in \cal{F}\LLa_{\lla}(H)=\cal{F}\LLa_{\mu}(H)$,
then the H\"ormander index $\sigma(\ell_0,\ell_1 ; \lambda ,\mu)$ 
has the following properties:
\begin{prop}
\begin{equation}
\sigma(\ell_0,\ell_1; \lla,\mu) = -\sigma(\ell_1,\ell_0; \lla,\mu),
\end{equation}
\begin{equation}
\sigma(\ell_0,\ell_1; \lla,\mu) = -\sigma(\ell_0,\ell_1; \mu,\lla),
\end{equation}
\begin{equation}\label{eq:cocycle1}
\sigma(\ell_0,\ell_1; \lla,\mu)+\sigma(\ell_1,\ell_2; \lla,\mu)
=\sigma(\ell_0,\ell_2; \lla,\mu).
\end{equation}
Let $\nu = W(\lla)$ also  with a unitary operator
$W= \Id + ~\text{copmact operator}$, 
then the cocycle conditions with respect to the last two components hold:
\begin{equation}\label{eq:cocycle2}
\sigma(\ell_0,\ell_1; \lla,\mu) +\sigma(\ell_0,\ell_1; \mu,\nu) =
\sigma(\ell_0,\ell_1; \lla,\nu).
\end{equation}
If we assume moreover $\ell_1 =V(\ell_0)$ with a
unitary operator $V$ of the form  $\Id + \,\text{compact operator}$, then
\begin{equation}
 \sigma(\ell_0,\ell_1; \lla,\mu) = -\sigma(\lla,\mu; \ell_0,\ell_1).
\end{equation}
\end{prop}
\begin{proof}
Four properties except last one
follow directly {}from the definition itself.

Let $\{U_s\}_{s\in[0,1]}$ and  $\{V_t\}_{t\in[0,1]}$ be 
such curves of unitary operators that each operator 
$U_s$ and $V_t$ are of the form of $\Id + ~\text{compact operator}$
and assume $U_0=\Id, U_1(\lla)=\mu, V_0=\Id$ and $V_1(\ell_0)=\ell_1$. 
Then for any $s\in[0,1]$,  ${\cal F}\LLa_{U_s(\lla)}(H) 
={\cal F}\LLa_{\lla}(H)$,
and for any $(s,t)$ $(U_s(\lla),V_t(\ell_0))$
is a Fredholm pair. So the two-parameter continuous family
of unitary operators 
$\{\cal{S}_{U_s(\lla)}(V_t(\ell_0))\}$ are in ${\cal U}_{{\cal F}}(H_J)$.  
Let us define a curve $\{{\bf c}(t)\}_{0\le t\le 4}$:
\begin{equation*}
{\bf c}(t)=\left\{\begin{array}{ll}
\cal{S}_{\lla}(V_t(\ell_0)) & ~\text{for}  ~0 ~\leq ~t ~\leq ~1,\\
\cal{S}_{U_{t-1}(\lla)}(\ell_1) & ~\text{for}  ~1 ~\leq ~t ~\leq ~2,\\
\cal{S}_{\mu}(V_{3-t}(\ell_0)) &~\text{for}  ~2 ~\leq ~t ~\leq ~3,\\
\cal{S}_{U_{4-t}(\lla)}(\ell_0)  & ~\text{for}  ~3 ~\leq ~t ~\leq ~4.
\end{array}\right.
\end{equation*}
The unitary Maslov index $\mathbf{M}(\{{\bf c}(t)\}_{0\le t\le 4})$ of this curve 
is zero, so
$$
\mathbf{M}(\{\cal{S}_{\lla}(V_t(\ell_0))\}_{t\in [0,1]})
-\mathbf{M}(\{\cal{S}_{\mu}(V_t(\ell_0))\}_{t\in [0,1]})
$$
$$=\mathbf{M}(\{\cal{S}_{U_t(\lla)}(\ell_0)\}_{t\in [0,1]})
-\mathbf{M}(\{\cal{S}_{U_t(\lla)}(\ell_1)\}_{t\in [0,1]}),
$$ and by Proposition \ref{prop:symmetry of Souriau map} 
this equal to
$$
=-\mathbf{M}(\{\cal{S}_{\ell_0}(U_t(\lla))\}_{t\in [0,1]})
+\mathbf{M}(\{(\cal{S}_{\ell_1}(U_t(\lla))\}_{t\in [0,1]}).
$$
Hence 
\begin{equation*}
 \sigma(\ell_0,\ell_1; \lla,\mu) = -\sigma(\lla,\mu; \ell_0,\ell_1).
\end{equation*}
\end{proof}

\begin{rem}
The H\"ormander index was first introduced in the paper \cite{Ho2}
to describe the phase transitions in the oscillatory integral
representation
of Fourier integral operators
or Lagrangian distributions for the global formulation
in terms of, so called, Maslov line bundle.  It was given 
also as a {\v Cech}
cocycle. Our definition above is given
in terms of the Maslov index, and we need not to assume the
transversality conditions between the first two Lagrangian subspaces and
the last two Lagrangian subspaces. The reason is, of course, the
Maslov index is defined for not only loops but also 
any paths.  In earlier papers written before 
the papers \cite{Go1}, \cite{RS} and
\cite{BF1} it was only considered for loops or with the assumption
that
the end points of paths do not meet with a particularly fixed
Maslov cycle. 
However
in order to construct the Maslov line bundle, it is enough
to consider the indexes for four Lagrangian subspaces satisfying
transversality conditions.

In the next subsection we construct 
an infinite dimensional analogue of the Maslov line bundle
which will be turn out to be a kind of 
the universal Maslov line bundle.
\end{rem}

\subsection{Universal covering space of the Fredholm-Lagrangian-Grassmannian}

In this section we characterize the universal covering space
$\tilde{{\cal F}\Lambda_\ell}(H)$ of
the Fredholm-Lagrangian-Grassmannian ${\cal F}\Lambda_{\ell}(H)$
in terms of the H\"ormander
index. We show the transition functions of the principal 
bundle $\pi :\tilde{{\cal F}\Lambda_\ell}(H)\to
{\cal F}\Lambda_\ell(H)$ 
are given by the H\"ormansder index.
Here we understand the space $\tilde{{\cal F}\Lambda_\ell}(H)$ 
consisting of pairs $([c], {\bf c}(1))$ 
of homotopy classes $[c]$ of paths $\{{\bf c}(t)\}$
in ${\cal F}\Lambda_\ell(H)$ with the 
common initial point ${\bf c}(0) = \ell^{\perp}$ and its end point ${\bf c}(1)$.

Let $\lla \in \Lambda(H)$ and assume $\lla\sim\ell$,
and we define  a map
$$
\phi_{\lla}:{\cal F}\Lambda_\ell(H)\times \mathbb{Z}
\longrightarrow \tilde{{\cal F}\Lambda _{\ell}}(H) 
$$
by
$$
\phi_\lla:(\theta,n)\longmapsto [{\bf c}(t)],
$$
where $\{{\bf c}(t)\}$ is a path joining ${\ell}^{\perp}$ and $\theta$,
and $\Mas(\{{\bf c}(t)\},\ell)= n$. 
Note that we know the homotopy class of such 
paths is uniquely determined.

By the definition of the topology on
the space $\tilde{{\cal F}\Lambda_\ell}(H)$,
it is immediate to show that the map is bijective, and not
continuous on the whole space of definition, but
\begin{prop}
The map $\phi_{\lla}$ restricted to the open subset 
$$({\cal F}\Lambda_{\ell}(H) \setminus \mathfrak{M}_{\lla}(H))
\times\mathbb{Z}
 = {\bf O}_{\lla}\times\mathbb{Z}
$$
is an isomorphism with the space
$$
{\pi}^{-1} ({\cal F}\Lambda_{\ell}(H)\setminus \mathfrak{M}_{\lla}(H)).
$$
\end{prop}

Now let $\lla \sim\ell$, $\mu\sim\ell$ and
let $\nu \in {\bf O}_{\lla}\cap{\bf O}_{\mu}$. Then
if $\phi_{\lla}(\nu,n) =\phi_{\mu}(\nu,m)$, then
$n-m =\sigma(\ell^{\perp},\nu;\lla,\mu)$ and so by the
cocycle condition (\ref{eq:cocycle2}) we have at once
\begin{prop}
The maps 
\begin{equation}\label{eq:transition function}
\begin{array}{r@{\,}l}
g_{\lla,\mu}:{\bf O}_{\lla}&\cap{\bf O}_{\mu} \to \mathbb{Z}\\
 \nu &\mapsto \sigma(\ell^{\perp},\nu;\lla,\mu)\\
\end{array}
\end{equation}
are the transition functions of the principal bundle
$\pi:\tilde{{\cal F}\Lambda_\ell}(H)\to{\cal F}\Lambda_\ell(H)$ 
with
the structure group $\pi_1({\cal F}\Lambda_\ell(H))\cong\mathbb{Z}$.
\end{prop}
                
{}From this fact we can define

\begin{defn}\label{defn:maslov line bundle}
We call the complex line bundle ${\cal L}_{\ell}$ on
${\cal F}\Lambda_\ell(H)$ defined by the
transition functions $\{h_{\lla,\mu}\}$ ($\lla,\mu \sim\ell$)
$$
h_{\lla,\mu}(\nu)= e^{\sqrt{-1}\frac{\pi}{2}\sigma(\ell^{\perp},\nu;\lla,\mu)}
$$
the universal Maslov line bundle.
\end{defn}

In fact, we have  the following:
let $H=H_0+H_1$ be an orthogonal direct sum by symplectic
subspaces with $\dim H_0 =2n <+\infty$, and we fix Lagrangian
subspaces $\ell_0\in H_0$ and $\ell_1$ in $H_1$.  Then we have an
embedding $ {\bf i}:\LLa(H_0) \to \LLa(H)$ 
\begin{equation}\label{eq:embedding of Lag}
\begin{array}{r@{\,}l}
{\bf i} :\LLa(H_0) \to  & {\cal F}\LLa_{\ell_0\oplus\ell_1}(H)\\
{\bf i} :\theta \mapsto & \theta\oplus\ell^{\perp}.\\
\end{array}
\end{equation}
Then the map ${\bf i}$ gives a relation between
the H\"ormander indexes on $\LLa(H_0)$ and 
${\cal F}\Lambda_{(\ell_0\oplus\ell_1)}(H)$:
for $\lla,\mu \in \LLa(H_0)$,  
$$
\sigma(\ell_0^{~\perp},\theta;\lla,\mu)
=\sigma((\ell_0\oplus\ell_1)^{\perp},\theta\oplus\ell_1^{\perp};
\lla\oplus\ell_1,\mu\oplus\ell_1).
$$
Hence
\begin{prop}
${\bf i}^*({\cal L}_{\ell_0\oplus\ell_1}) \cong ~\text
{the Maslov line
bundle on}~\LLa(H_0).$
\end{prop}

\begin{rem}
The collections of the vector spaces
$$
\coprod\limits_{\mu\in\mathcal{F}\Lambda_{\lambda}(H)}\,\, \lambda\cap\mu
$$
and
$$
\coprod\limits_{\mu\in\mathcal{F}\Lambda_{\lambda}(H)}\,\, 
H/(\lambda+\mu)
$$
are not apparently vector bundles, but
$$
\coprod\limits_{\mu\in\mathcal{F}\Lambda_{\lambda}(H)}
\left(\stackrel{\text{max}}\bigwedge \lambda\cap\mu \right)^*\otimes
\left(\stackrel{\text{max}}\bigwedge H/(\lambda+\mu)\right)
$$
has a line bundle structure. Here 
$\stackrel{\text{max}}\bigwedge \lambda\cap\mu$ 
means the highest degree exterior product. This
is isomorphic with the induced
bundle of the Quillen determinant line bundle
on the space of all Fredholm operators 
by the map $\mu\mapsto \mathcal{P}_{\lambda}+\mathcal{P}_{\mu}$
and also its complexification is isomorphic with the induced 
bundle by the map
$Id +\mathcal{S}_{\lambda}$ (\cite {Fu}). 
This is a trivial bundle, and the Maslov line bundle
is also trivial.
\end{rem}


\subsection{Bilinear forms and Maslov index} 

For ``differentiable curves'' in ${\cal U}_{\cal F}(H_J)$
satisfying a certain non-degeneracy condition, 
there is another way of describing
the ``unitary Maslov index''. 
We define a symmetric bilinear form which is analogous to 
Duistermaat \cite{Du} and Robbin-Salamon \cite{RS}.

Let $\{{\bf c}(t)\}$ be a ``$C^1$-path'' in $\cal{U}_{\cal F}(H_J)$. 
Here we mean $C^1$-path in the following sense :
there is a continuous family $\{D_t\}_{t\in I}$ of 
bounded operators $D_t\in {\cal B}(H_J)$ satisfying
\begin{equation}\label{eq:derivative}
\norm{\frac{1}{\delta}({\bf c}(t+\delta)-{\bf c}(t))- {t\cdot D_t}} = o(1)
\end{equation}
on the interval $I$.  We denote $D_t =\frac{d}{dt}{\bf c}(t)=\dot{\bf c}(t)$.

\begin{defn}\label{def:crossing form}
\begin{enumerate}
\item A parameter $t^*$ with $0\leq t^* \leq1$ is called a {\it crossing} for the
family $\{{\bf c}(t)\}$, if $\Ker({\bf c}(t^*)+\Id)\ne\{0\}$.
\item We define the {\it crossing form} at a crossing $t^*$ as a symmetric
  bilinear form on $\Ker({\bf c}(t^*)+\Id)$ by
\begin{equation*}
\tilde{Q}_{\frak{M}}(x,y)=\left.\frac{d}{dt}(x,~R_t(y))_J\right|_{t=t^*}
\quad\text{for $x,y\in\Ker({\bf c}(t^*)+\Id)$,}
\end{equation*}
where $\{R_t\}$ is a family of bounded selfadjoint operators
given by the relation 
${\bf c}(t)={\bf c}(t^*)e^{\sqrt{-1}R_t},~R_{t^*}=0$, i.e.,
$$
R_t= -\sqrt{-1}{\Log}~({\bf c}(t^*)^{-1}\circ {\bf c}(t)) 
~(\text{for ``$\Log$'' see Remark 
\ref{Log} below}).
$$  
Then $\dot{W}(t^*)= \sqrt{-1}{\bf c}(t^*)\circ\dot{R}_{t^*}$
\item We call a crossing $t^*$ is {\it regular}, 
if the form $\tilde{Q}_{\frak{M}}$
is non-degenerate on $\Ker({\bf c}(t^*)+\Id)$.
\end{enumerate}
\end{defn}

\begin{rem}
The logarithm above is defined by the integral 
\begin{equation}\label{Log}
{\Log}~M 
= \int\limits_{-\infty}^{0}\left\{(u - M)^{-1}-(u - 1)^{-1}\Id\right\}du
\end{equation}
for a bounded invertible operator $M\in {\cal B}(H_J)$ 
whose spectrum $\sigma(M)$ does not
contain any negative real numbers: $\sigma(M)\cap (-\infty,0]= \phi$.

The integral converges in the operator norm and 
the resulting family is again $C^1$-class, if
$\{M(t)\}$ is so. The derivative in the sense of (\ref{eq:derivative})
is given by the integral :
\begin{align}\label{eq:rep derivative}
&\frac{d}{dt}{\Log}~M(t)\\ 
&= \int\limits_{-\infty}^{0}
\left\{(u-M(t))^{-1}\circ \frac{d}{dt}M(t)\circ (u-1)^{-1}\right\}du.\notag
\end{align}

For our case $M(t) = {\bf c}(t^*)^{-1}\circ {\bf c}(t),~|t-t^*|\ll 1$, by a direct
calculation
\begin{align*}
&\frac{d}{dt}{\Log}~{\bf c}(t^*)^{-1}\circ {\bf c}(t)_{|t=t^*}\\
&=-\sqrt{-1} \int\limits_{-\infty}^{0}
\left\{(u-\Id)^{-2}\circ 
\frac{d}{dt}\left({\bf c}(t^*)^{-1}\circ {\bf c}(t)\right)_{|t=t^*}\right\}du\\
&=\dot{R}_{t^*}.
\end{align*}
\end{rem}
 
\begin{prop}\label{prop:M-sign}
Let $\{{\bf c}(t)\}$ be a path in $\cal{U}_{\cal{F}}(H_J)$ of class $C^1$
and $0<t^*<1$ a regular crossing. Then there exists a real $\delta>0$
such that
\begin{equation*}
\mathbf{M}(\{{\bf c}(t)\}_{|t-t^*|\le\delta})=\sign\tilde{Q}_{\frak{M}}.
\end{equation*}
\end{prop}

Before proving this proposition we give a lemma which describes
a behavior of eigenvalues closed to zero of a family of selfadjoint 
Fredholm operators under a certain non-degeneracy condition
(see \cite{Ka}):
\begin{lem}\label{lem:Kato1}
Let $\{A_t\}_{|t|\ll 1}$ be a $C^1$-class family of selfadjoint
Fredholm operators on a Hilbert space $H$.  
Assume that the symmetric bilinear form on $\Ker A_0$
$$
Q(x,y)=\frac{d}{dt}(x,~A_t(y))_{|t=0}
=(x,~\dot{A}_0(y)),~x,y\in \Ker A_0
$$
is non-degenerate. Then there exists
a positive number $\delta$ such that
for $0<t\le\delta$ there exist $p$ positive eigenvalues
and $q$ negative eigenvalues of the operator $A_t$, where
$p-q$ = $\sign Q$, $p+q=\dim\Ker A_0$.
Also for $-\delta\le t<0$ the opposite situations hold.
\end{lem}
\begin{proof}       
{}From the Fredholmness assumption of the continuous family $\{A_t\}$
there exist positive numbers $\delta$
and  $\varepsilon$ such that the projection operators $P_t$ for 
$|t|\le\delta$ defined by
\begin{equation}
P_t=\frac{1}{2\pi \sqrt{-1}}\int\limits_{|u|=\varepsilon} (u-A_t)^{-1}du
\end{equation}
have the constant rank equal to $\dim\Ker A_0$, and the range of each
$P_t$ = $\sum\limits_{|u|<\varepsilon}\Ker~(A_t-u)$.
By the approximation arguments we know that
the bilinear forms
$$
\left(\frac{1}{t}\cdot A_t\circ P_t(x),~P_t(y)\right),~x,y\in \Ker A_0
$$
and
$$
\left(\frac{d}{dt}{A_t}_{|t=0}(x),~y\right),~x, y\in\Ker A_0
$$
are uniformly close. In fact, for $x,y \in \Ker A_0$, 
\begin{align*}
&\left(\frac{1}{t}\cdot A_t\circ P_t(x),~P_t(y)\right)
-\left(\dot{A}_0(x),~y\right)\\
&= \left(\left(\frac{1}{t}
\cdot(A_t-A_0)-\dot{A}_0\right)\circ P_t(x),~P_t(y)\right)  \\
&+\left(\frac{1}{t}\cdot\left(P_t(x)-x\right),~A_0(P_t(y))\right) 
+\left(\dot{A}_0(P_t(x)-x),~P_t(y)\right)\\
&+\left(\dot{A}_0(x),~P_t(y)-y\right),
\end{align*}
and when $t\rightarrow 0$, 
$\norm{A_0(P_t(y))}\rightarrow 0$,
$\norm{\frac{1}{t}\cdot(P_t(x)-x)}$ is bounded,
and so these implies the assertion. Note here we used the fact that the
family $\{P_t\}$ is of class $C^1$. Hence  there exist 
$0<\delta_0\le\delta$ and  $0<\varepsilon_0\le \varepsilon$ and for 
$0<|t|\le\delta_0$ 
the signatures coincide and $A_t$ is an isomorphism between $P_t(H)$
= $P_t(\Ker A_0)$
which gives the existences of the $p+q$ eigenvalues of the operator
$A_t$, $0<e_1(t)\le e_2(t)\le \cdots\le e_p(t)\le\varepsilon_0$ and 
$-\varepsilon_0\le e_{-q}(t)\le e_{-q+1}(t)\le\cdots\le e_{-1}(t)<0$.
\end{proof}

\begin{proof}[Proof of Proposition \ref{prop:M-sign}]       
By the assumption there are a complex number $e^{\sqrt{-1}\theta_0}$
(close enough to $e^{\sqrt{-1}\pi}$, but $\ne e^{\sqrt{-1}\pi}$) 
and $\varepsilon >0$ such that for $|t-t^*|\le\varepsilon$ the 
operators $e^{\sqrt{-1}\theta_0}- {\bf c}(t)$ are invertible and
$$
\sum\limits_{|\theta|\le\varepsilon}
\dim\Ker ~({\bf c}(t) - e^{\sqrt{-1}(\pi+\theta)})<\infty.
$$

Put ${\bf c}(t+t^*)={\bf c}({t^*})e^{\sqrt{-1}R_{t}}$ and
let $A_t$ be a selfadjoint operator defined by the transformation
$$
A_t=
\sqrt{-1}(e^{\sqrt{-1}\theta_0}-{\bf c}(t+t^*))^{-1}(e^{\sqrt{-1}\theta_0}+{\bf c}(t+t^*))
-\sqrt{-1}\frac{e^{\sqrt{-1}\theta_0}-1}{e^{\sqrt{-1}\theta_0}+1},
$$
then $\{A_t\}_{|t|\le\varepsilon}$ 
is a $C^1$-class family of Fredholm operators.

The derivative $\dot{A}_0$ is given by a calculation using
the resolvent equation :
\begin{equation}\label{eq:derivative2}
\dot{A}_0=\left((e^{\sqrt{-1}\theta_0}-{\bf c}(t^*))^{-1}\right)^*
\circ 2\dot{R}_0\circ (e^{\sqrt{-1}\theta_0}-{\bf c}(t^*))^{-1}.
\end{equation}

This shows that the derivatives $\dot{A}_0$ and $2\dot{R}_0$
are conjugate, which gives us the coincidence of the signatures
of the two bilinear forms on $\Ker ({\bf c}(t^*)+ 1) =\Ker A_0$ :
for $x, y\in \Ker ({\bf c}(t^*)+ 1) =\Ker A_0$, 
\begin{align*}
&\left(\left((e^{\sqrt{-1}\theta_0}-{\bf c}(t^*))^{-1}\right)^*
\circ 2\dot{R}_0\circ (e^{\sqrt{-1}\theta_0}-{\bf c}(t^*))^{-1}(x),~y\right)\\
&=
\frac{2}{|e^{\sqrt{-1}\theta_0}+1|^2}(\dot{R}_0(x),~y) =(\dot{A}_0(x),~y).
\end{align*}
Hence by applying the preceding Lemma \ref{lem:Kato1} to the operator
family $\{A_t\}$ and returning back to the original 
family $\{{\bf c}(t)\}$ we have the desired numbers
of positive and negative eigenvalues of the family $\{{\bf c}(t)\}_{|t|\le\delta}$ for
sufficiently small $\delta$, which gives
$$
\mathbf{M}(\{{\bf c}(t)\}_{|t|\le\delta}) = \sign(\dot{A}_0)
=\sign(\dot{R}_0)
=\sign\tilde{Q}_{\frak{M}}.
$$
\end{proof}

\begin{rem}
For crossing $t^*=0$ or $1$, we only consider the
one-side differentiation in the definition of the crossing form.
In these cases we have
\begin{align*}
\mathbf{M}(\{{\bf c}(t)\}_{0\le t\le \delta})&=-q,\\
\mathbf{M}(\{{\bf c}(t)\}_{1-\delta\le t\le 1})&=p',
\end{align*}
where the signature of $\tilde{Q}_{\frak{M}}$ at $t^*=0$ is $(p,q)$ and
at $t^*=1$\, $(p',q')$.
\end{rem}

\begin{cor}\label{cor:2.11}
Let $\mu:\;I\rightarrow\cal{FL}_{\lambda}(H)$ be a $C^1$-class path
(so that $\cal{S}_{\lambda}\circ\mu(t)$ is a path in
$\cal{U}_{\cal{F}}(H_J)$ also of class $C^1$). Let
$0<t^*<1$ be a regular crossing of the curve 
$\{\cal{S}_{\lambda}\circ\mu(t)\}$.
Then there exist a $\delta>0$ such that
\begin{equation*}
\mathbf{Mas}(\{\mu(t)\}_{|t-t^*|\le\delta},\lambda)=\sign\tilde{Q}_{\frak{M}},
\end{equation*}
where $\tilde{Q}_{\frak{M}}$ 
denotes the crossing form of $\{\cal{S}_{\lambda}\circ\mu(t)\}$ 
at the time $t=t^*$.
\end{cor}
 

There is another bilinear form (see \cite{Du} and \cite{RS})
for describing the Maslov index
which will turn out to be more suitable for proving the 
{\it spectral flow formula} (see $\S 5$). 
It is based on a representation of $\mu$ as the graph
of a suitable bounded operator. Let $\mu:\;I\rightarrow
{\cal F}\LLa_{\lambda}(H)$ be a path in ${\cal F}\LLa_{\lambda}(H)$ of class
$C^1$ and let $0<t^*<1$ be a crossing of the curve 
$\{\cal{S}_{\lambda}\circ\mu(t)\}$, i.e., 
$\mu(t^*)\cap\lambda\ne\{0\}$. For $t$, $|t-t^*|\ll 1$,
$\mu(t)$ is transversal to $\mu(t^*)^{\perp}$
and in this neighborhood of $t^*$, each $\mu(t)$ can be written
as the graph of the bounded operator $A_t:\;
\mu(t^*)\rightarrow\mu(t^*)$, 
$\mu(t) = \{x+J\circ A_t(x)~|~x\in\mu(t^*)\}$. 
Note that the curve $\{A_t\}$
is also of class $C^1$. We consider the bilinear form
\begin{equation}\label{eq:second bilinear form}
Q_{\frak{M}}(x,y)=\left.\frac{d}{dt}\omega(x,J\circ A_t(y))\right|_{t=t^*}
\quad\text{for $x,y\in\mu(t^*)$.}
\end{equation}

In the above definition of the bilinear form $Q_{\frak{M}}$ we used
the fact that the inner product in the Hilbert space is compatible
with the symplectic form $\omega$, so that $\mu(t^*)^{\perp}$ is a
Lagrangian subspace, But this is not essential. In fact, let $\nu$ be a
Lagrangian subspace which is transversal to $\mu(t^*)$, then for
sufficiently small $|t-t^*|\ll 1$ Lagrangian subspaces $\{\mu_t\}$
are transversal to $\nu$. Then there is again a differentiable
family of bounded operators $\{\phi_t\}_{|t-t^*|\ll 1}$,
$\phi_t:\mu(t^*)\to \nu$, by which
we have $\mu(t)$ = graph of $\phi_t$ for each $t, ~|t-t^*|\ll 1$.

Now let $y\in \mu(t^*)$, then we have
\begin{equation}\label{eq:two-graph}
y+\phi_t(y) =~z_t +J\circ A_t(z_t),
\end{equation}
where $z_t = {\cal P}_{\mu(t^*)}(y+\phi_t(y))$ is a differentiable family
in $\mu(t^*)$ and $z_{t^*} =~y$.  Hence by differentiating the
both sides of (\ref{eq:two-graph}) we have
$$
{\frac{d}{dt}\phi_t}_{|t=t^*}(y) 
= {\cal P}_{\mu(t^*)}\left({\frac{d}{dt}\phi_t}_{|t=t^*}(y)\right)
+ J\circ {\frac{d}{dt}A_t}_{|t=t^*}(y).
$$
By this equality  
we have the invariance of the definition of the bilinear form 
$Q_{\frak{M}}$ from the auxiliary fixed Lagrangian subspace $\nu$:
\begin{prop}\label{prop:invariance from aux-Lag}
$$
Q_{\frak{M}}(x,~y) =
\frac{d}{dt}\omega(x, \phi_t(y))_{|t=t^*},~x,~y\in ~\mu(t^*).
$$
\end{prop}
\begin{proof}
Let $x,~y\in \mu(t^*)$, then
\begin{align*}  
&\frac{d}{dt}\omega(x, \phi_t(y))_{|t=t^*}\\
&=\omega
\left(x,~{\cal P}_{\mu(t^*)}\left({\frac{d}{dt}\phi_t}_{|t=t^*}(y)\right)
+ J\circ {\frac{d}{dt}A_t}_{|t=t^*}(y)\right)\\
&=\omega\left(x,~J\circ {\frac{d}{dt}A_t}_{|t=t^*}(y)\right)
=Q_{\frak{M}}(x,~y).
\end{align*}
\end{proof}

The bilinear form 
$Q_{\mathfrak{M}}$ 
is a symmetric bilinear form on $\mu(t^*)$ at each point $t^*$ 
solely defined by
the differentiable family $\{\mu_t\}$ itself
and we will show the following coincidence of the signature of
the bilinear forms.

\begin{prop}
On $\mu(t^*)\cap\lambda$, 
$\sign Q_{\frak{M}}=
\sign\tilde{Q}_{\frak{M}}$.
\end{prop}

\begin{proof}
We have two expression of the space $\mu(t)$:
\begin{enumerate}
\item
$\mu(t) = \{x+J\circ A_t(x)~|~x\in\mu(t^*)\}$, 
$A_t\in {\widehat{\cal B}}(\mu(t^*))$, $\{A_t\}$ is $C^1$-class.
\item
$\mu(t)=U_t(\lla^{\perp})$, where $\{U_t\}$
is a $C^1$-class family of unitary operators on $H_J$. 
\end{enumerate}

Put $U_{t+t^*}=U_{t^*}e^{\sqrt{-1}S_t}$ 
and ${\bf c}(t+t^*)={\cal S}_{\lla}(\mu(t+t^*)) = {\bf c}({t^*})e^{\sqrt{-1}R_t}$,
where $\{S_t\}$ and $\{R_t\}$ are $C^1$-class families
of selfadjoint operators on $H_J$. We represent $S_t =X_t + \sqrt{-1}Y_t$
with $X_t, Y_t \in {\cal B}(\mu(t^*))$, $X={\,^tX}$ and $Y=-{\,^tY}$. 

By differentiating  ${\bf c}(t)=U_t\circ\theta_{\lla}(U_t)$ at $t=t^*$ we have
\begin{align*}
&\frac{d}{dt}{\bf c}(t)_{|t=t^*}= \dot{\bf c}(t^*)\\
&=U_{t^*}\circ\sqrt{-1}\dot{S}_{0}\circ\theta_{\lla}(U_{t^*})
+U_{t^*}\circ\theta_{\lla}({U}_{t^*}\circ\sqrt{-1}\dot{S}_{0})\\
&=U_{t^*}(\sqrt{-1}(\dot{X}_0+\sqrt{-1}\dot{Y}_0)
+\sqrt{-1}(\dot{X}_0-\sqrt{-1}\dot{Y}_0))\circ\theta_{\lla}(U_{t^*}) \\
&=2\sqrt{-1}U_{t^*}\circ\dot{X}_0\circ\theta_{\lla}(U_{t^*})\\
&=\sqrt{-1}{\bf c}(t^*)\dot{R}_0.
\end{align*}

This identity says that the bilinear form $\tilde{Q}_{\frak{M}}$ 
on $H_J$ defined by $\dot{R}_0$ 
is unitary equivalent to the complexification of the bilinear form 
defined by the real selfadjoint
operator $2\dot{X}_0$ on $\mu(t^*)$.

Now by differentiating the equality 
$$
U_t(x)= 
{\cal P}_{\mu(t^*)}(U_t(x)) +J\circ A_t
({\cal P}_{\mu(t^*)}(U_t(x))),~x\in \lla^{\perp},
$$
we have
\begin{equation*}
{\cal P}_{\mu(t^*)}\circ \dot{S}_0(U_{t^*}(x))=\dot{A}_0(U_{t^*}(x)).
\end{equation*}
Note that we used here the equation 
$J\circ{\cal P}_{\mu(t^*)}+{\cal P}_{\mu(t^*)}\circ J =J$.

Let $x,~y\in \mu(t^*)$, then we have
\begin{align}\label{eq:bilinear2}
&\omega(x,J\circ\dot{A}_0(y))=<x,\dot{A}_0(y)>\\
&=<x, {\cal P}_{\mu(t^*)}\circ \dot{S}_0(y)>
= <x, \dot{X}_0(y)>.\notag
\end{align}

Hence the unitary equivalence (on the whole space $H_J$) 
of the bilinear forms defined
by the operators $\dot{R}_0$ and $2\dot{X}_0$ 
and the equation (\ref{eq:bilinear2}) (note that
the identity holds on $\mu(t^*)$) show the proposition. 
\end{proof}

\begin{rem}\label{orthogonal-inner-product} 
The unitary equivalence of the two bilinear forms
${Q}_{\frak{M}}$ and $\tilde{Q}_{\frak{M}}$ on $\mu(t^*)$
implies that the definition of 
the bilinear form  $\tilde{Q}_{\frak{M}}$
does not depend on the almost complex structure $J$ by which
we regard the real Hilbert space $H$ as a complex Hilbert space $H_J$.
This means we can freely replace the inner product by a suitable one.
For example, we can assume that 
any two transversal Lagrangian subspaces are orthogonal
(see Proof of Theorem \ref{thm:reduction2}).
\end{rem}

Now we have a similar formula with Proposition \ref{prop:M-sign}:

\begin{cor}\label{cor:Mas-formula2}
Let $\mu:\;I\rightarrow\mathcal{F}\Lambda_{\lambda}(H)$ be a $C^1$-class path.
Let
$0<t^*<1$ be a regular crossing of the curve. Then it is also regular
crossing of the curve
$\{\cal{S}_{\lambda}\circ\mu(t)\}$,
and there exists a positive $\delta>0$ such that
\begin{equation*}
\mathbf{Mas}(\{\mu(t)\}_{|t-t^*|\le\delta},\lambda)=\sign{Q}_{\frak{M}},
\end{equation*}
where ${Q}_{\frak{M}}$ 
denotes the crossing form of $\{\mu(t)\}$ 
at the time $t=t^*$.
\end{cor}

\begin{rem}
In the paper \cite{RS} the authors gave a definition of the Maslov
index (in the case of finite dimension) 
for such differentiable curves $\{{\bf c}(t)\}_{t\in[0,1]}$
that all their ``crossings''
are regular in terms of this bilinear form with
corrections at the end points by adding
the halves of the dimensions $\dim \lambda\cap {\bf c}(1)$
and $\dim \lambda\cap {\bf c}(0)$.  
\end{rem}

Finally in this subsection we give an example of a $C^1$-class
path with a regular crossing and calculate the Maslov index.
\begin{example}\label{ex:regular cross-path}
Let $F$ be a finite dimensional subspace in $J(\lla)$, and we define
a family of unitary operators such that
\begin{equation}
U(t)(x)=
\left\{\begin{array}{ll}
e^{\sqrt{-1}\pi t}\cdot x &\qquad x ~\in F\\
x& \qquad x~\in F^{\perp}\cap J(\lla)
\end{array}\right.
\end{equation}

For each $t,~t\in~[0,~1]$ let $\mu(t)= ~U(t)(J(\lla))$, then
$\mu(t)\in {\cal F}\LLa_{\lla}(H)$ and $t=1/2$ is an only non-trivial
crossing with $\lambda$ and is regular.  As is easily determined 
the crossing form is given by
\begin{equation}  
{Q}_{\frak{M}}(x,~y) = \pi <x,~y>,~x,~y\in J(F). 
\end{equation}
Hence we have
\begin{equation*}
\mathbf{Mas}(\{\mu(t)\}_{0\le t\le 1},\lambda)
=\sign{Q}_{\frak{M}}=\dim F.
\end{equation*}
Also for $0< \epsilon \ll 1$,
\begin{equation*}
\mathbf{Mas}(\{\mu(t)\}_{{1/2}\le t\le {1/2+\epsilon}},\lla)= 0,
\end{equation*}
and 
\begin{equation*}
\mathbf{Mas}(\{\mu(t)\}_{1/2-\epsilon\le t\le 1/2},\lambda)= \dim F.
\end{equation*}
\end{example}

\subsection{Maslov index for paths of Fredholm pairs of Lagrangian subspaces}
In this subsection we will denote the direct sum of the symplectic
Hilbert space $(H,\omega)$ and $(H, -\omega)$ with the notation 
$\Hi=H\boxplus H$ $\equiv$ $H_{\omega}\boxplus H_{-\omega}$. 
$\Hi$ is a symplectic Hilbert space with the symplectic form
$\Omega = \omega\, -\, \omega$, 
and the corresponding almost complex structure
$\J = J\boxplus -J$, so that we have 
$\Hi_{\J}=H_J\boxplus H_{-J}$.
Let $\{(\mu_t,\lambda_t)\}_{t\in I}$ be a continuous family of
Fredholm pairs of Lagrangian subspaces, then 
$\{\mu_t\boxplus\lambda_t\}$
is a curve in $\cal{F}\Lambda_{\Delta}(H_{\omega}\boxplus H_{-\omega})$,
where $\Delta$ is the diagonal of $H\boxplus H$. 
Of course it is natural to define the Maslov index of the curve
of Fredholm pairs $\{(\mu_t,\lambda_t)\}$ to be
$\Mas(\{\mu_t\oplus\lambda_t\},\Delta)$. 

\begin{prop}\label{prop:pair}
Let $\{\mu_t\}$ be a continuous curve in $\cal{F}\Lambda_{\lambda}
(H_{\omega})$, then
\begin{equation*}
\Mas(\{\mu_t\},\lambda)=\Mas(\{\mu_t\oplus\lambda\},\Delta).
\end{equation*}
\end{prop}
\begin{rem}
For loops this property  will be well-known. For arbitrary
continuous paths in the finite dimensional case 
this can be proved by making use of Proposition 
\ref{leray-kashiwara-maslov}($\S 3.1$),
but in the infinite dimensional case we 
have no such relations and we 
need a proof which is valid not only 
for loops but also
for any continuous paths.
\end{rem}

If we identify 
$\Hi_{\J}= \Delta+\Delta^{\perp}
=\Delta+\J(\Delta)\cong\Delta\otimes \C$, then 
$\tau_{\Delta}(a\boxplus b)=b\boxplus a$.
Let us decompose $H$ as $H= \lambda\oplus\lambda^{\perp}$ and let $\varphi:\;\Delta
\rightarrow\Delta^{\perp}$  be
\begin{equation*}
\varphi((x,y)\boxplus(x,y))=(-x,y)\boxplus(x,-y)
\end{equation*}
where we express elements in $\Delta$ by $(x,y)\boxplus (x,y)$,
\,$x+y\in\lambda+\lambda^{\perp}=H$. Then we have
\begin{equation*}
\graph\varphi=\lambda^{\perp}\boxplus\lambda.
\end{equation*}

Let $A=\J\circ\varphi:\;\Delta\rightarrow\Delta$ and $V:\;\Hi_{\J}
\rightarrow\Hi_{\J}$ by
\begin{equation*}
V=\frac{-\sqrt{-1}}{\sqrt{2}}-\frac{A\otimes\Id}{\sqrt{2}}
\end{equation*}
where we regard $A=A\otimes Id$ is complexified according 
to the identification $\Hi_J\cong\Delta\otimes\C$. Then we have
\begin{equation}
\sqrt{-1}(A\otimes\Id)((a,b)\boxplus (c,d)) = (c,-d)\boxplus(-a,b)
\end{equation}
for $(a,b)\boxplus (c,d)\in H_J\boxplus H_{-J} =
(\lambda+\lambda^{\perp})\boxplus (\lambda+\lambda^{\perp})$ 
and
\begin{equation}
V(\Delta^{\perp})=\lambda^{\perp}\boxplus\lambda.
\end{equation}

Now we define maps ${\bf{a}}_{\lambda}, \,{\bf{b}}_{\lambda}$ and
$P_{\lambda}$ as follows:
\begin{equation*}
\begin{array}{rcc}
{\bf{a}}_{\lambda}:\;\cal{U}_{\lambda}(H_J)&\longrightarrow&\cal{U}_{\Delta}(\Hi_{\J})\\
U & \mapsto & \tilde{U}\circ V
\end{array}
\end{equation*}
where $\tilde{U}=U\boxplus\Id\,:\,H_J\boxplus H_J \longrightarrow
H_J\boxplus H_{-J}$,
\begin{equation*}
\begin{array}{rcc}
{\bf{b}}_{\lambda}:\;\cal{U}_{\cal{F}}(H_J)&\longrightarrow &\cal{U}_{\cal{F}}
(\Hi_{\J})\\
W & \mapsto & \sqrt{-1}\cdot W\circ (A\otimes\Id),
\end{array}
\end{equation*}
and
\begin{equation*}
\begin{array}{rcc}
P_{\lambda}\,:\,\cal{F}\Lambda_{\lambda}(H)&\rightarrow&\cal{F}\Lambda_{\Delta}(\Hi)\\
\mu & \mapsto & \mu\boxplus\lambda.
\end{array}
\end{equation*}

\begin{lem}
The following diagram is commutative.
\begin{equation*}
\begin{CD}
\cal{U}_{\lambda}(H_J) @>{{\bf{a}}_{\lambda}}>>\cal{U}_{\Delta}(H_J\boxplus H_{-J})\\
@V{\rho_{\lambda}}VV  @VV{\rho_{\Delta}}V\\
\cal{F}\Lambda_{\lambda}(H) @>>{P_{\lambda}}>
\cal{F}\Lambda_{\Delta}(H_{\omega}\boxplus H_{-\omega})\\
@V{\cal{S}_{\lambda}}VV @VV{\cal{S}_{\Delta}}V\\
\cal{U}_{\cal{F}}(H_J)@>{{\bf{b}}_{\lambda}}>>\cal{U}_{\cal{F}}(H_J
\boxplus H_{-J})
\end{CD}
\end{equation*}
\end{lem}

\begin{proof}
It will be enough to prove 
$\cal{S}_{\Delta}\circ P_{\lambda}={\bf{b}}_{\lambda}\circ
\cal{S}_{\lambda}$. Since  
$\theta_{\Delta}(V) = V$, $V^2 = \sqrt{-1}\cdot A\otimes\Id$ and
$\theta_{\Delta}({\widetilde{U}}) = \Id\boxplus U^{*}$ we have
\begin{eqnarray*}
&&\cal{S}_{\Delta}\circ\rho_{\Delta}({\bf{a}}_{\lambda}(U))\\
&=&\tilde{U}\circ V\circ\tau_{\Delta}\circ(\tilde{U}\circ V)^*\circ
\tau_{\Delta}\\
&=&\tilde{U}\circ\sqrt{-1}(A\otimes\Id)\circ\theta_{\Delta}(\tilde{U})\\
&=&{\widetilde{U\circ\theta_{\lambda}(U)}}\circ \sqrt{-1}(A\otimes\Id),
\end{eqnarray*}
which prove the commutativity of the diagram.
\end{proof}

\begin{proof}[Proof of Proposition \ref{prop:pair}]
{}From the above lemma we can show that if $E$ is an eigenvalue of
$\cal{S}_{\Delta}(\rho_{\Delta}({\bf{a}}_{\lambda}(U)))$, then $-E^2$
is an eigenvalue of $\cal{S}_{\lambda}\circ\rho_{\lambda}(U)$.
Conversely if $l=e^{\sqrt{-1}\sigma}$ is an eigenvalue of $\cal{S}_{\lambda}\circ
\rho_{\lambda}(U)$, then only one of $\pm e^{\sqrt{-1}(\pi+\sigma)}$
is closed to $-1$. So if we have 
a continuous curve $\{\mu_t\} \subset \cal{F}\Lambda_{\lambda}(H)$,
then the numbers of eigenvalues of $\{\cal{S}_{\lambda}(\mu_t)\}$ and
$\{\cal{S}_{\Delta}(\mu_t\boxplus\lambda)\}$ which across
$e^{\sqrt{-1}\pi}$ 
coincide
in both directions. 
This proves the proposition.
\end{proof}

The next property will be also natural:

\begin{prop}
$\Mas(\{\mu_t\boxplus\lambda_t\},\Delta)=-\Mas(\{\lambda_t\boxplus
\mu_t\},\Delta)$.
\end{prop}

\section{Finite dimensional cases}

In the finite dimensional cases, 
the Maslov index for arbitrary continuous paths in the
Lagrangian-Grassmannian
manifolds was first defined in the paper [Go1] 
by noting the extendibility of the ``{\it Leray index}'' for
arbitrary pairs of points on the universal covering 
space of the Lagrangian Grassmannian by making use of 
the cocycle condition of the ``{\it Leray index}''.
Conversely, first we define Maslov index for 
arbitrary paths with
respect to a Maslov cycle as we gave above, 
then we can define the ``{\it Leray index}'' for arbitrary
pairs of points on the universal covering of the 
Lagrangian Grassmannian (Proposition \ref{leray-kashiwara-maslov}).

In the infinite dimensional case we could define the Maslov index 
for arbitrary paths
with respect to a Maslov cycle as we did in the above 
Definition \ref{def:M-index-path},
however we can not define ``{\it Kashiwara index}'' 
for arbitrary triples of Lagrangian subspaces 
like the finite dimensional case,
although we have a symmetric bilinear form similar 
to the finite dimensional case. We
can define it for mutually almost coincident triples, 
since then the symmetric bilinear
form is of finite rank. 
Nor we can define Leray index for arbitrary pairs of points
on the universal covering space of the 
Fredholm-Lagrangian-Grassmannian. 

In this section, following \cite{Go1}
we summarize 
the mechanism for defining the Maslov index for
paths in $\Lambda(H)$
of the finite dimensional symplectic vector space $H$,
and give an extension of the ``Kashiwara index'' 
to arbitrary triples of unitary
operators($\S 3.2$).

\subsection{Leray index and Kashiwara index}

Let $\ell_1$, $\ell_2$ and $\ell_3$ be three Lagrangian subspaces
and define the quadratic form $Q$ on the direct sum 
$\ell_1\oplus \ell_2\oplus \ell_3$ as follows:
\begin{equation}
Q(x,x',x'')=\omega(x,x')+\omega(x',x'')+\omega(x'',x),
~x\in\ell_1,~x'\in\ell_2,~x''\in\ell_3.
\end{equation}

The corresponding bilinear form is

\begin{align*}
&I_{\omega}(x,x',x''~; ~a,a',a'')\\
&\quad =\omega(x,a')+\omega(x',a)
+\omega(x,a'')+\omega(x'',a)+\omega(x',a'')+\omega(x'',a'),\\
&\qquad\qquad x,~a\in\ell_1,~x',~a'\in\ell_2,~x'',~a''\in\ell_3.
\end{align*}

The index  of this quadratic form is called as 
``{\it Kashiwaka index}''
or ``{\it cross index}'' 
of the triple of Lagrangian subspaces.
We denote it by $\sigma(\ell_1,~\ell_2,~\ell_3)$.


In the finite dimension cases, 
although $\Lambda_{\lambda}(H)=\Lambda(H)$
always, it should be noted that the Souriau map
$\mathcal{S}_{\lambda}:\Lambda(H)\to \mathcal{U}(H_J)$ itself  
depends on the chosen Lagrangian subspace $\lambda$.
Now let 
$$
\widetilde{\mathcal{U}}(H_J)
=\{(U,\alpha)\in\mathcal{U}(H_J)\times \mathbb{R}
~|~ \det U =e^{\sqrt{-1}\alpha}\}.
$$
be a realization of the universal covering of
the unitary group $\mathcal{U}(H_J)$, then the space
$\mathcal{S}_{\lambda}^{\,\,*}(\widetilde{\mathcal{U}}(H_J))
=\widetilde{\Lambda}(H)
=\{(\mu,\alpha)\in\Lambda(H_J)\times\mathbb{R}
~|~\det\mathcal{S_{\lambda}(\mu)}
=e^{\sqrt{-1}\alpha}\}$
is the universal covering of the Lagrangian Grassmannian
$\Lambda(H)$ with the projection map
$q_{\lambda}:\widetilde{\Lambda}(H)\to \Lambda(H)$.
Let $\widetilde{\ell_1}$ and $\widetilde{\ell_2}$ be two point
on $\mathcal{S}_{\lambda}^{\,\,*}(\widetilde{\mathcal{U}}(H_J))
=\widetilde{\Lambda}(H)$ 
and we assume
that $q_{\lambda}(\widetilde{\ell_1})=\ell_1$ and
$q_{\lambda}(\widetilde{\ell_2})=\ell_2$ are transversal, i.e.,
$\ell_1\cap\ell_2=\{0\}$.  Then 
$Id +\mathcal{S}_{\lambda}(\ell_1)\mathcal{S}_{\lambda}(\ell_2)^*$
is invertible, so we define ``{\it Leray index}'' 
$\mu(\widetilde{\ell_1},~\widetilde{\ell_2})$ of such a pair by
\begin{defn}
$$
\mu(\widetilde{\ell_1},~\widetilde{\ell_2})=
\frac{1}{2\pi}\left(\alpha_1-\alpha_2+\sqrt{-1}\Tr 
~\Log \left( -\mathcal{S}_{\lambda}(\ell_1)
\circ\mathcal{S}_{\lambda}(\ell_2)^*\right)\right),
$$
where $\Log$ is defined by (\ref{Log}).
\end{defn}

We have a fundamental relation of the Leray index and the 
Kashiwara index(cocycle condition of the Leray index):
\begin{prop}
Let $\widetilde{\ell_1}$,
$\widetilde{\ell_2}$ and $\widetilde{\ell_3}$
be three points on $\widetilde{\Lambda}(H_J)$ such that
each of the pairs 
$(q_{\lambda}(\ell_1), q_{\lambda}(\ell_2))$, 
$(q_{\lambda}(\ell_1), q_{\lambda}(\ell_3))$ 
and
$(q_{\lambda}(\ell_2), q_{\lambda}(\ell_3))$ 
is transversal, then
\begin{equation}
\mu(\widetilde{\ell_1},\widetilde{\ell_2})+
\mu(\widetilde{\ell_2},\widetilde{\ell_3})+
\mu(\widetilde{\ell_3},\widetilde{\ell_1})
=\sigma(q_{\lambda}(\ell_1),~q_{\lambda}(\ell_2),~q_{\lambda}(\ell_3)).
\end{equation}
\end{prop}

Then we define for any pairs 
$(\widetilde{\ell_1},~\widetilde{\ell_2})
\in\widetilde{\Lambda}(H)\times \widetilde{\Lambda}(H)$
(without transversality assumption between the pair 
$(q_{\lambda}(\widetilde{\ell}_1),~q_{\lambda}(\widetilde{\ell}_2))$)
the ``{\it Leray index}''
$\mu(\widetilde{\ell_1},~\widetilde{\ell_2})$
by the formula
\begin{equation}
\mu(\widetilde{\ell_1},\widetilde{\ell_2})
=\mu(\widetilde{\ell},\widetilde{\ell_2})
-\mu(\widetilde{\ell},\widetilde{\ell_1})
+\sigma(q_{\lambda}(\ell_1),~q_{\lambda}(\ell_2),~q_{\lambda}(\ell)).
\end{equation}
by taking an element $\widetilde{\ell}$ in 
$\widetilde{\Lambda}(H)$ 
such that 
$q_{\lambda}(\widetilde{\ell})$
is transversal to each $q_{\lambda}({\widetilde{\ell_i}})$ ($i=1,2$).

{\it The independence of this value from the choice of such 
$\widetilde{\ell}$ is proved by making use of the fact} 
\begin{equation}\label{cocycle-condition-sigma}
\partial \sigma(\ell_0,\ell_1,\ell_2,\ell_3)=
\sigma(\ell_1,\ell_2,\ell_3)-\sigma(\ell_0,\ell_2,\ell_3)
+\sigma(\ell_0,\ell_1,\ell_3)-\sigma(\ell_0,\ell_1,\ell_2)=0.
\end{equation}


Now we fix a $\lambda\in\Lambda(H)$ and
let $\{{\bf c}(t)\}_{t\in [0,1]}$ be a continuous curve in
$\Lambda(H)$, and take a lifting 
$\{\widetilde{{\bf c}}(t)\}_{t\in [0,1]}$
of the curve $\{{\bf c}(t)\}$. Then we have 
\begin{prop}\label{leray-kashiwara-maslov}
$$
\mathbf{Mas}(\{{\bf c(t)}\},\lambda)
=\mu(\widetilde{{\bf c}}(0),~\widetilde{{\bf c}}(1))
-\sigma(\lambda,~{\bf c}(1),~{\bf c}(0)).
$$
\end{prop}

\subsection{Complex Kashiwara index}

By the very definition of the Leray index
we see that it can be defined, by the same formula,
for any pairs of the points $((U,\alpha_1), (U,\alpha_2))$
in $\widetilde{U}(H_J)$ with the property
that $Id +U_1\circ U_2^{-1}$ is invertible:
$$
\mu((U,\alpha_1), (U,\alpha_2))
=\frac{1}{2\pi}\left(\alpha_1-\alpha_2+\sqrt{-1}\Tr 
~\Log \left( -U_1\circ U_2^{-1}\right)\right).
$$ 
Then for such triples the sum
$$
\mu((U,\alpha_1), (U,\alpha_2))+\mu((U,\alpha_2), (U,\alpha_3))
+\mu((U,\alpha_3), (U,\alpha_1))
$$
is independent of $(\alpha_1,~\alpha_2,~\alpha_3)$.
This enables us to extend the Kashiwara index for any triples of
unitary matrices. We explain it here.

We denote the sesquilinear extension of the symplectic
form $\omega$ by $\omega^{\mathbb{C}}$.

For three Lagrangian subspaces 
$L_i\in\Lambda^{\C}(H\otimes\C)$ ($i=1,2,3$) 
we define the bilinear form 
\begin{equation}
I_{\omega}^{\mathbb{C}}:
L_1\oplus L_2 \oplus L_3 \times
L_1\oplus L_2 \oplus L_3\longrightarrow \mathbb{C},
\end{equation}
by
\begin{align*}
&I_{\omega}^{\C}(z,z',z''~;~w,w',w'')\\
&\quad =\omega^{\C}(z,w')+\omega^{\C}(w,z')+
\omega^{\C}(z,w'')+\omega^{\C}(z'',w)+
\omega^{\C}(z',w'')\omega^{\C}(z'',w')\\
&\qquad z=x+\sqrt{-1}y,~w=a+\sqrt{-1}b\in L_1,
~z',~w'\in L_2,~z'',~w''\in L_3.
\end{align*}

Then this is an Hermite form on $L_1\oplus L_2 \oplus L_3$.
We denote the signature of this form by 
$$
\sigma^{\C}(L_1,~L_2,~L_3)
$$
and call it as 
``complex Kashiwara index'' 
or 
``complex cross index''.
If each 
$L_i=\lambda_i\otimes\C$ with $\lambda_i\in \Lambda(H)$,
then this is the sesquilinear extension 
of the bilinear form $I_{\omega}$
and their signatures coincides.

Next recall the isomorphism $\varPhi_{\lambda}$ (\ref{varPhi}):
  
$$\varPhi_{\lambda} : \mathcal{U}(H_J)\rightarrow 
\Lambda^{\C}(H\otimes\C),
$$
$$
\varPhi_{\lambda}(V) ~= 
~\text{the graph of the unitary operator}\, 
\, -{\mathfrak{k}}\circ V \circ \tau_{\lambda}\circ
{\mathfrak{K}}^{-1} \in \cal{U}(E_+,E_{-}).
$$

\begin{prop}
Let $(U_i,\alpha_i)\in \widetilde{\mathcal{U}}(H_J)$, $i=1,2,3$, i.e.,
$\det U_i=e^{\sqrt{-1}\alpha_i}$, and assume
$\det(U_i\circ U_j^{\,-1}+Id)\not= 0, ~i,~j =1,2,3$, then
\begin{align*}
&\mu((U_1,\alpha_1),(U_2,\alpha_2))+\mu((U_2,\alpha_2),(U_3,\alpha_3))
+\mu((U_3,\alpha_3),(U_1,\alpha_1))\\
&\quad =\sigma^{\C}
(\varPhi_{\lambda}(U_1),\varPhi_{\lambda}(U_2),\varPhi_{\lambda}(U_3)).
\end{align*}
Especially the value of $\sigma^{\C}$ does not depend
on the fixed $\lambda$.
\end{prop}

Now let $(U_1,\alpha_1)$ and $(U_2,\alpha_2)$ be any pair
of the points in $\widetilde{\mathcal U}(H_J)$ and choose
an element $(U,\alpha)$ such that 
\begin{equation}\label{U-transversal-condition}
\det(Id + U_1\circ U^{-1})\not=0~\text{and}~
\det(Id + U_2\circ U^{-1})\not=0.
\end{equation} 
Then we can define
the ``Leray index'' 
of the pair 
$(U_1,\alpha_1),~(U_2,\alpha_2)$ 
by the formula:
\begin{defn}
\begin{align*}
&\mu((U_1,\alpha_1),~(U_2,\alpha_2))\\
&\quad =\mu((U,\alpha),~(U_2,\alpha_2))-
\mu((U,\alpha),~(U_1,\alpha_1))
+\sigma^{\C}(\varPhi_{\lambda}(U_1),\varPhi_{\lambda}(U_2),
\varPhi_{\lambda}(U)).
\end{align*}
\end{defn}
The similar cocycle property of $\sigma^{\C}$
with (\ref{cocycle-condition-sigma}) guarantees that
this definition does not depend on the choice
of the element $(U,\alpha)$ which satisfies the condition
(\ref{U-transversal-condition}), and the function
is defined on the space 
$\widetilde{\mathcal{U}}(H)\times
\widetilde{\mathcal{U}}(H)$. 
Also it is the extension
of the Leray index defined on the space
$\widetilde{\Lambda}(H)\times \widetilde{\Lambda}(H)$ 
through the embedding
$$
\widetilde{\Lambda}(H)\times \widetilde{\Lambda}(H)
~\hookrightarrow
~\widetilde{\mathcal{U}}(H_J)\times\widetilde{\mathcal{U}}(H_J).
$$
Note that this embedding is defined 
by choosing a $\lambda\in\Lambda(H)$.

\begin{rem}
We do not state the invariance of the Maslov index and other indexes
(finite and infinite dimensions) under the unitary
and symplectic group actions. These will be proved by making
use of the properties of Souriau maps.
\end{rem}

\section{Polarization and a Reduction Theorem of the Maslov index}

For the proof of Theorem \ref{thm:pi_1} we employed the finite
dimensional reduction (= the diagram \ref{e:CD_lim})
of the Maslov index from infinite dimensions to finite dimensions.
In this section we prove a reduction theorem \ref{thm:reduction2}
of the Maslov index 
inside the infinite dimensions. 

\subsection{Symplectic transformation and Canonical relation}

First we remark a continuity of a symplectic transformation.

Let $H_i~(i=0,1)$ be two symplectic Hilbert spaces equipped with
compatible symplectic structures 
(symplectic forms $\omega_i$, inner products $<\bullet,\bullet>_i$ 
and almost complex structures $J_i$ $(i=0,1)$). 
As in $\S 2.5$ we consider the direct sum $H_0\boxplus H_1$ as a 
symplectic Hilbert space with the compatible symplectic form 
$$
\Omega((x,a),(y,b))=\omega_0(x,y)-\omega_1(a,b), ~(x,a),~(y,b)\in
H_0\boxplus H_1.
$$

Let $S:H_0 \to H_1$ be a linear map defined on the whole
space $H_0$ and we assume that $S$ keeps the symplectic forms:
\begin{equation}\label{symplectic-transformation}
\omega_1(S(x),~S(y)) = \omega_0(x,~y), 
~\text{for all}~x~\text{and}~y\in ~H_0.
\end{equation}

Then it is easy to see that $S$ is injective and
the graph $G_S=\{(x,S(x))|x\in H_0\}$ is an isotropic subspace.
Under this assumption (\ref{symplectic-transformation}),
the closure of the image $S(H_0)$ is a symplectic Hilbert space.  So
now we start assuming that $S$ has a dense image, then
\begin{prop}
The graph $G_S$ is a Lagrangian subspace in $H_0\boxplus H_1$.
Hence $S$ is a bounded operator by the closed graph theorem.
\end{prop}  
 
\begin{proof}
Let $(a,b)\in H_0\boxplus H_1$ and $\Omega((a,b), (x,S(x))=0$
for all $x\in H_0$. Then we have 
$\omega_0(x,a)=\omega_1(S(x), b)=\omega_1(S(x),S(a))$. Hence $S(a)=b$,
which shows that the graph $G_S$ is a Lagrangian subspace in
$H_0\boxplus H_1$.
\end{proof}

By this proposition, if a symplectic transformation $S$
between two symplectic Hilbert spaces is algebraically
isomorphic, then it must be topologically isomorphic. 

We call a Lagrangian subspace $C$ in the direct sum
$H_0\boxplus H_1$ with the symplectic form $\Omega= \omega_0 -\omega_1$
as a canonical relation(see for the global settings 
\cite{Ho2},\cite{Ho3}).  
The graph of a symplectic 
transformation defined on the whole space $H_0$ with the dense image 
is so a canonical relation.

In this section
we consider a particular canonical relation given as a 
graph of a closed symplectic transformation, that is,
let $S$ be a densely defined closed and symplectic
transformation, and in particular not continuous:
$$
S:\mathfrak{D}_S\to
H_1,~\omega_0(x,y)=\omega_1(S(x),~S(y)),~x
~\text{and}~y\in\mathfrak{D},
$$
where $\mathfrak{D}_S$ 
is a dense subspace in the symplectic Hilbert space $H_0$.
Then again we have 
\begin{prop}
Let us assume that $S$ has the dense image, 
then the graph $C=S_G$ is a Lagrangian subspace in $H_0 \boxplus H_1$.
\end{prop}
 
Let $\lambda \in\Lambda(H_0)$, then we have always
$\lambda\cap\mathfrak{D}_S\not=\{0\}$ and $S(\lambda)$ is an isotropic
subspace, but will not be a Lagrangian subspace in general.  

We will show a theorem that if we restrict this map $S$ to a
Fredholm-Lagrangian-Grassmannian then the image is a 
Lagrangian subspace and it reserves the Maslov index 
under additional assumptions.

\subsection{Polarization of symplectic Hilbert spaces} 

Again let $H$ be a symplectic Hilbert space.

\begin{defn}
We say that the symplectic Hilbert space is polarized, when
$H$ is decomposed into a direct sum of two Lagrangian subspaces
$$
H=\ell_+\oplus \ell_-.
$$
Or we say that 
the sum $H=\ell_+\oplus \ell_-$ is a polarization of $H$. 
\end{defn}

\begin{rem} In the polarization $H=\ell_+\oplus \ell_-$
the subspaces need not be orthogonal, however it is possible
by replacing the inner product ({\it symplectic form should not be changed
always} that the sum is orthogonal. Of course the new norm 
is equivalent to the initial one
(Remark (\ref{orthogonal-inner-product})).
\end{rem}

\begin{prop}
Let $H =\ell_-+\ell_+$ be a polarized symplectic Hilbert space.  
\begin{enumerate}
\item
Let $S$ be a closed subspace in $\ell_-$, and 
we take a complement $T$
of $S$ in $\ell_-$, $\ell_- = S+T$.  Put
$F=T^{\circ}\cap\ell_+$, 
then $S+F$ is a symplectic subspace,
in fact,
 $(S+F)^{\circ}= T+ S^{\circ}\cap\ell_+$ and
$$
S+F+T+S^{\circ}\cap\ell_+=H.
$$
\item
Let $S$ be a closed subspace in $\ell_-$ and $F$ be a closed subspace in
$\ell_+$.
Assume that $S+F$ is a symplectic subspace, then 
$F^{\circ}\cap\ell_- + S^{\circ}\cap\ell_+$ is a symplectic subspace,
$$
(S+F)^{\circ}=F^{\circ}\cap\ell_- + S^{\circ}\cap\ell_+.
$$
\end{enumerate}
\end{prop}

\begin{proof}
Proof of (a).  Put $G =S^{\circ}\cap\ell_+$, then
\begin{align*}
&(S+F)^{\circ}\\
& = S^{\circ}\cap(T^{\circ}\cap\ell_+)^{\circ}\\
&= T+S^{\circ}\cap\ell_+ = T+G,
\end{align*}
since $T\subset S^{\circ}$.  Next
by $F+G= T^{\circ}\cap\ell_+ + S^{\circ}\cap\ell_+
\subset \ell_+$
we have $(T^{\circ}\cap\ell_+ + S^{\circ}\cap\ell_+)^{\circ}$
= $(T+\ell_+)\cap({S}+\ell_+)$ = $\ell_+$, hence
$F+G =\ell_+$.  So $F+S +G+T =H$.

Proof of (b).  Put $T= F^{\circ}\cap\ell_-$ and
$G= S^{\circ}\cap\ell_+$, then $(T+G)^{\circ}$
= $(F+\ell_-)\cap(S+\ell_+)$ = $F+S$. Now it is enough to show that
$G+F=\ell_+$. Since $(G+F)\subset\ell_+$, 
$\ell_+\subset (F+S^{\circ}\cap\ell_+)^{\circ}$ =$F^{\circ}\cap(S+\ell_+)$
=$F^{\circ}\cap S+\ell_+ = \ell_+$. Here we used the fact that
$\ell_+\subset F^{\circ}$ and $F^{\circ}\cap S =\{0\}$.  This proves
(b).
\end{proof}

\begin{cor}\label{cor:new polarization}
Let $\lla_-+\lla_+ =H$ be a polarization and $S+T =\lla_-$ 
be a decomposition by closed subspaces.  Put $F=T^{\circ}\cap\lla_+$
and $G=S^{\circ}\cap\lla_+$, then we have a new polarization of 
$H,~H= (T+F)+(S+G)$.
\end{cor}

\begin{cor}\label{cor:existence of inner product3}
Let $H =\ell_-+\ell_+$ be a polarized symplectic Hilbert space,
$S$ a closed subspace in $\ell_-$ and $F$ a closed subspace in
$\ell_+$. Assume that $S+F$ is a symplectic subspace as 
in the proposition above (b), then we can introduce
a compatible inner product with the symplectic form 
which satisfies that all for isotropic subspaces 
$S$, $F$, $T=F^{\circ}\cap\ell_-$
and  
$G=S^{\circ}\cap\ell_+$ are mutually orthogonal.
\end{cor}
We use this corollary in the proof
of Theorem \ref{thm:reduction2}.
\begin{rem}
The operation $S \mapsto S^{\circ}\cap\ell_+$
is idempotent. In fact,  $(S^{\circ}\cap\ell_+)^{\circ}\cap\ell_-$
= $(S +\ell_+)\cap\ell_- = S$.
\end{rem}

Let us consider the following situation:

{\it 
\begin{enumerate}   
\item[{\bf [CP1]}]\label{ass:polarizations}
Let $H$ and $B$ be symplectic Hilbert spaces with a compatible 
symplectic structure $~\omega_H$ and $\omega_B$, respectively.
We assume both are polarized with Lagrangian subspaces $\lla_{\pm}$
and $\ell_{\pm}$: 
$$
B=\lla_+ +\lla_-,
$$
$$
H=\ell_+ +\ell_-.
$$
\item[{\bf [CP2]}]\label{ass:injective maps}
There are continuous
injective maps
\begin{align*}    
&{\bf i}_+ :\ell_+ \to \lambda_+\\
&{\bf i}_- :\lambda_- \to \ell_- 
\end{align*} 
having ``{\em dense images}''.
\item[{\bf [CP3]}]\label{ass:injective and half symplectic}
For any $x\in\ell_+$ and $b\in \lambda_-$, 
$$
\omega_B({\bf i}_+(x),~b) = \omega_H(x ,~{\bf i}_-(b)).
$$
\end{enumerate}}

\begin{rem}
An example satisfying these conditions
is given in $\S 5.1$, Example \ref{example:5}.
\end{rem}

Let $\mu$ $\in {\cal F}\LLa(B)$ and put
$$
{\bg}(\mu)=
\{(x,~y)\in \ell_+\oplus \ell_-~|~\exists ~b\in 
\lla_-, ~({\bf i}_+(x),~b)\in\mu~\text{and}~y={\bf i}_-(b)\}.
$$

The subspace ${\bg}(\mu)$ is always isotropic, but need not be always
Lagrangian nor closed.  For example, 
if we take $\lla_-$ as $\mu$, then
${\bg}(\lla_-) = {\bf i}_-(\lla_-)$ is dense but not closed in
$\ell_-$, so it is not a
Lagrangian subspace. However if we restrict the map ${\bg}$
to  ${\cal F}\LLa_{\lla_-}(B)$, then we have
\begin{thm}\label{thm:reduction2}
\begin{enumerate} 
\item
For $\mu\in{\cal F}\LLa_{\lla_-}(B)$, ${\bg}(\mu)\in{\cal
  F}\LLa_{\ell_-}(H)$.
\item 
The map 
${\bg}:{\cal F}\LLa_{\lla_-}(B)\to{\cal F}\LLa_{\ell_-}(H)$
is continuous (and more strongly it is differentiable).
\item
Let $\{{\bf c}(t)\}$ be a continuous curve in 
${\cal F}\LLa_{\lla_-}(B)$,
then 
\begin{equation}\label{eq:invariance of Maslov}
\Mas(\{{\bf c}(t)\},~\lla_-) =\Mas(\{{\bg}({\bf c}(t))\},~\ell_-). 
\end{equation}
\end{enumerate}
\end{thm}
{\it We prove this theorem in the next subsection.} 
\bigskip

Let $\mathfrak{D}= {\bf i}_+(\ell_+)+\lambda_-$ and 
${\bf S}:\mathfrak{D}\to H,
~{\bf S}:{\bf i}_+(x)+b \mapsto x+{\bf i}_-(b)$, then
by the assumption {\bf [CP3]} the map ${\bf S}$
is a symplectic transformation and we have
\begin{prop}
Let $\mathbb{H} = H\boxplus B$ 
be the symplectic Hilbert space
with the symplectic form $\Omega= \omega_H - \omega_B$ and
let $C$ =
$$
 \{(x,~y)\boxplus(a,~b)\in\mathbb{H}~|~x\in \ell_+,~b\in
\lambda_-,~y={\bf i}_-(b), ~a={\bf i}_+(x)\}.
$$
Then $C$  is the graph of the map ${\bf S}$, 
is a Lagrangian subspace and ${\bf S}(\mu) = {\bg}(\mu)$.
\end{prop}

\begin{proof}It is easy to show that $C$ is isotropic.
So we only prove the following:
let $(u,~v)\boxplus(k,~k')\in \mathbb{H}$ satisfying
$$
\Omega((u,~u')\boxplus(k,~k'),~(x,~{\bf i}_-(b))
\boxplus({\bf i}_+(x),~b))=0,
$$
for any  $(x,~{\bf i}_-(b))\boxplus(~{\bf i}_+(x),~b)\in C$.

Then we have
$$\omega_H(u,~{\bf i}_-(b)) + \omega_H(u',~x) -\omega_B(k,~b)
-\omega_B(k',~{\bf i}_+(x))=0. 
$$
{}From this equation and Assumption {\bf [CP3]},
we have
$u'-~{\bf i}_-(k')=~0$ and 
$k-~{\bf i}_+(u) =0$, which show
$(u,~u')\boxplus(k,~k')\in C$ and $C$ is a Lagrangian subspace.

Now we see the coincidence ${\bf S}(\mu) = {\bg}(\mu)$ by the
definitions.  
\end{proof}

\subsection{Proof of Theorem \ref{thm:reduction2}}

Let $S$ be a finite dimensional subspace in $\lla_-$, 
then there is a
finite dimensional subspace $L$ in $\ell_+$ such 
that ${\bf i}_-(S)+ L$ 
is a symplectic subspace in $H$ and by Assumption {\bf [CP3]}  
$S+{\bf i}_+(L)$ is also a
symplectic subspace in $B$.  
There are many possibility to choose such a subspace $L$. 
We fix one of them by introducing a compatible inner product in $H$
with respect to which $\ell_{\pm}$ are orthogonal. So we put 
$L = J_H({\bf i}_-(S))$, 
where $J_H$ is the almost complex structure
defined by the compatible inner product. 



Hereafter we put the symplectic subspaces
$S+{\bf i}_+(L)$ = $B_S$ and ${\bf i}_-(S)+L$ = $H_S$
with $L = J_H({\bf i}_-(S))$, 
corresponding to a finite dimensional subspace $B$ in $\lambda_-$.

Then we have a symplectic isomorphism 
\begin{equation*}
\begin{array}{r@{\,}l}
{\bf i}^S={{\bf i}_-}_{|{B_S}}\oplus 
{({\bf i}_+)^{-1}}_{|{{\bf i}_{+}(L)}}: B_S=S+{\bf i}_+(L)
\stackrel{\sim}{\rightarrow} 
& H_S={\bf i}_-(S)+L\\
x+~y\qquad \mapsto& {\bf i}_-(x)+ ({\bf i}_+)^{-1}(y).
\end{array}
\end{equation*}

Next we remark that if $\theta$ 
is a Lagrangian subspace in $B_S$ then
$(S^{\perp}\cap\lla_-)+\theta$ is a Lagrangian subspace in $B$.
Let us denote the Lagrangian subspace of the form 
$(S^{\perp}\cap\lla_-) + \theta$ by $\lla(S,~\theta)$
and denote the subset of such Lagrangian subspaces that are
``transversal''
to $\lla_+$ by 
$$
\Lambda_S^{\,(0)}
=\{\lla(S,~\theta)=(S^{\perp}\cap\lambda_-) +\theta
\,~|\, ~\lla(S,~\theta)\cap\lambda_+=\{0\},~\theta\in\Lambda(B_S)\}.
$$
Then we have 
$$
\bigcup\limits_{B\subset\lambda_-,~\dim B<\infty}
\bigcup\limits_{\lla\in\Lambda_S^{\,(0)}}
{\cal F}\LLa_{\lla}^{(0)}(B)={\cal F}\LLa_{\lla_-}(B). 
$$


For $\mu\in {\cal F}\LLa_{\lla_-}(B)$ put $S=\lla_-\cap\mu$,
then by the way above 
there is a Lagrangian
subspace $\theta$ in $B_S$ such that
the Lagrangian subspace of the form 
$(S^{\perp}\cap\lla_-)+\theta$ is 
transversal both to 
$\mu$ and $\lla_+$ 
by Proposition \ref{prop:existence of transversal subspaces}.  
 
Let $\lla(S,\theta)\in\Lambda_S^{\,(0)}$, then we can
define new
polarizations of $H$ and $B$ which satisfy
Assumptions {\bf [CP1]}, {\bf [CP2]} and {\bf [CP3]}
by making use of $\lambda(S,\theta)$.
Theses are simply obtained by replacing the 
Lagrangian subspaces $\lla_-$ and $\ell_-$ by
$\lla(S,~\theta)$ and 
$\left(({\bf i}_-(S))^{\perp}\cap \ell_-\right)+ {\bf i}^S(\theta)$
respectively with the same $\lambda_+$ and $\ell_+$.
We repalce the map ${\bf i}_-$  
${{\bf i}_-}_{|{S^{\perp}\cap\lla_-}} \oplus {{\bf i}^S}_{|{\theta}}$.
Note that this change of polarizations do not change the map ${\bg}$. 

Now let $\lla\in\Lambda_S^{\,(0)}$, 
then any $\mu$ in ${\cal F}\LLa_{\lla}^{(0)}(B)$
is expressed as a graph of a continuous map $\phi: \lla_+ \to \lla$.
Hence on ${\cal F}\LLa_{\lla}^{(0)}(B)$ the map ${\bg}$ is expressed
in the form
$$
{\bg}(\mu)= ~\text{the graph of the map}~{\bf i}^S\circ \phi \circ {\bf i}_+.
$$
So it will be apparent of the continuity and also of 
the differentiability
of the map ${\bg}$ on $\mathcal{F}\Lambda_{\lambda}^{\,(0)}(B)$, 
if we know ${\bg}(\mu)\in {\cal F}\LLa_{\ell_-}(H)$.
 
To show the last assertion
it is enough to prove the case when $\lla$ is $\lla_-$ and also it
will be enough to prove $\gamma(\mu)^{\circ}\subset \gamma(\mu)$. 
So let $x\in \ell_+$ be an 
arbitrary element and assume 
that an element $a+b\in\ell_++\ell_-$ satisfies
$$ 
\omega_H(x+{\bf i}_-\circ \phi\circ{\bf i}_+(x),~a+b)=0.
$$
Then we have 
$$
\omega_H(x,~b)+\omega_H({\bf i}_-\circ \phi\circ{\bf i}_+(x),~a)=0.
$$
Hence by Assumption {\bf [CP3]}

\begin{align*}  
&\omega_H(x,~b)+\omega_H({\bf i}_-\circ \phi\circ{\bf i}_+(x),~a)\\
&=\omega_H(x,~b)+\omega_B(\phi\circ{\bf i}_+(x),~{\bf i}_+(a))\\
&=\omega_H(x,~b)+\omega_B(\phi\circ{\bf i}_+(x),~{\bf i}_+(a))
+\omega_B({\bf i}_+(x),~\phi\circ{\bf i}_+(a))\\
&\qquad\qquad\qquad\qquad\qquad\qquad\quad\quad
           -\omega_H(x,~{\bf i}_-\circ\phi\circ{\bf i}_+(a))\\
&=\omega_H(x,~b) +\omega_B({\bf i}_+(x)+\phi\circ{\bf i}_+(x),~{\bf
  i}_+(a)+\phi\circ{\bf i}_+(a))\\
&\qquad\qquad\qquad\qquad\qquad\qquad\quad\quad
           -\omega_H(x,~{\bf i}_-\circ\phi\circ{\bf i}_+(a))\\
&=\omega_H(x,~b)-\omega_H(x,~{\bf i}_-\circ\phi\circ{\bf i}_+(a))=0.
\end{align*}
{}From this we have
$$
b-{\bf i}_-\circ\phi\circ{\bf i}_+(a)=0,
$$
which shows that 
$(a,b)=(a,{\bf i}_-\circ\phi\circ{\bf i}_+(a))\in \gamma(\mu)$.
Hence ${\bg}(\mu)$ is a Lagrangian subspace.
These prove Theorem \ref{thm:reduction2} (a) and (b).

For the proof of  Theorem \ref{thm:reduction2} (c), 
we compute  Maslov indexes of particular paths
explicitly along the following steps (F1) $\sim$ (F4)
(see also similar arguments in \cite{CLM1}):

\begin{enumerate}
\item [(F1)] Let $\{{\bf c}(t)\}$ be a path such that all
${\bf c}(t)$ is transversal to $\lla_-$, 
then ${\bg}({\bf c}(t))$ is also transversal to
$\ell_-$ for any $t$. 
Hence $\Mas(\{{\bf c}(t)\},\lla_-)=\Mas(\{{\bg}({\bf c}(t))\},\ell_-)= 0$.
\item [(F2)] Let $L$ be a subspace in $\ell_+$ with 
$\dim L =1$ 
and let $\{{\bf c}(t)\}$ be a loop defined in 
Example \ref{ex:regular cross-path} of $\S 3.4$ 
for $F= {\bf i}_+(L)$. We assume the symplectic structure in $B$
is compatible.  So $F+J_B(F)$ is symplectic and 
$L+ {\bf i}_-\circ J_B(F)$ is also symplectic in $H$.  Let
$\tilde{F}=F^{\perp}\cap\lla_+$, 
then the path $\{{\bg}({\bf c}(t))\}$ is expressed
as
\begin{align*}
&{\bg}({\bf c}(t))\\
&=\{
\cos(\pi t) \cdot x 
+~\sin (\pi t) \cdot {\bf i}_-\circ J_B\circ {\bf
  i}_+(x)+ z~| ~x\in L,~z\in ({\bf i}_+)^{-1}(\tilde{F})\}.
\end{align*}
The path $\{{\bg}({\bf c}(t))\}$ has only one non-trivial
crossing at $t=1/2$.
We show that this is regular and determine 
the signature of the crossing form at $t=1/2$.

Let $A={\bf i}_-\circ J_B\circ {\bf i}_+$ and we take a suitable
subspace $\tilde{K}$
in $\ell_-$ such that $\tilde{K}\cap A(L)=\{0\}$, 
$\tilde{K}+({\bf i}_+)^{-1}(\tilde{F})$ is symplectic, and
$L+\tilde{K}=\nu$ is a Lagrangian subspace
transversal to $\{{\bg}({\bf c}(1/2))\}$, so that
we can express $\{{\bg}({\bf c}(t))\}_{|t-1/2|\ll 1}$ as graphs
of linear maps
\begin{equation*}
\begin{array}{r@{\,}l}
f_t:{\bg}({\bf c}(1/2))=A(L)&+({\bf i}_+)^{-1}(\tilde{F})  \to \nu\\
f_t:u\quad &+\quad z \qquad\mapsto\quad\cot(\pi t) \cdot A^{-1}(u).
\end{array}
\end{equation*}
Now we determine the crossing form at $t=1/2$.  Let $x,~y\in L$, then
\begin{align*}
&\frac{d}{dt}\omega_H(A(x),~f_t(A(y)))_{|t=1/2}\\
&=\omega_H\left({\bf i}_-\circ J_B\circ
{\bf i}_+(x),~\frac{-\pi}{\sin^2(\pi t)}y\right)_{|t=1/2}\\
&=-\pi \omega_B(J_B\circ {\bf i}_+(x),~{\bf i}_+(y))\\
&=\pi <{\bf i}_+(x),~{\bf i}_+(y)>^B.
\end{align*}
{}From this equality,
both of the Maslov indexes $\Mas(\{{\bf c}(t)\},~\lla_-) =1$ and
$\Mas(\{{\bg}({\bf c}(t))\},\ell_-) =1$.

\item [(F3)] Let $\mu$ in ${\cal F}\LLa_{\lla_-}(B)$ and
assume $\mu$ is contained in the Maslov cycle
${\mathfrak{M}}_{\lla_-}(B)$ with $\dim ~\mu\cap\lla_-$ = $N>0$.

Let $\{{\bf c}(t)\}$ 
be the path in ${\cal F}\LLa_{\lla_-}(B)$ 
defined in Example \ref{ex:regular cross-path} in $\S 3.4$
for $F={\bf i}_+\circ J_H\circ {\bf i}_-(\mu\cap\lla_-)$.

Note that here we used
Corollary \ref{cor:existence of inner product3}
for the existence of a compatible inner product in the symplectic
Hilbert space $B$. So by the corresponding almost complex structure $J_B$,
the isotropic subspaces $\mu\cap\lla_-$ is written as
$\mu\cap\lambda_- =
J_B({\bf i}_+\circ J_H\circ {\bf i}_-(\mu\cap\lla_-))$.
Then we can construct a path
in Example \ref{ex:regular cross-path} of $\S 3.4$ in terms of
unitary operators.

Now again the path $\{{\bg}({\bf c}(t))\}$ has only one non-trivial
crossing at $t=1/2$ with $\ell_-$ and by the same way 
as in (F2) we know that
the crossing form is positive definite
on $\mu\cap\lla_-$.

\item [(F4)] Finally we can prove the coincidence of the Maslov
indexes (\ref{eq:invariance of Maslov})
for arbitrary continuous paths. 
Since if the given path $\{{\bf c}(t)\}$ is a loop, 
then they coincide because of the fact that
they coincides for a generator of the fundamental group of the
space ${\cal F}\LLa_{\lla_-}(B)$. If it is not a loop, then  
by joining the paths in (F3) from the end point which is in 
the Maslov cycle ${\mathfrak{M}}_{\lla_-}(B)$ 
and we make this catenated path to a loop again 
by joining a path in (F1).
Now we know the Maslov indexes of these loops coincides and
Maslov indexes of added paths are all coincident, 
so that this prove Theorem \ref{thm:reduction2} (c).
\end{enumerate}

\section{Closed symmetric operators and Cauchy data spaces}

In this section we discuss Lagrangian subspaces in the
symplectic Hilbert space $\beta$ explained in 
Example (\ref{example:2}).

\subsection{Cauchy data space} 

Let $L$ be a real Hilbert space,  
and $A$ be a closed densely defined symmetric operator 
with the domain $\mathfrak{D}_m\subset L$. 
We denote the domain of the adjoint
operator
$A^*$ by $\mathfrak{D}_M$. As explained in the 
Example (\ref{example:2})
the factor space ${\bf{\beta}}=\mathfrak{D}_M/\mathfrak{D}_m$ is a
symplectic Hilbert space. 

Even if we add a bounded selfadjoint operator $B$ to the operator $A$,
we have the same domain $\mathfrak{D}_M$ 
of the adjoint operator $(A+B)^*$, the graph norms defined
on $\mathfrak{D}_M$ are equivalent and moreover the
symplectic forms defined by the operator $A$ and $A+B$
coincide. 

In any case we denote by $\gamma$ 
the natural projection map $\gamma: \mathfrak{D}_M \to {\bf{\beta}}$.
It will be
clear that $\gamma(\Ker A^*)$ is an isotropic subspace, but
we need some assumptions on the operator $A$ 
to show the closedness of it.  The space  
$\gamma(\Ker A^*)$ is called as the ``Cauchy data space''.  

\begin{prop}
Let $\mathfrak{D}$ be a subspace such that 
$\mathfrak{D}_m\subset \mathfrak{D}\subset\mathfrak{D}_M$.
Then the restriction
of the adjoint operator $A^*$ to the domain $\mathfrak{D}$
is selfadjoint, if and only if, the factor space $\gamma(\mathfrak{D})$
is a Lagrangian subspace in ${\bf{\beta}}$.
\end{prop}

{}From now on we assume that
\begin{enumerate}
\item [{[{\bf E1}]}:]
$A$ has at least one selfadjoint ``Fredholm'' extension, that is,
there exists a subspace $\mathfrak{D}$ (closed in the graph norm
topology) such that $A_{\mathfrak{D}}=A^*_{\,\,\,|\mathfrak{D}}$ is 
selfadjoint and has the finite dimensional kernel and the
image $A_{\mathfrak{D}}(\mathfrak{D})$ is closed in $L$ and is of
finite codimension.
\item [{[{\bf E2}]}:]
$\Ker A^*\cap \mathfrak{D}_m= \{0\}$.
\end{enumerate}

\begin{rem}
By the assumption [{\bf E2}], $A^*:\mathfrak{D}_M \to L$ 
is surjective. The condition [{\bf E2}]
requires that the unique continuation property holds,
for the case of differential operators. 
Both of these conditions [{\bf E1}] and [{\bf E2}] are satisfied by elliptic
operators of Dirac type on compact manifolds. For such operators
the unique continuation property holds for any hypersurfaces. 
The space $\mathfrak{D}_m$
will be the minimal domain of the definition, i.e.,
the subspace of the first order Sobolev space with the vanishing
boundary values, and the Cauchy data
space will be realized as a subspace in the distributions 
on the boundary manifold.
\end{rem}

\begin{prop}
Under the assumptions {{\bf [E1]}} and {\bf [E2]},
\begin{enumerate}
\item 
$\gamma(\Ker A^*)$ is a Lagrangian subspace,
\item
$\gamma(\Ker A^*)$ and $\gamma(\mathfrak{D})$ is a Fredholm pair.
\end{enumerate}
\end{prop}

\begin{proof}
Since $A^*(\mathfrak{D})$ is a closed finite codimensional subspace,
we have $(A^*)^{-1}(A^*(\mathfrak{D})) = \Ker A^* +\mathfrak{D}$
is a closed subspace in $\mathfrak{D}_M$(equipped with the graph norm
topology), hence $\gamma(\Ker A^* +\mathfrak{D})=\gamma(\Ker A^*) 
+\gamma(\mathfrak{D})$ is closed in ${\bf{\beta}}$.

Again since $\Ker A^* +\mathfrak{D}$ is closed and
$\dim (\Ker A^*\cap\mathfrak{D})< +\infty$ we know
that $\Ker A^* +\mathfrak{D}_m$ is also close in $\mathfrak{D}_M$, and
so $\gamma(\Ker A^*)=\gamma(\Ker A^* +\mathfrak{D}_m)$ 
must be close in  ${\bf{\beta}}$, and that it is a closed isotropic subspace.

Now we have relations:
\begin{enumerate}
\item
$\dim \gamma(\Ker ^*)\cap\gamma(\mathfrak{D})= \dim \Ker
A^*\cap\mathfrak{D}$,
\item
$\dim L/(\Ker A^*+ \mathfrak{D}) 
= \dim\Ker A^* \cap\mathfrak{D}$.
\end{enumerate}
So we have
$$
\dim \gamma(\Ker A^*)\cap\gamma(\mathfrak{D})=
\dim \gamma(\Ker A^*)^{\circ}\cap\gamma(\mathfrak{D})=
\dim {\bf{\beta}}/(\gamma(\Ker A^*)+\gamma(\mathfrak{D}))<\infty,
$$
and hence 
$$\gamma(\Ker A^*)\cap\mathfrak{D}
=\gamma(\Ker A^*)^{\circ}\cap\mathfrak{D},
$$
$$
\gamma(\Ker A^*)+\mathfrak{D}
=\gamma(\Ker A^*)^{\circ}+\mathfrak{D}.
$$
{}From these equalities $\gamma(\Ker A^*)$ is a Lagrangian subspace
and $\gamma(\Ker A^*)\in
\mathcal{F}\Lambda_{\gamma(\mathfrak{D})}({\bf{\beta}})$.
\end{proof}

\begin{cor}\label{soft}
Under the assumptions {\bf [E1]} and {\bf [E2]} the extension of $A$
on $\mathfrak{D}_m + \Ker A^*$,
$A^*_{|\mathfrak{D}_m + \Ker A^*}$,
is a selfadjoint operator.  
\end{cor}

\begin{rem}
The extension in the above corollary (\ref{soft}) is called 
``Soft extension''. This is also an interesting extension, although it
is far from Fredholm operators. For example,
in the paper \cite{Gu} the asymptotic
behavior of non-zero eigenvalues was investigated
for a symmetric elliptic operator of even order on a bounded domain.
\end{rem}

\begin{example}\label{example:4}
Let 
\begin{equation}\label{ordinary-differential-operator}
A=\mathbb{J}\left(\frac{{d}}{{d t}}+B\right)
\end{equation} 
be an ordinary differential operator acting 
on $C_0^{\infty}(0,1)\otimes\mathbb{R}^{2N}$, where
$\mathbb{J}=
\left(\begin{array}{cc}
0& I\\[0.1cm]      
-I&0
\end{array}\right)
$
and $B$ is a $2N\times 2N$ symmetric matrics.
In this case $\beta=\mathfrak{D}_M/\mathfrak{D}_m \cong
\mathbb{R}^{2N}\boxplus\mathbb{R}^{2N}$ with the symplectic
form $\mathbb{J}\boxplus-\mathbb{J}$.
The cases treated in \cite{Fl} reduce to this case (see also \cite {OF}).
\end{example}

\begin{example}\label{example:5}
We describe an example of the Cauchy data space which can be
realized in the distribution space on a manifold.

Let $M$ be a manifold with boundary $\Sigma=\partial M$, and
$A$ be a first order symmetric elliptic operator 
acting on a space of  smooth 
sections of a real smooth vector bundle $E$.
Here we mean that the operator is symmetric, when it is
symmetric on the space of smooth sections whose supports
do not intersect with the boundary manifold $\Sigma$.

We assume that the unique continuation property holds
for this operator $A$ with respect to 
the boundary hypersurface $\Sigma$.

The minimal domain of the definition $\mathfrak{D}_m$
is the subspace of the first order Sobolev space
consisting of sections with vanishing boundary values. 
Even for this case, it is not easy to determine the
domain of the adjoint operator. It is a little bit
bigger than the whole first order Sobolev space.
The Cauchy data space $\beta$ is included in the
Sobolev space of order $-1/2$ on $\Sigma$ (\cite{Ho1}).

Then we assume that $A$ has a product form near the boundary 
hypersurface $\Sigma$ in the following sense:

Let $\mathcal{N}\cong [0,1]\times \Sigma$ is a neighborhood
of $\Sigma$ and on this neighborhood, the operator 
$A$ takes the form
$$
A= \sigma\left(\frac{\partial}{\partial t}+ B\right),
$$
where $\sigma$ is a bundle automorphism of the restriction
of $E$ to $\mathcal{N}$, and is independent from the coordinate of
the normal direction $t\in [0,1]$. It is also  
skew symmetric and satisfies $\sigma^2=-Id$.  The operator $B$
is selfadjoint, first order elliptic operator on the vector
bundle $E_{\,|\Sigma}$, also independent from the normal 
variable $t$ and satisfies the relation  $\sigma\circ B+ B\circ
\sigma=0$ by the symmetric assumption. Now we can characterize
the Cauchy data space in the following form:

Let $\{e_k\}_{k\in\mathbb{Z}\backslash\{0\}}$, 
$e_k>0$ for $k>N_0$, $e_k<0$ for $k<-N_0$ and $e_k=0$ for $|k|\leq N_0,
k\not=0$, 
be the eigenvalues of the boundary
operator $B$ and denote by $\{\phi_k\}_{k\in\mathbb{Z}\backslash\{0\}}$ 
the corresponding
orthonormal eigensections. Then we define the spaces by
$$
H^+ =\{\sum\limits_{k<0~\text{finite sum}}c_j\phi_j \},
$$
$$
H^- =\{\sum\limits_{k>0~\text{finite sum}}c_j\phi_j \} 
$$
Let $\lambda_{\pm}$ be the closure of $H^{\pm}$ with respect to
the $\pm 1/2$ order Sobolev norm respectively, then the direct sum
$\lambda_+ +\lambda_- =\beta$ 
(\cite{Ho1}, \cite{APS}).

Then let $L_{\pm}$ be the closures of $H^{\pm}$ with respect to the
$L_2$-norm, then we have two symplectic Hilbert spaces
$L_2(M)$ and $\beta$ satisfying the conditions {\bf [CP1]}, {\bf
  [CP2]} and {\bf [CP3]} in $\S 4.2$.
\end{example}

\begin{rem}
In the above example (\ref{example:5}), 
if the boundary of the manifold $M$
is divided into two components $\Sigma_0$ and $\Sigma_1$, 
then the space of boundary values $\beta$ is also divided
into the sum $\beta=\beta_0 \oplus \beta_1$, where $\beta_i$
is in the Sobolev space of order $-1/2$ on $\Sigma_i$ $(i=0,1)$.
Now the Cauchy data space $\gamma(\Ker A^*)$ 
defines a closed symplectic transformation $S: \mathfrak{D}\to \beta_1$,
where $\mathfrak{D}= \{x\in\beta_0\,|\, \exists ~y\in\beta_1,~(x,y)\in
\gamma(\Ker A^*)\}$ $S(x)=y,~(x,y)\in \gamma(\Ker A^*)$. We should
note that this follows from the unique continuation property, i.e.,
$(\beta_0\oplus\{0\})\cap \gamma(\Ker A^*)=\{0\}$, and
$(\{0\}\oplus\beta_1)\cap \gamma(\Ker A^*)=\{0\}$.
Also a selfadjoint Fredholm extension is given by the
Atiyah-Patodi-Singer boundary condition for the case of the operators
with product form near the boundary.  Or more generally, even if it is
not of a product form near the boundary, there are such extensions
by global elliptic boundary conditions(see for example \cite{Ra}).
\end{rem}

\subsection{Continuity of Cauchy data spaces} 

Let $A$ be the symmetric operator as above satisfying the
conditions {\bf [E1]} and {\bf [E2]}.
Let $\{B_t\}_{t\in [0,1]}$ be a continuous
family of bounded selfadjoint operators on the Hilbert space $L$.
If each operator $A+B_t$ for $t\in [0,1]$ satisfies
the conditions {\bf [E1]} ``with a common domain'' $\mathfrak{D}$,
i.e., $(A+B_t)^*_{\,\,\,|\mathfrak{D}}=A_{\mathfrak{D}}+B_t$ ($t\in[0,1]$) 
is selfadjoint and a Fredholm operator,
and {\bf [E2]}, then we have a family
of Lagrangian subspaces $\{\gamma(\Ker (A+B_t)^*)\}_{t\in [0,1]}$
in ${\bf{\beta}}$ and each of them and
$\gamma(\mathfrak{D})$ is a Fredholm pair.

\begin{prop}\label{continuity-of-Cauchy-data-space}
The family $\{\gamma(\Ker (A+B_t)^*)\}_{t\in [0,1]}$
is a continuous family.  Hence it is a continuous
path in $\mathcal{F}\Lambda_{\gamma(\mathfrak{D})}({\bf{\beta}})$.
\end{prop}

\begin{proof}
It will be enough to prove at $t=0$. So let 
$T_t:\mathfrak{D}_M \to L\oplus \Ker(A+B_0)^*$ be a map
defined by $T_t(x)= (A+B_0)^*(x)\oplus \pi_0(x)$, where $\pi_0$
is a projection operator $\pi_0:\mathfrak{D}_M \to \Ker(A+B_0)^*$.
Then since $T_0$ is an isomorphism, for a sufficiently small
$\epsilon > 0$ the maps $T_t$ for $0\leq t\leq \epsilon$ 
are also isomorphisms. Hence we have 
$\Ker (A+B_t)^* = (T_t)^{-1}(\{0\}\oplus \Ker(A+B_0)^*)$,
and that the family $\{\Ker (A+B_t)^*\}_{0\leq t\leq \epsilon}$
is continuous at $t=0$ since the family 
$\{(T_t)^{-1}\}_{0\leq t\leq \epsilon}$ is a continuous family.
\end{proof}

\begin{rem}
If the operator $A^*_{\,\,\,|\mathfrak{D}}=A_{\mathfrak{D}}$ has 
a compact resolvent and a Fredholm operator, 
then for any selfadjoint bounded operator $B$
the sum $A+B$ satisfies the condition {\bf [E1]}.
\end{rem}

\begin{rem}
In the case of the paper \cite{Fl}, the family (= the family of operators
of the form (\ref{ordinary-differential-operator})) 
has varying domains 
$\{\mathfrak{D}\}$ 
where the operator is realized as a selfadjoint
operator according to the each value of the parameter. But
in this case the operator family can be transformed into the above case
of a fixed domain of the definition for the selfadjoint realization
by a continuous family of unitary operators. It would
not be clear whether we can do such transformations
for the family of elliptic operators 
in the higher dimensions (Example (\ref{example:5})).
\end{rem}

\subsection{Spectral flow and Maslov index}     

Finally we just formulate an equality between ``Spectral flow'' and 
``Maslov index'' arising from the family 
of operators explained 
in the previous subsection.

Let $\mathcal{F}(L)$ be the space of bounded
Fredholm operators defined on a Hilbert space $H$. 
It is a classifying space for the K-group.
Then the non-trivial component of the 
subspace $\widehat{\mathcal{F}}(L)$ consisting of selfadjoint Fredholm
operators, we denote it by $\widehat{\mathcal{F}}_*(L)$, 
is a classifying space for the $K^{-1}$-group (in the complex case) 
and $K^{-7}$-group(in the real case).  Both of their fundamental groups
are isomorphic to $\mathbb{Z}$(\cite{AS}).  These isomorphisms are 
given by an integer, so called, the spectral flow 
for a family of selfadjoint Fredholm operators(\cite{APS}).  
This integer is also defined 
for continuous path of selfadjoint Fredholm operators
without any assumptions at the end points(\cite{Ph}).
We do not state the definition, but is given in a similar way as
the definition of the Maslov index we gave in this article, 
or rather it should be thought of the initiating method 
which was given in the paper \cite{Ph} to define the spectral flow
based on the basic spectral property 
of the Fredholmness of the operators.  

{\it Let $L$, $A$ and $\{B_t\}_{t\in [0,1]}$ be as above, that is,
the family $\{A+B_t\}$ acting on the Hilbert space $L$ 
satisfys the conditions {\bf [E1]} with a common
subspace $\mathfrak{D}$ 
on which the operators $(A+B_t)^*_{\,\,\,|\mathfrak{D}}$ 
are selfadjoint and Fredholm.
Then we see that $(A+B_t +s)^*_{\,\,\,|\mathfrak{D}}=A_{\mathfrak{D}}+B_t+s$ 
is also a Fredholm operator for
sufficiently small $|s|<< 1$. Now instead of the condition 
{\bf [E2]} we assume a stronger property:

{\bf [E2']} There exists an $\epsilon > 0$ such that for each $t\in[0,1]$
and $|s|<\epsilon$,
$\Ker (A+B_t + s)^*\bigcap\mathfrak{D}_m=\{0\}$.}

\begin{rem}
Of course this condition is satisfied by Dirac type operators.
\end{rem}

Under these assumptions {\bf [E1]} and {\bf [E2']}, and
with the common domain of the definition $\mathfrak{D}$
for the selfadjoint Fredholm realization, we have
\begin{thm}
The spectral flow for the family
$\{A_{\mathfrak{D}}+B_t\}_{t\in[0,1]}$
and the Maslov index of the path
of Cauchy data spaces $\{\gamma(\Ker(A+B_t)^*)\}_{t\in[0,1]}$
with respect to the Maslov cycle $\gamma(\mathfrak{D})$
coincides.
\end{thm}

We do not give a proof of this theorem. First it was proved in
\cite{Fl} that a coincidence between ``Spectral flow'' and
``Maslov index of boundary data'' for a family of ordinary differential
operators(Example (\ref{example:4})). 
In this case the family of ordinary differential operators 
arises as the family of the Euler equations
of the symplectic action integral which is defined by two
transversally intersecting Lagrangian submanifolds
in a symplectic manifold and the Maslov index in this case is of 
the finite dimension(see also \cite{OF}). Then
it was proved
in \cite{Yo} on three dimensional manifolds and generalized to
higher dimensions in \cite{Ni} for a family 
of Dirac operators $\{A_t\}_{t\in [0,1]}$ 
with invertible operators at the end points $t=0,1$. 
In these cases the Maslov indexes are that in the infinite dimension.
We reproved
the theorem in the above general form in \cite{FO2}. There we
also proved a general addition formula for the spectral flow
when we decompose a manifold into two parts. To prove it we apply our
reduction theorem in $\S 3$ of the Maslov index in the infinite
dimensions.  Such types of the formula
were also investigated in several authors or believed
to hold in a more general contexts(\cite{Ta}, \cite{DK}, \cite{CLM2}).
We tried to make clear the meaning of the condition 
that the operators in the family are of the form, 
so called, product form near the separating boundary
manifold (Example (\ref{example:5})). This kind of restriction for the family
will correspond to a condition assumed in the
Mayer-Vietoris exact sequence of the singular homology theory.








                \appendix
              \setcounter{secnumdepth}{1}
                \setcounter{thm}{0}
                \setcounter{section}{0}
                \setcounter{equation}{0}

\section*{Appendix}

In this appendix  we gather up some of fundamental facts without
proofs from the theory of 
Functional Analysis, on which our arguments heavily rely.
Because the objects we will deal with are infinite dimensional spaces and
their homotopical properties. 


Our Hilbert spaces will
be mostly real separable Hilbert spaces and the theorems
we sum up here are valid for both real and complex cases if we do not
state particularly.  So let $H$
be a separable Hilbert space with the
inner product by $<\bullet\,,\,\bullet>$ and
as usual we denote the norm of the element $x \in H$ by 
$\Vert x\Vert = \sqrt{<x\,,\,x>}$.

\section{Topology of operator spaces}
\begin{thm}[Kuiper's Theorem]\label{thm:Kuiper}
Let $H$ be an infinite dimensional
real (complex or quaternionic) separable Hilbert space,
then
the group of linear isomorphisms, we denote it by $GL(H)$,
is contractible to a point. Note that the topology of $GL(H)$
is defined by the norm convergence and it is a topological group with
this topology.
\end{thm}

\begin{cor}
Let $H$ be an infinite dimensional
real (complex or quaternionic) separable Hilbert space,
then the subgroup of $GL(H)$ consisting of linear isomorphisms
which preserves the inner product is also contractible to a point.
We will denote them by $\cal{O}(H)$ (orthogonal group) for the real
case, 
$\cal{U}(H)$ (unitary group) for the
complex case and $Sp(H)$ (symplectic group) for the quaternionic case.
\end{cor}

\begin{rem}
In the real case, the groups $GL(H)$ and $\cal{O}(H)$ are connected
(path-wise) unlike the finite dimensional case.
\end{rem}

\begin{thm}[Palais's Theorem]\label{thm:Palais}
Let $B$ be a Banach space and we assume there is a sequence
of projection operators $\{\pi_n\}_{n=1}^{\infty}$ 
onto finite dimensional subspaces $L_n =\pi_n(B)$
such that $L_n \subset L_{n+1}$ and for each $x\in B$,$\{\pi_n(x)\}$
converges to $x$ in the sense of norm, that is, 
$\{\pi_n\}_{n=1}^{\infty}$ 
converges to the identity operator in the strong sense.
Then for each  open set $O\subset B$, 
the injection map $j:ind\lim\limits_{\rightarrow} \pi_n(O)
\rightarrow O$ is a homotopy equivalence.
\end{thm}

Let $H$ be a real (or complex) Hilbert space, and by fixing
a complete orthonormal basis $\{x_n\}_{n=1}^{\infty}$, we have
inclusions of finite dimensional subspaces $E_n$, where $E_n$ is
spanned by $\{x_i\}_{i=1}^{n}$. Also from these inclusions of
subspaces
we have inclusions of the general linear groups $\GL(n,\mathbb{R})$
(or $\GL(n,\mathbb{C})$):
\[\GL(n,\mathbb{R}) \subset   \GL(n+1,\mathbb{R}) 
\]and 
\[\GL(n,\mathbb{C}) \subset   \GL(n+1,\mathbb{C}) 
\]
in an obvious way.  Then we have also inclusions 
$\GL(n,\mathbb{R}) \rightarrow \GL_K(H)$ 
(in the complex case $\GL(n,\mathbb{C}) \rightarrow \GL_K(H)$),
where we denote 
\[\GL_K(H) =\{g\in \GL(H)\,|\, g ~\text{is of
  the form} ~Id + ~\text{compact operator}\}
\]
corresponding to the each case. 

\begin{prop}
The inclusion maps
\[
j:ind\lim\limits_{\rightarrow} \GL(n,\mathbb{R}) \rightarrow \GL_K(H) 
\] for the real case and
\[
j:ind\lim\limits_{\rightarrow} \GL(n,\mathbb{C})\rightarrow \GL_K(H)
\] for the complex case,
are homotopy equivalences.
\end{prop}

\section{Spectral notions}
Let $A$ be a densely defined closed operator (bounded or not bounded) 
on a Hilbert space $H$.  Let $\lambda \in\mathbb{C}$, then
$\lambda$ is called a resolvent of the operator $A$, if $A-\lambda$
has a bounded inverse defined on the whole space $H$. We denote
the set of all resolvents by $\rho(A)$.  The complement
$\mathbb{C}\,\backslash~\rho(A)$ is called spectrum of $A$ and we denote
it by $\sigma(A)$. Let $\lambda \in\sigma(A)$, then if $A-\lambda$
has a densely defined inverse, but not continuous, then $\lambda$
is called a continuous spectrum and we denote the subset consisting
of continuous spectra by ${\bf C}_{\sigma}(A)$. Again let $\lambda$ be in
$\sigma(A)$
and assume that $A-\lambda$ is not invertible, that is $\{ x\in H\,|\,
(A-\lambda)(x)=0\}\not= 0$, then such $\lambda$ is called an
eigenvalue or a point spectrum.
We denote the set of eigenvalues by ${\bf P}_{\sigma}(A)$.  The element in
the complement in $\sigma(A)$ of the union 
${\bf P}_{\sigma}(A)\cup {\bf C}_{\sigma}(A)$
is called residual spectrum, and we denote them by
${\bf R}_{\sigma}(A)$.  Let $\lambda\in 
{\bf R}_{\sigma}(A)$, then
$A-\lambda$ has an inverse, but the image ${\it Im}(A-\lambda)$ is not
dense.

Now let $A$ be a selfadjoint operator (bounded or not bounded), then
we know that there are no residual spectrum of $A$, that is, the
spectrum $\sigma(A) = {\bf C}_{\sigma}(A)\cup{\bf P}_{\sigma}(A)$ and
$\sigma(A)\subset \mathbb{R}$.
  
We denote by $\sigma_{ess}(A)$ a subset of $\sigma(A)$, each of which
element
is called ``essential spectrum'', 
if $\lambda$ is an eigenvalue of infinite multiplicity or a continuous
spectrum.
If  $A$ is
bounded selfadjoint, then $\sigma(A)$ is compact and
$\norm{A} = \sup\limits_{t\in\sigma(A)}|t|$.

Let $\{E_t\}_{\{t\in \mathbb{R}\}}$ be a family of orthogonal projections
defined on a Hilbert space $H$ satisfying following properties 
({\it Sp 1}), ({\it Sp 2}), ({\it Sp 3}) and ({\it Sp 4}), 
then we call
$\{E_t\}_{\{t\in \mathbb{R}\}}$ a spectral measure:
\begin{align*}
(Sp \,\,1)&\quad  E_t(H) \subset E_s(H)\quad\text{for}\quad t\leq s.\\
(Sp \,\,2)&\quad  E_t \quad\text{is right strong
  continuous},\,\text{that is},\,
\text{for each}\quad x\in H,\\
&\quad \lim\limits_{0<\delta\rightarrow 0} E_{t+\delta}(x) = E_t(x)\\
(Sp \,\,3)&\quad \lim\limits_{t\rightarrow\infty} E_t(x) =x\quad\text{for
  each}\quad x\in H\\
(Sp \,\,4)&\quad \lim\limits_{t\rightarrow-\infty} E_t(x) =0\quad\text{for
  each}\quad x\in H.
\end{align*}

\begin{thm}[Spectral Decomposition Theorem]\label{thm:spec decompo}
Let $A$ be a selfadjoint operator (bounded or not bounded) 
defined on a Hilbert space $H$. Then
there is a unique spectral measure $\{E_t\}_{\{t\in\mathbb{R}\}}$
such that 
\[
A= \int\limits_{-\infty}^{\infty}\,t\, dE_t.
\]
\end{thm}

\begin{rem}
The domain $\mathfrak{D}$ of the operator $A$ is characterized as 
\[\mathfrak{D}=\{x\in H \,|\, \int\limits_{-\infty}^{+\infty}|t|d
\norm{E_t(x)}^2 < \infty\}. 
\]
\end{rem}


\section{Fredholm operators}
Let $H$ be a Hilbert space (or Banach space) and
let $T$ be a densely defined closed operator with the domain
$\mathfrak{D}$.
We call a closed operator $T$ is  a Fredholm operator, if it satisfies
\begin{align*}
&\dim Ker (T) \quad\text{is finite},\\
&\text{the image}\quad \it{Im}~(T)=T(\mathfrak{D})\quad\text{is closed},\\
&\dim Coker~(T)=H/ \it{Im}~(T)\quad\text{is finite}.
\end{align*}

\begin{rem}
For bounded Fredholm operators $T$ we can prove that the image $T(H)$
is closed from the finite codimensionality of it by making use 
of the open mapping theorem.
\end{rem}

Let $\mathcal{F}(H)$ be the space of all ``{\it bounded}'' Fredholm operators
defined on a Hilbert space $H$. 

\begin{prop}\label{prop:open}
The space $\cal{F}(H)$ is an open
subset in the space of all bounded operators $\cal{B}(H)$ with the
topology of the norm convergence. 
\end{prop}

Let $\cal{B}(H)$ be the algebra of all bounded operators on $H$
and $\cal{K}(H)$ the two-sided ideal of compact operators, then
the quotient algebra $\cal{B}(H)/\cal{K}(H)$ is called Calkin algebra.
If $\pi$ denotes the natural projection map 
$\pi:\cal{B}(H)\rightarrow \cal{B}(H)/\cal{K}(H)$, then

\begin{prop}\label{prop:Calkin}
$\pi^{-1}\left((\cal{B}(H)/\cal{K}(H))^*\right) =\cal{F}(H),$
where $(\cal{B}(H)/\cal{K}(H))^*$ denotes the group consisting 
of the invertible elements in the Calkin algebra $\cal{B}(H)/\cal{K}(H)$.
\end{prop}

For a Fredholm operator (closed or bounded) $T$ 
we denote
by $ind ~(T)$ the difference 
\[
ind ~(T)=\dim \Ker T - \dim \Coker T,
\] 
and call it the "\textit{Fredholm index}" of the operator.  Especially
for  bounded Fredholm operators $T\in\cal{F}(H)$ we have 

\begin{thm}\label{thm:Fredholm index}
\[
ind:\cal{F}(H)\rightarrow \mathbb{Z}
\]
is a locally constant function, in fact, it distinguishes the
connected components (= path-wise connected) of the space $\cal{F}(H)$.  
\end{thm}

\begin{rem}
\begin{enumerate}
\item 
If $H$ is finite dimensional, then the quantity $ind~(T)$ always
vanishes. So this has an only meaning in the infinite dimension.
\item
In the paper \cite{CL} a similar result for the connected components
is proved for
the space of all closed Fredholm operators.  The topology
for such a space is introduced by embedding it into the space
of bounded operators on the product space $H\times H$ as 
orthogonal projection operators onto the graphs.
\end{enumerate}
\end{rem}

\begin{thm}\label{thm:invariance under compact perturbation}
Let $K$ be a compact operator on $H$ and $T$ be a bounded Fredholm
operator, then $T+K$ is also a Fredholm operator and
\[
ind ~(T+K) = ind ~T.
\]
\end{thm}

\section{Existence of a compatible symplectic structure}

\begin{prop}\label{prop:appen1}
Let $(H,~(\bullet,~\bullet))$ be a real Hilbert space and $\omega$ a bounded and 
non-degenerate skew symmetric bilinear form on $H$. 
Then we can replace the inner product by another one $<\bullet,~\bullet>$
such that $(H,~\omega,~<\bullet,~\bullet>, J)$ is a compatible symplectic
Hilbert space.
\end{prop}

\begin{proof}  
Let $A$ be the operator defined by 
$$
\omega(x,y)= \left(A(x),~y\right).
$$
Then $A$ is bounded, skew-symmetric and invertible.  
Put $|A| = \sqrt{{\,^tA}\circ A}$,
and the new inner product by $<x,~y>= (|A|(x),~y)$. 
By this inner product
we can express $\omega(x,~y)= <J(x),~y>$, where $J=|A|^{-1}\circ A$.
Now $J^2 =|A|^{-1}\circ A\circ |A|^{-1}\circ A 
=|A|^{-2}\circ A^2 =-Id$, and also 
$<J(x),~J(y)>$ = $(|A|\circ |A|^{-1}\circ A(x),~|A|^{-1}\circ A(y))$ =
$(|A|^{-1}{\,^tA}\circ A(x),~y)$ = $(|A|(x),~y)=<x,~y>$. 
\end{proof}

\begin{cor}\label{cor:orhogonal polarization}
Let $H$ be a symplectic Hilbert space and we assume that $H$ is
polarized by two Lagrangian subspaces $\lla$  nd $\mu$ :
$H=\lla\oplus\mu$.
Then there is an inner product with which the symplectic structure is
compatible
and the decomposition is orthogonal.
\end{cor}
\begin{proof}
In the above proof we can assume that 
the subspaces $\lla$ and $\mu$ are orthogonal
with  respect to the inner product $(\bullet,\bullet)$. 
Then the operator
$A$, $(A(x),~y)=\omega(x,~y)$
maps $\lla$ to $\mu$ and $\mu$ to $\lla$. Hence ${\,^tA}\circ A$ and
its square root keep the subspaces $\lla$ and $\mu$ invariantly, 
so that 
with the new inner product $(|A|(x),~y)$ 
the Lagrangian subspaces $\lla$ and $\mu$ are orthogonal.
Now by the same
way as Proposition \ref{prop:appen1} 
the new inner product $(|A|(x),~y)$ gives us the compatible symplectic
structure. 
\end{proof}





\begin{thebibliography}{CLM1}

\bibitem[Ar]{Ar} {\it V. I. Arnold}, 1967, 
Characteristic class entering in quantization
conditions, {\it Functional Anal. Appl. {\bf 1}}(9167), 1-13.

\bibitem[APS]{APS}
{\sc M. F. Atiyah, V.K. Patodi and I. M. Singer},
Spectral asymmetry and Riemannian geometry: 
I,  {\it Math. Proc. Camb. Phil. Soc., {\bf 77}}(1975), 43-69.
II, {\it Math. Proc. Camb. Phil. Soc., {\bf 78}}(1975), 405-432.
III, {\it Math. Proc. Camb. Phil. Soc., {\bf 79}}(1976), 71-99.


\bibitem[AS]{AS}
{\sc M. F. Atiyah and I. M. Singer},
Index theory for skew-adjoint Fredholm operators, {\it Inst.
Hautes \'Etudes Sci. Publ. Math. {\bf 37}} (1969), 5-26.

\bibitem[BF1]{BF1} {\sc B. Boo\ss--Bavnbek and K. Furutani}, The
Maslov index -- a functional analytical definition and the
spectral flow formula, {\it Tokyo J. Math. {\bf 21}} (1998),
1--34.

\bibitem[BF2]{BF2} \bysame, \bysame, Symplectic functional analysis
and spectral invariants, in: B. Booss-Bavnbek, K.P.
Wojciechowski (eds.), ``Geometric Aspects of Partial
Differential Equations'', Amer. Math. Soc. Series {\it
Contemporary Mathematics}, vol. 242, Providence, R.I., 1999,
pp. 53--83.

\bibitem[BFO]{BFO}
{\sc B. Booss-Bavnbek, K. Furutani, and N. Otsuki},
Criss--cross reduction of the Maslov index and a proof of the
Yoshida--Nicolaescu Theorem, {\it Tokyo J. Math. 
vol. {\bf  24}, No. 1(2001)}

\bibitem[CL]{CL} {\sc H. O. Cordes and J. P. Labrousse}, 
The invariance of the index in the metric of closed operators, 
{\it J. Math. and Mechanics {\bf 12}}(1963), 693-719.

\bibitem[CLM1]{CLM1} {\sc S.E. Cappell, R. Lee, and E.Y. Miller},
On the Maslov index, {\it Comm. Pure Appl. Math. {\bf 47}}
(1994),  121--186.

\bibitem[CLM2]{CLM2} {\sc ---}, {\sc ---}, {\sc ---}, Selfadjoint
elliptic operators and manifold decompositions Part I: Low
 eigenmodes and stretching, {\it Comm. Pure Appl. Math. {\bf
 49}} (1996), 825--866. Part II: Spectral flow and Maslov index,
{\it Comm. Pure Appl. Math. {\bf 49}} (1996), 869--909. Part
III: Determinant line bundles and Lagrangian intersection, {\it
Comm. Pure Appl. Math. {\bf 52}} (1999), 543--611.


\bibitem[CP]{CP} 
{\sc A. Carey AND J. Phillips}, Unbounded Fredholm Modules and Spectral Flow,
{\it Canad. J. Math. {\bf 50}}(1998), no. 4, 673-718.




\bibitem[DK]{DK} {\sc M. Daniel and P. Kirk, with an appendix
by K.P. Wojciechowski}, A general splitting formula for the
spectral flow, {\it Michigan Math. Journal {\bf 46}}(1999),
589--617.

\bibitem[Du]{Du}{\sc J.J. Duistermaat}, On the Morse index in
variational calculus, {\it Adv. in Math. {\bf 21}} (1976), 173--195.


\bibitem[Fl]{Fl}{\sc A. Floer}, A relative Morse index for the
  symplectic action, {\it Comm. Pure Appl. Math. {\bf 41}}(1988), 393-407.

\bibitem[FO1]{FO1}{\sc K. Furutani and N. Otsuki}, 
Spectral flow and
intersection numbers, {\it J. Math. Kyoto Univ. {\bf 33}}(1993), 261-283.

\bibitem[FO2]{FO2}{\sc K. Furutani and N. Otsuki}, 
Maslov index in the infinite dimension and 
a splitting formula for a spectral flow,
{\it Japanese Journal of Mathathematics. Vol. 28, No. 2}(2002), 215-243\\

\bibitem[Fu]{Fu}{\sc K. Furutani},
On the Quillen determinant,
{\it to appear in Journal of Geometry and Physics}. 

\bibitem[Ge]{Ge}{\sc E. Getzler}, The odd Chern character in cyclic
homology and spectral flow, {\it Topology {\bf 32}}(1993), 489-507.

\bibitem[Go1]{Go1}{\sc M. de Gosson}, La d{\'e}finition de l'indice de Maslov
sans hypoth{\`e}se de transversalit{\'e},
{\it C. R. Acad. Sci. Paris {\bf 310}, S\'erie I}
(1990), 279--282.

\bibitem[Go2]{Go2}{\sc ---}, The structure of
$q$--symplectic geometry, {\it J. Math. Pures Appl. {\bf 71}}
(1992), 429--453.



\bibitem[Gu]{Gu}{\sc G. Grubb}, Spectral asymptotics for the
``soft'' selfadjoint extension of asymmetric elliptic
differential operator, {\it Journal of Operator Thoery,
  Vol. 10},(1983), 9-20.

\bibitem[Ho1]{Ho1}{\sc L. Hoermander},
Pseudo-differential operators and non-elliptic boundary
problems, {\it Ann. of Math. {\bf 83}}
(1966), 129--209.

\bibitem[Ho2]{Ho2}{\sc ---},
Fourier integral operators I, {\it Acta Math,  }

\bibitem[Ho3]{Ho3}{\sc ---},
The Analysis of Linear Partial Differential Operators III,
Springer, Berlin, 1985.

\bibitem[Ka]{Ka}{\sc T.Kato}, 
\textit{Perturbation Theory for Linear Operators},
2nd ed., Springer-Verlag, Berlin, 1980.

\bibitem[KL]{KL}{\sc P. Kirk and M. Lesch}, 
The $\eta-$invariant, Maslov index, and
spectral flow for Dirac-type operators on manifolds with boundary,
(preprint).

\bibitem[KS]{KS}{\sc N. J. Kalton and R. C. Swanson}, 
A symplectic Banach space
with no Lagrangian subspaces, {\it Transaction of the Amer. Math. 
Soc.,{\bf 273}, No 1}
(1982), 385-392.


\bibitem[Le]{Le}{\sc J. Leray},
``Analyse Lagrangi\'enne et
m\'ecanique quantique: Une structure math\'ematique
apparent\'ee aux d\'eveloppements asymptotiques et \`a
l'indice de Maslov", 
S\'erie Math. Pure et Appl., I.R.M.P.,
Strasbourg, 1978 (English translation 1981, MIT Press).

\bibitem[Ni]{Ni}{\sc L. Nicolaescu},
The Maslov index, the spectral flow, and decomposition of
manifolds, {\it Duke Math. J. {\bf 80}} (1995), 485--533.

\bibitem[OF]{OF}{\sc N. Otsuki and K. Furutani}, Spectral flow
and Maslow index arising from Lagrangian intersections, {\it
Tokyo J. Math. {\bf 14}} (1991), 135--150.



\bibitem[Ph]{Ph}{\sc J. Phillips}, Self--adjoint Fredholm
operators and spectral flow, {\it Canad. Math. Bull. {\bf 39}}
(1996), 460--467.

\bibitem[Ra]{Ra}{\sc J. V. Ralston}, Deficiency indices of symmetric operators
with elliptic boundary conditions, {\it Comm. Pure Appl. Math.
{\bf 23}} (1970), 221--232.

\bibitem[RS]{RS}{\sc J. Robbin and D. Salamon},
The Maslov index for paths,
{\it Topology {\bf 32}} (1993),  827-844.

\bibitem[So]{So}{\sc J.M. Souriau}, Construction explicite de l'indice
de Maslov, {\it in}: ``Group Theoretical Methods in Physics'', Springer
Lecture Notes in Physics vol. 50, Berlin, 1975, pp. 117--148.

\bibitem[Sw]{Sw}
{\sc R. C. Swanson},
Fredholm intersection
theory and  elliptic boundary deformation problems, I,
{\it J. of Diff. Equations  {\bf 28}} (1978),  189--201.

\bibitem[Ta]{Ta}{\sc C. H. Taubes}, Casson's invariant and gauge
  theory, {\it J. Differential Geometry {\bf 31}}(1990), 547-599.

\bibitem[Yo]{Yo}{\sc T. Yoshida}, Floer homology and splittings
of manifolds, {\it Ann. of Math. {\bf 134}} (1991),  277--323.
\end{thebibliography}
\end{document}